\documentclass[3p,authoryear]{elsarticle2}

\usepackage[utf8]{inputenc} 
\usepackage[T1]{fontenc}
\usepackage{dsfont}
\usepackage{amsthm}
\usepackage{amsmath,amsfonts,amssymb}
\usepackage[mathscr]{eucal}
\usepackage{natbib}
\setcitestyle{notesep={; },round,aysep={},yysep={;}}
\usepackage{color}
\usepackage{xcolor}
\usepackage{alphalph,etoolbox}
\patchcmd{\subequations}{\alph{equation}}{\alphalph{\value{equation}}}{}{}
\usepackage{tikz}
\usepackage{tikz-3dplot}
\usepackage[ruled, lined, linesnumbered]{algorithm2e}
\usepackage{subcaption}
\usepackage{enumitem}
\usepackage{float}
\usepackage{booktabs,adjustbox}
\usepackage{supertabular}
\usepackage{multicol}
\usepackage{array,multirow,graphicx}

\newcount\n
\def\tablebody{}
\makeatletter
\loop\ifnum\n<100
\advance\n by1
\protected@edef\tablebody{\tablebody
	\textbf{\number\n.}& shortText
	\tabularnewline
}
\repeat

\makeatletter
\let\mcnewpage=\newpage
\newcommand{\TrickSupertabularIntoMulticols}{%
	\renewcommand\newpage{%
		\if@firstcolumn
		\hrule width\linewidth height0pt
		\columnbreak
		\else
		\mcnewpage
		\fi
	}%
}
\makeatother

\usepackage{hyperref}
\hypersetup{pdfborder={0 0 0}, colorlinks=true, citecolor=red}
\hypersetup{pdfauthor=whatever}

\theoremstyle{definition}

\newcommand{\red}[1]{{\color{red}#1}}

\newcommand{\blue}[1]{{\color{blue}#1}}
\newcommand{\gray}[1]{{\color{gray}#1}}

\newcommand{\black}[1]{{\color{black}#1}}

\newcommand{\mygreen}[1]{{\color{mygreen}#1}}
\definecolor{mygreen}{rgb}{0,0.6,0}
\definecolor{lightTeal}{RGB}{153,200,206}
\definecolor{green}{RGB}{22,125,134}

\renewcommand{\emptyset}{\varnothing}
\newcommand{\pluseq}{\mathrel{+}=}
\newcommand{\minuseq}{\mathrel{-}=}

\newcommand{\ULD}{ULD}
\newcommand{\ULDgroup}{ULD group}
\newcommand{\centerOfGravity}{CoG}
\newcommand{\problem}{MBSBPP}
\newcommand{\RGS}{RGS}
\newcommand{\GRASP}{GRASP}

\newcommand{\sorting}[1]{\text{\texttt{#1}}}
\newcommand{\variable}[1]{\text{\texttt{#1}}}

\newcommand{\itemSet}{I}
\newcommand{\itemIndex}{i}
\newcommand{\itemWeight}{w}
\newcommand{\itemRotatability}{r}
\newcommand{\itemTiltability}{t}
\newcommand{\itemStackability}{\vartheta}
\newcommand{\itemSize}{s}
\newcommand{\loadedItems}{L}
\newcommand{\numLoadedItems}{|\loadedItems|}
\newcommand{\loadedItemIndex}{\ell}
\newcommand{\loadedItemIndexTwo}{\jmath}
\newcommand{\reducedLoadedItems}{R}
\newcommand{\itemIndexTwo}{k}
\newcommand{\loadingPoint}{c}

\newcommand{\dimension}{d}
\newcommand{\dimensionTwo}{\theta}
\newcommand{\dimensionThree}{\eta}
\newcommand{\dimensions}{D}

\newcommand{\ULDSet}{C}
\newcommand{\ULDIndex}{c}
\newcommand{\ULDCapacity}{W}
\newcommand{\ULDVolume}{\mathcal{V}}
\newcommand{\edgeWidth}{e}
\newcommand{\verticalEdgeOffset}{\delta}
\newcommand{\paddingHeight}{\hslash}
\newcommand{\minimumItemOverlap}{o}

\newcommand{\vertexSet}{V}

\newcommand{\facetSet}{F}
\newcommand{\facetIndex}{f}

\newcommand{\ULDgroups}{U}
\newcommand{\ULDgroupsIndex}{u}
\newcommand{\smallesNumberFittingULDgroups}{m}

\newcommand{\size}{s}

\newcommand{\iterator}{j}
\newcommand{\numRGSiterations}{M}
\newcommand{\sortingCriterion}{\mathcal{C}}

\newcommand{\substructureAllowed}{\varsigma}
\newcommand{\useSubstructure}{b}
\newcommand{\solution}{S}

\newcommand{\extremePoint}{e}
\newcommand{\extremePoints}{E}
\newcommand{\orientations}{O}
\newcommand{\allOrientations}{\Theta}
\newcommand{\orientation}{\sigma}
\newcommand{\position}{p}
\newcommand{\blockingItemSet}{\mathcal{B}}
\newcommand{\numBlockingItems}{|\blockingItemSet|}
\newcommand{\blockingItemIndex}{b}

\newcommand{\similarItemHeight}{h}
\newcommand{\similarItemsSet}{S_{\similarItemHeight, \itemStackability}}
\newcommand{\listItemsOrientations}{\mathcal{L}}

\newcommand{\surfaceAreaHeight}{\mathcal{H}}

\newcommand{\itemGroup}{M}
\newcommand{\newItemGroup}{\bar{\itemGroup}}
\newcommand{\allItemGroups}{\mathcal{M}}
\newcommand{\allNewItemGroups}{\bar{\allItemGroups}}
\newcommand{\numItemGroups}{m}

\newcommand{\degreeOfRandomization}{\rho}
\newcommand{\selectedDegreeOfRandomization}{\delta}
\newcommand{\randomNumber}{y}

\newcommand{\boundingBox}{B}
\newcommand{\normalVector}{n}
\newcommand{\planeEquationOffset}{a}

\newcommand{\gridCellSize}{\bar{\size}}
\newcommand{\gridCellLimit}{n}
\newcommand{\gridCellIndex}{\ell}
\newcommand{\gridCellIndexMin}{\ell^{min}}
\newcommand{\gridCellIndexMax}{\ell^{max}}

\newcommand{\shiftValue}{\nu}

\newcommand{\volumeUtilizationImportance}{\mathcal{I}_v}
\newcommand{\weightBalanceImportance}{\mathcal{I}_w}
\newcommand{\volumeUtilizationScore}{\mathcal{S}_v}
\newcommand{\weightBalanceScore}{\mathcal{S}_w}
\newcommand{\solutionScore}{\mathcal{S}}
\newcommand{\maximumCenterOfGravityDeviation}{\mathit{cog}^{\max}}
\newcommand{\centerOfGravityDeviation}{\mathit{cog}^\text{dev}}
\newcommand{\remainingItems}{R}
\newcommand{\neighboringItems}{Q}

\newcommand{\utilization}{u}
\newcommand{\numItemTypes}{$\mathcal{J}$}
\newcommand{\numItems}{|\itemSet|}
\newcommand{\avgNumItems}{$\overline{|\itemSet|}$}
\newcommand{\avgNumLoadedItems}{$\overline{\numLoadedItems}$}
\newcommand{\avgUtilization}{$\bar{\utilization}$}
\newcommand{\minUtilization}{$\utilization^\text{min}$}
\newcommand{\maxUtilization}{$\utilization^\text{max}$}
\newcommand{\solutionTime}{t}
\newcommand{\avgTime}{$\bar{\solutionTime}$}

\newcommand{\maxTime}{$\solutionTime^\text{max}$}
\newcommand{\numULDs}{$|\ULDSet|$}
\newcommand{\avgNumULDs}{$\overline{|\ULDSet|}$}
\newcommand{\medianUtilization}{$\utilization^\text{med}$}
\newcommand{\numULDsCogViolated}{$\text{G}^\text{vio}$}
\newcommand{\percentageCogViolated}{\numULDsCogViolated (\%)}
\newcommand{\numSubstructureUsed}{\#sub.}
\newcommand{\total}{Total}

\newcommand{\numInstance}{No.}
\newcommand{\overallUtilization}{$\utilization^\text{total}$}
\newcommand{\containerID}{ID}
\newcommand{\CoGDevX}{$\centerOfGravityDeviation_1$}
\newcommand{\CoGDevY}{$\centerOfGravityDeviation_2$}

\newcommand{\stackCumVol}{\sorting{s--cv}}
\newcommand{\stackHighVol}{\sorting{s--hv}}
\newcommand{\cumVol}{\sorting{cv}}
\newcommand{\highVol}{\sorting{hv}}
\newcommand{\randomSorting}{\sorting{r}}

\newcommand{\containerLoading}{CL}
\newcommand{\airCargo}{AC}
\newcommand{\knapsackProblem}{KP}
\newcommand{\binPacking}{BP}
\newcommand{\palletLoading}{PL}

\newcommand{\timeQuotientNoGrid}{$\frac{\solutionTime_{NG}}{\solutionTime_{D}}$}
\newcommand{\utilizationQuotientNoGrid}{$\frac{\utilization_{NG}}{\utilization_{D}}$}
\newcommand{\timeQuotientNoBlocking}{$\frac{\solutionTime_{NB}}{\solutionTime_{D}}$}
\newcommand{\utilizationQuotientNoBlocking}{$\frac{\utilization_{NB}}{\utilization_{D}}$}
\newcommand{\timeQuotientOneExtremePoint}{$\frac{\solutionTime_{CR}}{\solutionTime_{D}}$}
\newcommand{\utilizationQuotientOneExtremePoint}{$\frac{\utilization_{CR}}{\utilization_{D}}$}
\newcommand{\timeQuotientNoExtremePointShift}{$\frac{\solutionTime_{NM}}{\solutionTime_{D}}$}
\newcommand{\utilizationQuotientNoExtremePointShift}{$\frac{\utilization_{NM}}{\utilization_{D}}$}

\newcommand{\noGrid}{\texttt{No grid}}
\newcommand{\noBlocking}{\texttt{No blocking}}
\newcommand{\oneExtremePoint}{\texttt{Crainic et al}}
\newcommand{\noExtremePointShift}{\texttt{No moving}}

\newcommand{\ie}{\emph{i.e.}}

\newcommand{\smallFilledSquare}{\mathbin{\rule{0.12cm}{0.12cm}}}

\newcommand{\mycuboid}[8]{ 
	\filldraw[cube,fill=#7,opacity=#8] (#1,#2,#6) -- (#1,#5,#6) -- (#4,#5,#6) -- (#4,#2,#6) -- cycle;
	\filldraw[cube,fill=#7,opacity=#8] (#4,#2,#3) -- (#4,#5,#3) -- (#4,#5,#6) -- (#4,#2,#6) -- cycle;
	\filldraw[cube,fill=#7,opacity=#8] (#4,#5,#6) -- (#4,#5,#3) -- (#1,#5,#3) -- (#1,#5,#6) -- cycle;
}
\tdplotsetmaincoords{70}{100}

\begin{document}

\begin{frontmatter}
\title{A Fast Optimization Approach For A Complex Real-Life 3D Multiple Bin Size Bin Packing Problem}

\author[schenker]{$\text{Katrin He\ss{}ler}$}
\ead{katrin.hessler@dbschenker.com}

\author[schenker]{Timo Hintsch}
\ead{timo.hintsch@dbschenker.com}

\author[schenker]{Lukas Wienkamp}
\ead{lukas.wienkamp@dbschenker.com}

\address[schenker]{Global Data \& AI, F.LCA - Global Technology and Data, Schenker AG, Kruppstra\ss{}e 4, 45128 Essen, Germany.\label{schenker}}

\journal{European Journal of Operational Research}

\begin{abstract}
We investigate a real-life air cargo loading problem which is a variant of the three-dimensional Variable Size Bin Packing Problem with special bin forms of cuboid and non-cuboid unit load devices (\ULD s).
Packing is constrained by additional practical restrictions, such as load stability, (non-)stackable items, and weight distribution constraints.
To solve the problem, we present an insertion heuristic embedded into a Randomized Greedy Search.
The solution space is limited by only considering certain candidate points (so-called extreme points), which are promising positions to load an item.
We extend the concept of extreme points proposed in the literature and allow moving extreme points for non-cuboid \ULD s.
A special sorting of the items is suggested, which combines a layered structure and free packing.
Moreover, we propose dividing the space of each \ULD\ into smaller cells to accelerate the collision, non-floating, and stackability check while loading items.
In a computational study, we analyze individual algorithm components and show the effectiveness of our method on adapted real-life instances from the literature.
\end{abstract}

\end{frontmatter}

\newpage

\section{Introduction}
\label{sec:intro}

In the air freight industry, there are several steps from the arrival of individual shipments at the terminal to the fully loaded and ready-to-take-off aircraft. 
One major step of the operational planning is to load cargo onto special air cargo pallets or containers, so-called \emph{unit load devices} (\ULD s), which are then loaded into the aircraft.
Air cargo pallets must be covered with a net and straps to secure their cargo.
For containers no nets or straps are required.
Both \ULD\ types, pallets and containers, can be modeled as bins with a shape that varies with respect to the designated position in the aircraft.
When loading cargo into \ULD s, the aim is to achieve a high utilization of the \ULD s and at the same time ensure load stability.
Assuming that the cargo consists of a set of cuboid items, the problem can be modeled as a three-dimensional \emph{Multiple Bin Size Bin Packing Problem} (\problem), according to the typology of \cite{WaescherEtAl2007}. 

The three-dimensional \problem\ is a generalization of the Bin Packing Problem \citep{GareyJohnson1979} and, therefore, NP-hard.
The aim of this work is to develop a solution algorithm that delivers high-quality solutions in short computational times, suitable to be used in operational planning.
Due to the problem complexity, we present a greedy-based insertion heuristic and refrain from an exact approach.
The problem we consider is a slight variation of the problem presented in \cite{PaquayEtAl2016}, in which the authors propose a mixed-integer program.
Other algorithmic approaches are a Relax-and-Fix algorithm and derived matheuristics \citep{PaquayEtAl2018b} and a two-phase insertion heuristic \citep{PaquayEtAl2018}.

The contributions of this paper are the following:
\begin{itemize}
\item We present a new variant of the three-dimensional \problem\
that takes into account additional constraints and problem characteristics that occur in real-world air transportation situations:
the stability of the load, avoiding unused space between items, (non-)stackable items, weight distribution, specially shaped, non-cuboid \ULD s, the usage of padding material, reserved edge space on pallets for cargo net locks on which loading is prohibited, and the usage of a substructure to allow for loading large items at the bottom
exceeding the \ULD\ edge.
The problem extends the one presented in \cite{PaquayEtAl2016}, which does not consider the last three points.
\item We suggest a new greedy-based insertion algorithm that sequentially loads the available \ULD s.
Each \ULD\ is loaded by a Randomized Greedy Search that repeatedly calls an insertion heuristic based on so-called extreme points \citep{CrainicEtAl2008}, which are promising candidate positions to load an item.
The framework is similar to the one proposed by \cite{GajdaEtAl2022}.
However, our approach differs from theirs and other approaches as it ensures very fast iterations thanks to various acceleration techniques:
A three-dimensional grid, dividing the space of each \ULD\ into smaller cells, is used to accelerate the collision, non-floating, and stackability check while loading items.
Moreover, instead of evaluating all generated extreme points with a merit function, a first-fit approach is applied.
The resulting very low runtime of the insertion heuristic allows a large number of iterations to be carried out.
Note that \citep{CrainicEtAl2008} already tested first-fit heuristics but did not perform several iterations with randomly sorted items.
A high solution quality is achieved by using an extended set of extreme points and a special sorting of the items which, together with the first-fit approach, combines a layered structure and free packing.
Our new approach outperforms the two-phase insertion heuristic of \cite{PaquayEtAl2018} with respect to utilization while offering similar solution quality regarding the weight balance.
Moreover, our algorithm shows good results for classical three-dimensional Bin Packing Problem instances \citep{BischoffRatcliff1995}.
\item New insights on extreme points, that can also be applied to other three-dimensional packing problems, are presented:
We allow the moving of extreme points to handle specially shaped, non-cuboid \ULD s.
This ensures a high loading density underneath the beveled \ULD\ sides at the top.
Moreover, we extend the concept of extreme points proposed by \cite{CrainicEtAl2008} by adapting the projection routine that generates new extreme points.
We show that the generation is quadratic in the number of loaded items, which is worse than the linear complexity of generating the extreme point set of \cite{CrainicEtAl2008}.
However, the algorithm is still very fast, as the worst-case complexity is almost never reached.
\end{itemize}

The remainder of the paper is organized as follows:
In the next section, we formally define the three-dimensional \problem.
In Section~\ref{sec:literature_review}, we give an overview of the existing literature on  three-dimensional \problem s.
The insertion heuristic to load a single \ULD\ is introduced in Section~\ref{sec:insertion_heuristic}.
In Section~\ref{sec:loading_ULD}, the heuristic is embedded into a Randomized Greedy Search.
To load  multiple \ULD s, the framework is extended further in Section~\ref{sec:heuristic}.
An analysis of individual algorithm components and computational results on various benchmark sets are then reported in Section~\ref{sec:comp_results}.
Section~\ref{sec:conclusion} draws final conclusions and discusses future research directions.

\section{Problem definition}
\label{sec:problem_definition}

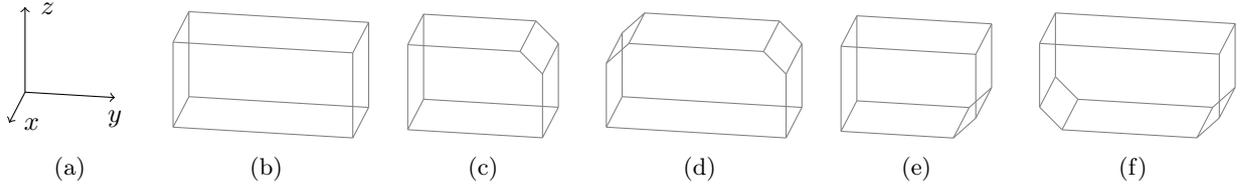
\begin{figure}
	\newcommand\widthFactor{0.15}
	\begin{subfigure}[t]{\widthFactor\textwidth}
		\centering
		\begin{tikzpicture}
			[tdplot_main_coords,
			cube/.style={black},
			scale=0.6]
			\draw[->] (0,0,0) -- (2,0,0);
			\draw[->] (0,0,0) -- (0,2,0);
			\draw[->] (0,0,0) -- (0,0,2);
			\node[] at (2,0.5,0) {$x$};
			\node[] at (0,2,-0.5) {$y$};
			\node[] at (0,0.5,2) {$z$};
		\end{tikzpicture}\\
		\caption{}\label{fig:coordinate_system}
	\end{subfigure}
	\begin{subfigure}[t]{\widthFactor\textwidth}
		\centering
		\begin{tikzpicture}
			[tdplot_main_coords,
			cube/.style={black},
			scale=0.6]
			\draw[gray] (0,0,0) -- (2,0,0);
			\draw[gray] (0,0,0) -- (0,4,0);
			\draw[gray] (0,0,0) -- (0,0,2);
			\draw[gray] (2,0,0) -- (2,4,0);
			\draw[gray] (2,0,0) -- (2,0,2);
			\draw[gray] (0,0,2) -- (2,0,2);
			\draw[gray] (0,0,2) -- (0,4,2);
			\draw[gray] (0,4,0) -- (2,4,0);
			\draw[gray] (0,4,0) -- (0,4,2);
			\draw[gray] (0,4,2) -- (2,4,2);
			\draw[gray] (2,0,2) -- (2,4,2);
			\draw[gray] (2,4,0) -- (2,4,2);
		\end{tikzpicture}\\
		\caption{}\label{fig:ULD_shapes1}
	\end{subfigure}
	\hfill
	\begin{subfigure}[t]{\widthFactor\textwidth}
		\centering
		\begin{tikzpicture}
			[tdplot_main_coords,
			cube/.style={black},
			scale=0.6]
			\draw[gray] (0,0,0) -- (2,0,0);
			\draw[gray] (0,0,0) -- (0,3,0);
			\draw[gray] (0,0,0) -- (0,0,2);
			\draw[gray] (2,0,0) -- (2,3,0);
			\draw[gray] (2,0,0) -- (2,0,2);
			\draw[gray] (0,0,2) -- (2,0,2);
			\draw[gray] (0,3,0) -- (2,3,0);
			\draw[gray] (0,3,0) -- (0,3,1.5);
			\draw[gray] (0,3,1.5) -- (2,3,1.5);
			\draw[gray] (2,3,0) -- (2,3,1.5);
			\draw[gray] (2,3,1.5) -- (2,2.5,2);
			\draw[gray] (2,2.5,2) -- (0,2.5,2);
			\draw[gray] (0,3,1.5) -- (0,2.5,2);
			\draw[gray] (0,2.5,2) -- (0,0,2);
			\draw[gray] (2,2.5,2) -- (2,0,2);
		\end{tikzpicture}\\
		\caption{}\label{fig:ULD_shapes:one_upper_cut}
	\end{subfigure}
	\hfill
	\begin{subfigure}[t]{\widthFactor\textwidth}
		\centering
		\begin{tikzpicture}
			[tdplot_main_coords,
			cube/.style={black},
			scale=0.6]
			\draw[gray] (0,0,0) -- (2,0,0);
			\draw[gray] (0,0,0) -- (0,4,0);
			\draw[gray] (0,0,0) -- (0,0,1.5);
			\draw[gray] (2,0,0) -- (2,4,0);
			\draw[gray] (2,0,0) -- (2,0,1.5);
			\draw[gray] (0,0,1.5) -- (2,0,1.5);
			\draw[gray] (0,4,0) -- (2,4,0);
			\draw[gray] (0,4,0) -- (0,4,1.5);
			\draw[gray] (2,4,0) -- (2,4,1.5);
			\draw[gray] (2,0,1.5) -- (2,0.5,2);
			\draw[gray] (2,0.5,2) -- (2,3.5,2);
			\draw[gray] (2,3.5,2) -- (2,4,1.5);
			\draw[gray] (0,0,1.5) -- (0,0.5,2);
			\draw[gray] (0,0.5,2) -- (0,3.5,2);
			\draw[gray] (0,3.5,2) -- (0,4,1.5);
			\draw[gray] (2,0.5,2) -- (0,0.5,2);
			\draw[gray] (2,3.5,2) -- (0,3.5,2);
			\draw[gray] (2,4,1.5) -- (0,4,1.5);
		\end{tikzpicture}\\
		\caption{}\label{fig:ULD_shapes:two_upper_cuts}
	\end{subfigure}
	\hfill
	\begin{subfigure}[t]{\widthFactor\textwidth}
		\centering
		\begin{tikzpicture}
			[tdplot_main_coords,
			cube/.style={black},
			scale=0.6]
			\draw[gray] (0,0,0) -- (2,0,0);
			\draw[gray] (0,0,0) -- (0,2.5,0);
			\draw[gray] (0,0,0) -- (0,0,2);
			\draw[gray] (2,0,0) -- (2,2.5,0);
			\draw[gray] (2,0,0) -- (2,0,2);
			\draw[gray] (0,0,2) -- (2,0,2);
			\draw[gray] (0,3,0.5) -- (0,3,2);
			\draw[gray] (0,3,2) -- (2,3,2);
			\draw[gray] (2,3,0.5) -- (2,3,2);
			\draw[gray] (0,0,2) -- (0,3,2);
			\draw[gray] (2,0,2) -- (2,3,2);
			\draw[gray] (0,2.5,0) -- (2,2.5,0);
			\draw[gray] (2,2.5,0) -- (2,3,0.5);
			\draw[gray] (2,3,0.5) -- (0,3,0.5);
			\draw[gray] (0,3,0.5) -- (0,2.5,0);
		\end{tikzpicture}\\
		\caption{}
	\end{subfigure}
	\hfill
	\begin{subfigure}[t]{\widthFactor\textwidth}
		\centering
		\begin{tikzpicture}
			[tdplot_main_coords,
			cube/.style={black},
			scale=0.6]
			\draw[gray] (0,0.5,0) -- (2,0.5,0);
			\draw[gray] (0,0.5,0) -- (0,3.5,0);
			\draw[gray] (0,0,0.5) -- (0,0,2);
			\draw[gray] (2,0.5,0) -- (2,3.5,0);
			\draw[gray] (2,0,0.5) -- (2,0,2);
			\draw[gray] (0,0,2) -- (2,0,2);
			\draw[gray] (0,0,2) -- (0,4,2);
			\draw[gray] (0,3.5,0) -- (2,3.5,0);
			\draw[gray] (0,4,0.5) -- (0,4,2);
			\draw[gray] (0,4,2) -- (2,4,2);
			\draw[gray] (2,0,2) -- (2,4,2);
			\draw[gray] (2,4,0.5) -- (2,4,2);
			\draw[gray] (2,3.5,0) -- (2,4,0.5);
			\draw[gray] (2,4,0.5) -- (0,4,0.5);
			\draw[gray] (0,4,0.5) -- (0,3.5,0);
			\draw[gray] (2,0.5,0) -- (2,0,0.5);
			\draw[gray] (2,0,0.5) -- (0,0,0.5);
			\draw[gray] (0,0,0.5) -- (0,0.5,0);
		\end{tikzpicture}\\
		\caption{}
	\end{subfigure}
	\caption{The underlying coordinate system and different \ULD\ shapes relevant in practice.}\label{fig:ULD_shapes}
\end{figure}

We first provide some practical insights before formally defining the problem.
The underlying coordinate system is depicted in Figure~\ref{fig:coordinate_system}.
Note that the length, width, and height correspond to $x$, $y$, and $z$, respectively.

A set of cuboid items must be loaded into a set of \ULD s.
Some items can be rotated orthogonally and/or tilted while others cannot.
There are six possible orientations (rotated: yes/no, tilted: no/across the $x$-axis/across the $y$-axis) for rotatable and tiltable items \cite[Figure 2]{PaquayEtAl2018}.
Tiltable items are always rotatable.
Furthermore, some fragile items, for example those whose outside is a cardboard, are not stackable, \ie, nothing is allowed to be placed on top.

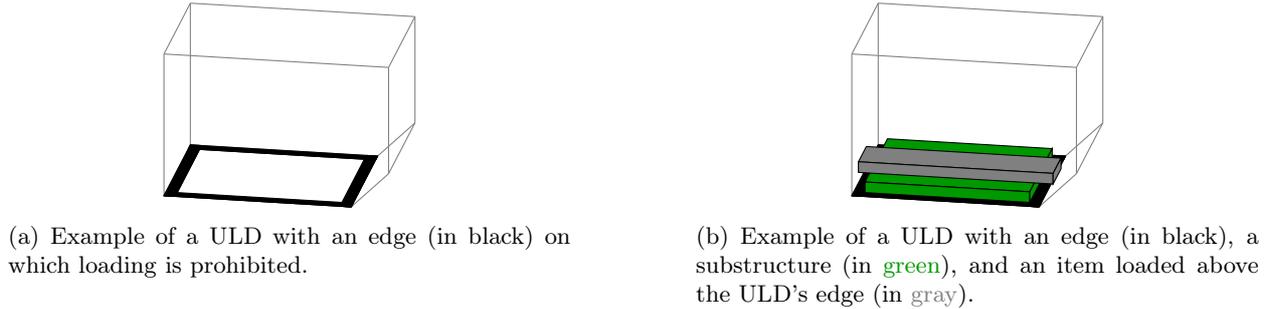
\begin{figure}
	\begin{subfigure}[t]{0.45\textwidth}
		\centering
		\begin{tikzpicture}
			[tdplot_main_coords,
			cube/.style={black},
			scale=1.0]
			\draw[gray] (0,0,0) -- (2,0,0);
			\draw[gray] (0,0,0) -- (0,2.5,0);
			\draw[gray] (0,0,0) -- (0,0,2);
			\draw[gray] (2,0,0) -- (2,2.5,0);
			\draw[gray] (2,0,0) -- (2,0,2);
			\draw[gray] (0,0,2) -- (2,0,2);
			\draw[gray] (0,3,0.5) -- (0,3,2);
			\draw[gray] (0,3,2) -- (2,3,2);
			\draw[gray] (2,3,0.5) -- (2,3,2);
			\draw[gray] (0,0,2) -- (0,3,2);
			\draw[gray] (2,0,2) -- (2,3,2);
			\draw[gray] (0,2.5,0) -- (2,2.5,0);
			\draw[gray] (2,2.5,0) -- (2,3,0.5);
			\draw[gray] (2,3,0.5) -- (0,3,0.5);
			\draw[gray] (0,3,0.5) -- (0,2.5,0);
			\mycuboid{0}{0}{0}{2}{.15}{0}{black}{1};
			\mycuboid{0}{0}{0}{.15}{2.5}{0}{black}{1};
			\mycuboid{1.85}{0}{0}{2}{2.5}{0}{black}{1};
			\mycuboid{0}{2.35}{0}{2}{2.5}{0}{black}{1};
		\end{tikzpicture}\\
		\caption{Example of a \ULD\ with an edge (in black) on which loading is prohibited.}
	\end{subfigure}
	\hfill
	\begin{subfigure}[t]{0.45\textwidth}
		\centering
		\begin{tikzpicture}
			[tdplot_main_coords,
			cube/.style={black},
			scale=1.0]
			\draw[gray] (0,0,0) -- (2,0,0);
			\draw[gray] (0,0,0) -- (0,2.5,0);
			\draw[gray] (0,0,0) -- (0,0,2);
			\draw[gray] (2,0,0) -- (2,2.5,0);
			\draw[gray] (2,0,0) -- (2,0,2);
			\draw[gray] (0,0,2) -- (2,0,2);
			\draw[gray] (0,3,0.5) -- (0,3,2);
			\draw[gray] (0,3,2) -- (2,3,2);
			\draw[gray] (2,3,0.5) -- (2,3,2);
			\draw[gray] (0,0,2) -- (0,3,2);
			\draw[gray] (2,0,2) -- (2,3,2);
			\draw[gray] (0,2.5,0) -- (2,2.5,0);
			\draw[gray] (2,2.5,0) -- (2,3,0.5);
			\draw[gray] (2,3,0.5) -- (0,3,0.5);
			\draw[gray] (0,3,0.5) -- (0,2.5,0);
			\mycuboid{0}{0}{0}{2}{.15}{0}{black}{1};
			\mycuboid{0}{0}{0}{.15}{2.5}{0}{black}{1};
			\mycuboid{1.85}{0}{0}{2}{2.5}{0}{black}{1};
			\mycuboid{0}{2.35}{0}{2}{2.5}{0}{black}{1};
			\mycuboid{.15}{.15}{0}{1.85}{2.35}{.15}{mygreen}{1};
			\mycuboid{0.9}{0}{.15}{1.5}{2.6}{.3}{gray}{1};
		\end{tikzpicture}\\
		\caption{Example of a \ULD\ with an edge (in black), a substructure (in \mygreen{green}), and an item loaded above the \ULD's edge (in \gray{gray}).}
	\end{subfigure}
	\caption{Two examples of \ULD s with an edge and with/without a substructure.}\label{fig:edge_substructure}
\end{figure}

In case the cargo is loaded onto a pallet (and not into a closed container), the available space for loading is not only defined by the outer shape of the  \ULD\ but further restricted.
The reason for this is that the cargo must be covered with a net which is attached to the edge of the pallet via net locks.
It is not possible to load cargo directly onto these locks.
However, if an item is sufficiently far above the \ULD\ edge, for example because it is stacked on top of other items, it may overlap with the edge as access to the locks is guaranteed.
To allow for the loading of large items at the bottom that can only be loaded if placed above the \ULD's edge, a substructure can be used to artificially raise the \ULD's base area.
Figure~\ref{fig:edge_substructure} illustrates the \ULD's edge and the use of a substructure.

Load stability is crucial to avoid damages during transportation.
To prevent items from sliding, it is preferable to avoid holes between items in horizontal direction.
To fill small vertical gaps and allow for sufficient support of items, padding material can be placed on top of (stackable) items and the \ULD\ floor.
We refer to the part of the item that is supported by padding material as \emph{indirectly supported}.
The part of the item that rests on another stackable item or the floor is \emph{directly supported}.

\begin{figure}
	\begin{subfigure}[t]{0.45\textwidth}
		\centering
		\begin{tikzpicture}
			[tdplot_main_coords,
			cube/.style={black},scale=0.8]
			\draw[gray] (0,0,0) -- (5,0,0);
			\draw[gray] (0,0,0) -- (0,5,0);
			\draw[gray] (0,0,0) -- (0,0,2.5);
			\draw[gray] (5,0,0) -- (5,5,0);
			\draw[gray] (5,0,0) -- (5,0,2.5);
			\draw[gray] (0,0,2.5) -- (5,0,2.5);
			\draw[gray] (0,0,2.5) -- (0,5,2.5);
			\draw[gray] (0,5,0) -- (5,5,0);
			\draw[gray] (0,5,0) -- (0,5,2.5);
			\draw[mygreen, line width = 1.5mm] (3.5,5,0) -- (1.5,5,0);
			\draw[mygreen, line width = 1.5mm] (5,3.5,0) -- (5,1.5,0);
			\mycuboid{1.5}{1.5}{0}{3.5}{3.5}{0}{mygreen}{1};
			\draw[gray] (0,5,2.5) -- (5,5,2.5);
			\draw[gray] (5,0,2.5) -- (5,5,2.5);
			\draw[gray] (5,5,0) -- (5,5,2.5);	
		\end{tikzpicture}\\
		\caption{Feasible region for the \centerOfGravity.}
	\end{subfigure}
	\hfill
	\begin{subfigure}[t]{0.45\textwidth}
		\centering
		\begin{tikzpicture}
			[tdplot_main_coords,
			cube/.style={black},scale=0.8]
			\draw[gray] (0,0,0) -- (5,0,0);
			\draw[gray] (0,0,0) -- (0,5,0);
			\draw[gray] (0,0,0) -- (0,0,2.5);
			\draw[gray] (5,0,0) -- (5,5,0);
			\draw[gray] (5,0,0) -- (5,0,2.5);
			\draw[gray] (0,0,2.5) -- (5,0,2.5);
			\draw[gray] (0,0,2.5) -- (0,5,2.5);
			\draw[gray] (0,5,0) -- (5,5,0);
			\draw[gray] (0,5,0) -- (0,5,2.5);
			\draw[mygreen, line width = 1.5mm] (3.5,5,0) -- (1.5,5,0);
			\draw[mygreen, line width = 1.5mm] (5,3.5,0) -- (5,1.5,0);
			\mycuboid{1.5}{1.5}{0}{3.5}{3.5}{0}{mygreen}{1};
			\mycuboid{0}{0}{0}{1}{1}{1}{white}{1};
			\draw[red,dotted,very thick] (0.5,0.5,0.5) -- (0.5,0.5,0);
			\draw[red,dotted,very thick] (0.5,0.5,0) -- (0.5,5,0);
			\draw[red,dotted,very thick] (0.5,0.5,0) -- (5,0.5,0);
			\foreach \Point in {(0.5,0.5,0), (0.5,5,0), (5,0.5,0)}{
				\node[dotted] at \Point {\red{\textbullet}};
			}
			\foreach \Point in {(0.5,0.5,0.5)}{ 
				\node[dotted] at \Point {\black{$\smallFilledSquare$}};
			}
			\draw[gray] (0,5,2.5) -- (5,5,2.5);
			\draw[gray] (5,0,2.5) -- (5,5,2.5);
			\draw[gray] (5,5,0) -- (5,5,2.5);	
		\end{tikzpicture}\\
		\caption{The \centerOfGravity\ is marked with a black square~$\smallFilledSquare$. The \centerOfGravity\ in $x$- and $y$-direction is marked with a \red{red} bullet~\red{\textbullet}, respectively. In this example, the \centerOfGravity\ is outside the feasible region in $x$- and $y$-direction.}
	\end{subfigure}
	\caption{Example of how to determine the \centerOfGravity\ in $x$- and $y$-direction.}\label{fig:weight_balance_3d}
\end{figure}
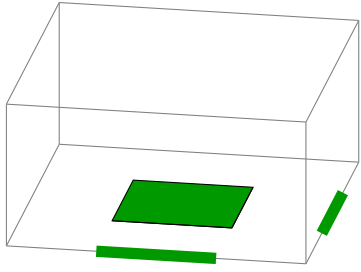
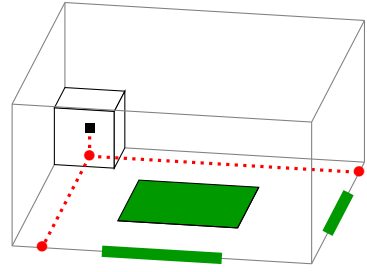

It is required that the \ULD's \emph{center of gravity} (\centerOfGravity) in $x$- and $y$-direction falls into a pre-defined area around its geometric center (see Figure~\ref{fig:weight_balance_3d}).
We assume that the weight of all items is equally distributed so that the center of gravity of each item is its geometrical center.

\bigskip

We can formally define the problem as follows.

\section*{Input}

\paragraph{Items}
A set  $\itemSet$ of cuboid items. For each item~$\itemIndex\in\itemSet$, we have
\begin{itemize}[noitemsep,topsep=0pt,label={--}]
	\item size $\itemSize^\itemIndex = (\itemSize^\itemIndex_1, \itemSize^\itemIndex_2, \itemSize^\itemIndex_3)$ with length~$\itemSize_1^\itemIndex$,  width~$\itemSize_2^\itemIndex$, height~$\itemSize_3^\itemIndex$,
	\item weight~$\itemWeight_\itemIndex$,
	\item is rotatable $\itemRotatability_\itemIndex \in \{0, 1\}$,
	\item is tiltable $\itemTiltability_\itemIndex \in \{0, 1\}$,
	\item is stackable $\itemStackability_\itemIndex \in \{0,1\}$.
\end{itemize}

\paragraph{ULDs}
A set $\ULDSet$ of \ULD s. For each \ULD~$\ULDIndex\in\ULDSet$, we have
\begin{itemize}[noitemsep,topsep=0pt,label={--}]
	\item set of vertices $\vertexSet_\ULDIndex$,
	\item set of facets $\facetSet_\ULDIndex$ where each facet $\facetIndex \in \facetSet_\ULDIndex$ consists of the vertices defining that facet ($\facetIndex \subseteq \vertexSet_\ULDIndex$),
	\item weight capacity $\ULDCapacity_\ULDIndex$,
	\item volume capacity $\ULDVolume_\ULDIndex$,
	\item edge width $\edgeWidth_\ULDIndex$,
	\item vertical edge offset $\verticalEdgeOffset_\ULDIndex$ (minimum vertical distance between an item and the \ULD\ edge to allow for edge overlap),
	\item use of substructure allowed $\substructureAllowed_\ULDIndex \in \{0, 1\}$.
\end{itemize}

\paragraph{Packing parameters}
Additional parameters defining what constitutes a feasible solution:
\begin{itemize}[noitemsep,topsep=0pt,label={--}]
	\item maximum padding height $\paddingHeight$ (maximum height of a gap that can be filled with padding material).
	\item minimum item overlap $\minimumItemOverlap$ (minimum portion of an item's base area that must be supported by other items, padding material, or the floor).
	\item maximum center of gravity deviation from a \ULD 's geometric center $\maximumCenterOfGravityDeviation$.
\end{itemize}

\section*{Constraints}

\paragraph{ULD fit}
Each loaded item must lay within the \ULD.

\paragraph{ULD edge}
Items must not be placed on the ULD edge and can only overlap with the edge at a height of at least $\verticalEdgeOffset$.

\paragraph{Collisions}
No two loaded items must collide with each other.

\paragraph{Floating}
Each loaded item must sufficiently be supported by either the ULD floor, other loaded items, or padding material.
Here, sufficiently supported means that
\begin{enumerate}[noitemsep,topsep=0pt,label=(\roman*)]
	\item the item has any level of direct support (\ie, at least one part of the item rests on another item or the floor) and
	the percentage of the item's base area that is directly or indirectly supported (\ie, rests on another item, padding material, or the floor) amounts to at least the minimum item overlap~$\minimumItemOverlap$, or
	\item the item's four bottom corner points rest on items or the \ULD\ floor.
\end{enumerate}

\paragraph{Stacking}
Neither items nor padding material must be loaded on top of a non-stackable item.

\paragraph{Weight distribution}
The \centerOfGravity\ of the \ULD\ load must not deviate more than $\maximumCenterOfGravityDeviation$ from the \ULD's geometric center.

\section*{Objective}

Respecting all packing restrictions, the task is to (1) maximize the volume of loaded items and (2) determine a volume-minimal set of \ULD s.
Both objectives are optimized in a strictly lexicographic sense.

\section*{Supported ULD shapes}

The definition of a \ULD\ via a set of vertices and facets allows for arbitrary shapes to be expressed.
However, only few of them are relevant in practice.
Hence, we make the following assumptions on the supported \ULD\ shapes, which are also exploited by the algorithm.
Each \ULD\ consists of at least six facets, with exactly two facets parallel to the $x$-$y$, $x$-$z$, and $y$-$z$ plane, respectively.
We refer to these facets as \emph{non-tilted facets}.
In addition, some \ULD s may have a special shape with up to two \emph{tilted facets} to fit into the fuselage of the aircraft.
These are facets that are not parallel to the $x$-$y$, $x$-$z$, or $y$-$z$ plane but orthogonal to the $y$-$z$ plane.
No neighboring tilted facets are allowed, \ie, corner vertices of tilted facets must differ.
Also, we assume that each \ULD\ is convex.
Figure~\ref{fig:ULD_shapes} gives an overview of outer \ULD\ shapes that are relevant in practice.

\section*{Further notation}

We introduce the following notation, that will be used throughout the remainder of the paper.
\begin{itemize}[noitemsep,topsep=0pt,label={--}]
\item We denote an item's $\itemIndex$ orientation by
\begin{align*}
\orientation_\itemIndex = (\text{tilt, rotation}) \in \allOrientations, \text{where } \allOrientations = \{\text{not tilted, tilted across $x$-axis, tilted across $y$-axis}\} \times \{0, 1\}.
\end{align*}
Here, the tilt is always applied first, \ie, refers to the item's original size/orientation, and the rotation to the already tilted item.
When we apply a certain orientation to an item $\itemIndex$, we update its size $\itemSize^\itemIndex = (\itemSize^\itemIndex_1, \itemSize^\itemIndex_2, \itemSize^\itemIndex_3)$ accordingly.
\item We identify each facet~$\facetIndex \in \facetSet_\ULDIndex$ with its plane equation~$\normalVector^\facetIndex_1 x+\normalVector^\facetIndex_2 y+\normalVector^\facetIndex_3 z=\planeEquationOffset^\facetIndex$, where the normal vector~$\normalVector^\facetIndex$ is assumed to point towards the inside of the \ULD.
\item We denote the index set of dimensions/directions as $\dimensions = \{1, 2, 3\}$ corresponding to $\{x, y, z\}$.
\item When we load an item $\itemIndex$, its \emph{position} $\extremePoint^\itemIndex$ refers to its corner point closest to the coordinate origin, \ie, the corner with the minimum $x$-, $y$- and $z$-coordinates. We refer to the corner point furthest away from the coordinate origin, \ie, the one with the maximum $x$-, $y$- and $z$-coordinates, as \emph{end position}.
\item Whenever we focus on only one item $\itemIndex \in \itemSet$ or one \ULD\ $\ULDIndex \in \ULDSet$, we drop the index $\itemIndex$/$\ULDIndex$ from all item/\ULD\ related parameters.
\end{itemize}

\section{Literature review}
\label{sec:literature_review}

In the literature, several variants of the three-dimensional \problem\ have been discussed.
In the following, we first provide an overview of constraints relevant in practice and the corresponding three-dimensional packing problems.
Afterwards, we focus on heuristic solution approaches for such rich packing problems.
With regard to packing problems with practical constraints in the air freight industry, we refer to \cite{BortfeldtWaescher2013,ZhaoEtAl2014,Brandt2017,BrandtNickel2019} for a summary of the topic and related problems.

\cite{BischoffRatcliff1995} provide a first overview of packing requirements relevant in practice.
Among others, they mention container weight capacity, weight distribution, grouping of items, load stability, and multiple item destinations (multi-drop).
The paper focuses on container loading problems, however, several of the restrictions are also relevant for pallet building problems and for loading \ULD s in the air freight industry.
Table~\ref{tab:overview_packing_constraints} provides an overview of various problem definitions.
We have limited the selection to problems considering weight distribution and/or stackability constraints.

\setlength{\tabcolsep}{5.0pt}
\begin{table}[htbp]
	\centering
	\caption{Overview of real-world three-dimensional packing problem definitions and constraints.}
	\begin{tabular}{ll|lllllllllll}
		\hline
		category & problem feature &  \rotatebox[origin=l]{90}{\cite{BischoffRatcliff1995}} & \rotatebox[origin=l]{90}{\cite{RatcliffBischoff1998}} & \rotatebox[origin=l]{90}{\cite{DaviesBischoff1999}} & \rotatebox[origin=l]{90}{\cite{TernoEtAl2000}} & \rotatebox[origin=l]{90}{\cite{ChanEtAl2006}} & \rotatebox[origin=l]{90}{\cite{BaldiEtAl2012}} &  \rotatebox[origin=l]{90}{\cite{JunqueiraEtAl2012}} & \rotatebox[origin=l]{90}{\cite{CeschiaSchaerf2013}} & \rotatebox[origin=l]{90}{\cite{PaquayEtAl2016}} & \rotatebox[origin=l]{90}{\cite{TrivellaPisinger2016}} & \rotatebox[origin=l]{90}{our approach} \\
		\hline
		type &       &  \containerLoading & \containerLoading &  \containerLoading & \palletLoading & \airCargo    & \knapsackProblem &  \containerLoading &  \containerLoading  &   \airCargo    &  \binPacking &   \airCargo    \\[1ex]
		bin/\ULD & weight capacity & x$^*$ &  &    &  x     &     x$^\ddag$    &     &       &    x   &   x    &       & x \\
		& weight distribution & x$^*$ &  & x   &  x & x & x     & x &  & x &   x    & x \\
		& multiple bin sizes &       &       &       &       &   x    &      &       &    x   &   x    &       &  x\\
		& non-cuboid \ULD s &  &  &       &       &   x    &       &       &       &     x  &       & x \\
		& \ULD\ edge &  &       &       &       &       &       &       &       &       &       & x \\[1ex]
		item & orientation constraints & x & x & x     & x &   x    &   x    &   x    &   x    &   x    &       & x \\
		& stackability & x &  &      &  &       &       &       &       &    x   &       &  x\\
		& load bearing capacity &       &   x    &    &    x   &       &     &    x   &    x   &       &       & \\[1ex]
		packing & multi-drop &  x     &       &    &       &       &       &       &   x    &       &       &  \\
		& grouping &   x$^*$    &       &   & x      &    x   &      &       &       &       &       &  \\
		& filling material &  &       &       &       &       &       &       &       &       &       &  x\\
		\hline
		\multicolumn{13}{l}{\containerLoading:~Container Loading, \palletLoading:~Pallet Loading, \airCargo:~Air Cargo, \knapsackProblem:~Knapsack Problem, \binPacking:~Bin Packing.}\\
		\multicolumn{13}{l}{$^*$:~Mentioned but not considered in solution approach. \quad $^\ddag$:~Indirectly considered by objective function.}
	\end{tabular}%
	\label{tab:overview_packing_constraints}%
\end{table}%
\setlength{\tabcolsep}{6.0pt}

Weight capacity constraints are relevant in most real-world problems.
Weight balance constraints, which ensure the center of gravity falls (as close as possible) to an ideal point or a certain area, are considered in the context of container loading problems \citep{DaviesBischoff1999}, pallet building problems \citep{TernoEtAl2000}, air cargo loading problems \citep{ChanEtAl2006,PaquayEtAl2016}, the three-dimensional knapsack problem \citep{BaldiEtAl2012}, and bin packing problem \citep{TrivellaPisinger2016}.
In practice, these restrictions are important, for example, when lifting a container with a crane or transporting pallets with a forklift.
Note that weight balance constraints should not be confused with axle weight limits of trucks \citep{PollarisEtAl2016} or an as evenly as possible distributed weight along truck axles \citep{GajdaEtAl2022}.
Multiple bin sizes are often considered in air cargo loading problems.
The same holds for special non-cuboid bin shapes that are only relevant in air cargo loading problems \citep{PaquayEtAl2016,ChanEtAl2006}.
To the best of our knowledge, we are the only ones considering the \ULD\ edge.
Some approaches also consider very specific cost functions \citep{ChanEtAl2006}.
However, we do not go into further detail because our approach simply maximizes the volume of loaded items.

Due to its content, an item may only allow a limited number of orientations.
In \cite{TrivellaPisinger2016,JunqueiraEtAl2012}, rotating items is not allowed.
Moreover, it can be specified whether an item is stackable, \ie, not fragile.
A more detailed stackability definition is proposed by \cite{RatcliffBischoff1998} that define a maximum load bearing capacity depending on the orientation of the item.
Other load bearing capacity definitions are introduced in \cite{TernoEtAl2000,JunqueiraEtAl2012,CeschiaSchaerf2013}.
Multi drop scenarios, \ie, situations in which items must be delivered to different destinations, generate further restriction on feasible packings.
Here, loading patterns that allow for access to the items for the earlier destinations without moving other items around are required.
This typically results in solutions following a vertical wall-/stack-building approach to ensure that items with the same destination can be placed on top of each other to avoid spreading them out across the bin floor \citep{BischoffRatcliff1995,CeschiaSchaerf2013}.
Ensuring that pre-defined groups of items are loaded in the same bin is considered in \cite{TernoEtAl2000,ChanEtAl2006}.
The problem definition that comes closest to ours is the one of \cite{PaquayEtAl2016}, which consider the same problem without the \ULD\ edge and filling material.

Most of the heuristic approaches rely on the seminal works of \cite{MartelloEtAl2000,CrainicEtAl2008} that introduce ways to generate promising positions to load items.
\cite{MartelloEtAl2000} propose the concept of so-called corner points and, based on that, develop an exact branch-and-bound algorithm to fill a single bin.
Later on, \cite{CrainicEtAl2008} extend the idea of corner points to the concept of extreme points allowing for the creation of additional candidate points for loading items.
While in \cite{MartelloEtAl2000} all corner points have to be recalculated when an additional item is loaded, \cite{CrainicEtAl2008} suggest a procedure that only needs to calculate additional extreme points upon loading an item.
Hence, despite the potential larger number of candidate points, the complexity of updating the extreme points is only  $\mathcal{O}(\numLoadedItems)$ whereas recalculating the corner points has $\mathcal{O}(\numLoadedItems^2)$ complexity (where $\numLoadedItems$ denotes the number of already loaded items).
\cite{CrainicEtAl2008} present a constructive heuristic based on the extreme point concept.

Several other heuristics rely on local search.
\cite{LodiEtAl2002,LodiEtAl2004} present a unified tabu search framework \emph{TSpack}, that is applicable to various multi-dimensional bin-packing problem variants.
\cite{FaroeEtAl2003} propose a guided local search that iteratively removes one bin from a feasible solution.
It uses memory to guide the search to promising regions of the solution space.
A two-level tabu search \emph{TS$^2$pack} is proposed by \cite{CrainicEtAl2009}, where the first level aims to reduce the number of bins and the second level optimizes the packing of the bins.
A local search approach for a multi-drop multi-container loading problem is introduced by \cite{CeschiaSchaerf2013}.
\cite{TrivellaPisinger2016} present a multi-level local search heuristic for a load-balanced multi-dimensional bin-packing problem.
The algorithm uses an implicit representation of multi-dimensional packings by means of interval graphs \citep{FeketeSchepers2004} and iteratively improves the load balancing using different search levels.

\cite{ParrenoEtAl2008} present a \emph{Greedy Randomized Adaptive Search Procedure} (\GRASP) that is based on a constructive block heuristic using the concept of maximal free space.
Later, the authors propose a hybrid \GRASP\ and Variable Neighborhood Descent algorithm \citep{ParrenoEtAl2010}, where the constructive phase is based on their previous work.
In \cite{AlvarezValdesEtAl2013}, the best solutions found by a GRASP algorithm are combined into a Path Relinking procedure to intensify and diversify the search.
Other heuristic approaches for three-dimensional packing problems can be found in \cite{ChanEtAl2006,EgebladPisinger2009,BaldiEtAl2012}.

The solution approach that is closest to ours is the randomized constructive heuristic of \cite{GajdaEtAl2022} that combines items, sorts the combined items, partially perturbs the sorting randomly, and finally constructs the packing.
Within the constructive packing phase, the location at which an item is placed is determined by iterating over potential points, which are a subset of extreme points.
The supposedly best location is selected using a merit function considering the proportion of the contact surface with the underlying item.
The algorithm performs well on large-scale industry instances.

Returning to the problem definition of \cite{PaquayEtAl2016}, which is closest to ours, two heuristic approaches to this problem have been proposed: a MIP-based constructive heuristic \citep{PaquayEtAl2018b} and a tailored two-phase constructive heuristic \citep{PaquayEtAl2018b}.
We will compare our algorithm with the latter in Section~\ref{sec:comp_results_paquay}.

Another common approach for three-dimensional packing problems with practical constraints is to model it as a two-level problem where the lower level creates layers and the higher level combines them.
Packing horizontal layers is often desirable, especially to assure a high load stability.
Approaches are presented in \cite{MackBortfeldt2010,ElhedhliEtAl2019,GzaraEtAl2020,CalzavaraEtAl2021}.

\section{Insertion Heuristic}\label{sec:insertion_heuristic}
In this section, we introduce the insertion heuristic to load a single \ULD.
Its sole objective is to maximize the utilized volume of the given \ULD.
The weight balance will only be considered in Section~\ref{sec:loading_ULD}, when we embed the heuristic into a \emph{Randomized Greedy Search} (\RGS).
The framework will be extended further to allow for the loading of multiple \ULD s in Section~\ref{sec:heuristic}.

Our algorithm is based on the seminal extreme point insertion heuristic introduced by \citet{CrainicEtAl2008}, who extended the corner point heuristic of \cite{MartelloEtAl2000}.
The fundamental idea is to load the items one by one according to a predefined sequence, beginning at one of the \ULD’s corners.
Newly loaded items are tried to be placed close to already loaded ones in a way that avoids fragmenting the remaining space.
To that end, whenever an item is loaded, a set of so-called extreme points is computed.
These are candidate points at which further items can be loaded.

Aside from supporting \ULD s with blocked edge space for cargo net locks, one of the main difference of our approach to the one of \citet{CrainicEtAl2008} (and \citet{PaquayEtAl2018}) is an extension of the extreme points concept by generating more points (Section~\ref{sec:generating_extreme_points}) and allowing some of them to be moved (Section~\ref{sec:choose_extreme_point}).
Additionally, we do not evaluate all generated extreme points to decide where the next item is loaded but perform a first fit approach.
As a result, our insertion heuristic is very fast so that it can be called multiple times by the \RGS.

On a high level, Algorithm~\ref{alg:insertion_heuristic} outlines the insertion heuristic we propose.
The selection of the sorting criterion and whether a substructure is used is dealt with by the \RGS.
After initializing the set of extreme points, the \ULD\ is adapted to handle the \ULD's edge and a possible substructure.
Then, the items are grouped and sorted resulting in an ordered list of items and orientations.
In this ordered list, each item can appear multiple times with different orientations.
Afterwards, the next item to be loaded and a set of orientations are selected iteratively, and the item is loaded at the first extreme point that does not violate any loading restriction.
In the following, we will elaborate in greater detail
\begin{itemize}[leftmargin=1.8cm]\setlength\itemsep{-0.5ex}
	\item[(line \ref{alg:adapt_ULD})] in Section~\ref{sec:adapt_ULDs}, how the \ULD\ is adapted to handle the \ULD's edge and a possible substructure,
	\item[(lines \ref{alg:create_sorted_groups}--\ref{alg:choose_item})] in Section~\ref{sec:choose_item}, how items are grouped and sorted, and how the next item and the set of orientations are chosen,
	\item[(line \ref{alg:choose_extreme_point})] in Section~\ref{sec:choose_extreme_point}, how the next extreme point is chosen and potentially moved,
	\item[(line \ref{alg:load_item})] in Section~\ref{sec:can_be_loaded}, how to check whether an item can be loaded at an extreme point, and
	\item[(line \ref{alg:update_extreme_points})] in Section~\ref{sec:generating_extreme_points}, how the set of extreme points is updated.
\end{itemize}

\begin{algorithm}
	\DontPrintSemicolon
	\SetKw{Continue}{continue}
	\SetKw{Break}{break}
	\SetKwInput{Input}{Input}
	\SetKwInput{Output}{Output}
	\newcommand\mycommfont[1]{\footnotesize\ttfamily\textcolor{mygreen}{#1}}
	\SetCommentSty{mycommfont}
	\SetNoFillComment
	\Input{Set of items $\itemSet$ and a \ULD~$\ULDIndex\in\ULDSet$, \linebreak
		sorting criterion $\sortingCriterion$, \linebreak
		degree of randomization~$\degreeOfRandomization$, \linebreak
		whether a substructure is used $\useSubstructure$}
	\Output{Loaded \ULD}
	Define set of extreme points $\extremePoints=\{(0,0,0)\}$\;
	Adapt \ULD\ (dependent on $\useSubstructure$)\;\label{alg:adapt_ULD}
	Determine ordered list~$\listItemsOrientations$ of items and orientations by means of sorting criterion~$\sortingCriterion$\;\label{alg:create_sorted_groups}
	\For{$(\itemIndex, \orientations)\in \listItemsOrientations$}{\label{alg:choose_item}
		\While{$\extremePoint=getNextExtremePoint(\extremePoints)$ exists  and $\itemIndex$ is not loaded\label{alg:choose_extreme_point}}{
			\For{orientation $\orientation \in \orientations$}{\label{alg:tilts_and_rotations}
				\If{$\itemIndex$ can be loaded at the (potentially moved) extreme point $\extremePoint$ with orientation $\orientation$}{\label{alg:load_item}
					Load $\itemIndex$ at the (potentially moved) extreme point $\extremePoint$ with orientation $\orientation$\;
					Remove $\itemIndex$ from $\listItemsOrientations$\;
					Update set of extreme points $\extremePoints$\;\label{alg:update_extreme_points}
					\Break\;
				}
			}
		}
	}
	\Return loaded ULD\;
	\caption{Insertion heuristic}\label{alg:insertion_heuristic}
\end{algorithm}

\subsection{\ULD\ adaption}\label{sec:adapt_ULDs}

For \ULD s with edge width $\edgeWidth>0$ and vertical edge offset $\verticalEdgeOffset>0$, the available space for loading is not only defined by the outer shape of the \ULD, but needs to be further restricted.
To ensure that the space above the \ULD\ edge up until the height~$\verticalEdgeOffset$ is blocked, we load four non-stackable dummy items of height $\verticalEdgeOffset - 1$ on the edge of the \ULD\ floor.
Here, the non-stackability guarantees that the dummy items cannot be used as support for any other actual item loaded on top.
Setting the height of the dummy item to $\verticalEdgeOffset - 1$ allows items to overlap with the edge at height $\verticalEdgeOffset$ without touching the dummy items (and thus violating the stackability constraint).

Using a substructure can also be modelled by dummy items.
Next to the four dummy items blocking the space above the edge, we additionally add one stackable dummy item of height~$\verticalEdgeOffset$ covering the whole \ULD\ floor except the edges.

\subsection{Sorting and selecting the next item}\label{sec:choose_item}

\cite{CrainicEtAl2008} propose different two-level sorting criteria that are based on the volume, height, and area of the items.
We adopt some of their ideas and extend them by introducing different kinds of item groups to ensure that similar or even identical items are loaded
in close proximity for increased packing density.
Here, a particular focus lies on ensuring that items with the same height are loaded next to each other to create packing layers that allow for maximum support
of further items loaded on top.
Additionally, we consider an item's stackability when selecting the next item to avoid loading non-stackable items close to the \ULD\ floor and thus blocking the
loadable space above. 
\ \\

In line~\ref{alg:create_sorted_groups} of Algorithm~\ref{alg:insertion_heuristic}, we create a list of items and orientations defining the order in which items will be loaded.
This is done by creating two sets of item groups, ordering them according to a provided sorting criterion, and then determining the corresponding item orientations.

\subsubsection{Item groups}\label{sec:item_groups}
We begin by clustering all items into groups of \textit{identical items}.
Two items $\itemIndex, \itemIndexTwo$ are considered identical if
$(i)$ they have the same weight $\itemIndex = \itemIndexTwo$,
$(ii)$ they have the same loading characteristics (rotatability, tiltability, stackability) $\itemRotatability_\itemIndex = \itemRotatability_\itemIndexTwo$, $\itemTiltability_\itemIndex = \itemTiltability_\itemIndexTwo$,  $\itemStackability_\itemIndex = \itemStackability_\itemIndexTwo$, and
$(iii)$ they have the same set of dimensions $\{\itemSize^\itemIndex_1, \itemSize^\itemIndex_2, \itemSize^\itemIndex_3\} = \{\itemSize^\itemIndexTwo_1, \itemSize^\itemIndexTwo_2, \itemSize^\itemIndexTwo_3\}$.

Then we compute groups of \textit{similar items} based on the possible item heights and the item stackabilities.
Two items  $\itemIndex, \itemIndexTwo$ are considered similar if
$(i)$ they can be tilted in such a way that they have the same height and
$(ii)$ they have the same stackability characteristic $\itemStackability_\itemIndex = \itemStackability_\itemIndexTwo$

We create a group of similar items $\similarItemsSet$ for each possible item height and  stackability indication~$\itemStackability\in\{0,1\}$:
\begin{align*}
	\similarItemsSet = \left\{\itemIndex \in \itemSet : \itemStackability_\itemIndex=\itemStackability \wedge \left[ \size_3^\itemIndex = \similarItemHeight \vee  ( \text{item $\itemIndex$ is tiltable } \wedge (\size_1^\itemIndex = \similarItemHeight \vee \size_2^\itemIndex = \similarItemHeight) )\right] \right\}.
\end{align*}
By definition, all members of a given identical item group belong to the same similar item groups.
Hence, we will treat each $\similarItemsSet$ as a set of identical item groups instead of a set of items, \ie, $\similarItemsSet \subseteq \mathcal{P}(\itemSet)$ instead of $\similarItemsSet \subseteq \itemSet$, where $\mathcal{P}(\itemSet)$ denotes the power set of $\itemSet$.
Note that depending on the dimensions and tiltability, an identical item group can belong to up to three different sets~$\similarItemsSet$.

\subsubsection{Sorting}\label{sec:sorting}

We select one of the following five sorting criteria and sort $(i)$ the set of all similar item groups and $(ii)$ within each similar item group $\similarItemsSet$ the set of identical item groups accordingly in non-ascending order:
\begin{sloppypar}
\begin{itemize}[noitemsep,topsep=0pt,label={--}]
\item \sorting{cumulated volume}: The cumulative volume of all items in a group.
\item \sorting{highest volume}: The volume of the group member with the highest volume.
\item \sorting{stackability--cumulated volume}: Items are sorted lexicographically based on \sorting{stackability} first (stackable items are preferred over non-stackable ones) and on \sorting{cumulated volume} second.
\item \sorting{stackability--highest volume}: Analogous to \sorting{stackability--cumulated volume} but using \sorting{highest volume} as second order criterion.
\item \sorting{random}.
\end{itemize} 
\end{sloppypar}
Additionally, we introduce the parameter $\degreeOfRandomization \in (0, 1]$, which allows to define a degree of randomization applied to the item sorting.
Based on that we rearrange both sorted sets of item groups.
Given a sorted list of item groups $\allItemGroups = [\itemGroup_1, \itemGroup_2, \ldots, \itemGroup_\numItemGroups]$, we create a new sorting $\allNewItemGroups = [\newItemGroup_1, \newItemGroup_2, \ldots, \newItemGroup_\numItemGroups]$.
Item groups are drawn one-by-one from the original sorting $\allItemGroups$, where the $\iterator$-th item group $\newItemGroup_\iterator$ corresponds to the item group at position $\left\lceil \randomNumber_\iterator^{\frac{1}{\degreeOfRandomization}}(\numItemGroups-\iterator+1) \right\rceil$ in $\allItemGroups\setminus [\newItemGroup_1, \newItemGroup_2, \ldots, \newItemGroup_{\iterator-1}]$, where $\randomNumber_\iterator \in [0,1]$ is a uniformly distributed random number.
This is an adaption of the random distribution proposed by \cite{Shaw1997}, which is often used in metaheuristics.
The higher the parameter~$\degreeOfRandomization$, the more randomized the sorting is.
For $\degreeOfRandomization\rightarrow0$, the randomness is close to zero, while the maximum value of $\degreeOfRandomization=1$ corresponds to complete randomness.

For the sorting policies \sorting{stackability--cumulated volume} and \sorting{stackability--highest volume}, we maintain the strict prioritization of stackable item groups over non-stackable ones, \ie, the randomization is applied separately to the sorted (sub-)sets of stackable and non-stackable item groups.

If we now unpack the sorted item groups, we end up with an ordered list of the individual items defining the loading order:
\begin{align*}
\underbrace{\overbrace{\itemIndex_{\ell_1},\itemIndex_{\ell_1+1}, \ldots}^\text{IdenticalItems 1}, \overbrace{\itemIndex_{\ell_2}, \itemIndex_{\ell_2+1}, \ldots}^\text{IdenticalItems 2}, \ldots}_\text{SimilarItems 1},
\underbrace{\overbrace{\itemIndex_{\ell_\iterator},\itemIndex_{\ell_\iterator +1}, \ldots}^\text{IdenticalItems $\iterator$}, \ldots}_\text{SimilarItems 2}, \ldots
\end{align*}
In this ordered list, each item can appear up to three times, depending on the dimensions and tiltability.

\subsubsection{Determining item orientation}
As the idea is to load items with the same height next to each other, we need to ensure that we not only consecutively process all items contained in the same similar item group, but also that the items's orientations are aligned with the heights defined by the respective similar item groups.
That means, for an item $\itemIndex \in \similarItemsSet$, we need to determine the possible orientations~$\orientations$ such that $\itemSize^\itemIndex_3 = \similarItemHeight$.
To do so, it suffices to consider one tilt:
If the item is not tiltable, there is nothing to show.
If the item is tiltable, we consider the following three cases.
\begin{enumerate}[noitemsep,topsep=0pt,label=(\roman*)]
\item The item size is the same in all dimensions. All tilts yield the same item size, we can just pick one.
\item The item size is different in all dimensions. There can only be one tilt where the item height matches the similar item group's height~$\similarItemHeight$.
\item Two dimensions of the item size are the same.
\begin{itemize}[noitemsep,topsep=0pt,label=--]
\item The one unique dimension matches the height $\similarItemHeight$. There is only one applicable tilt.
\item The two identical dimensions match the height~$\similarItemHeight$. We can just pick one of the two tilts where the item height matches $\similarItemHeight$. Since we assume each tiltable item to be rotatable (see Section~\ref{sec:problem_definition}), for either tilt we can achieve the same item length and width.
\end{itemize}
\end{enumerate}
Now, the sorted item list from Section~\ref{sec:sorting} combined with the corresponding relevant item orientations, constitutes list $\listItemsOrientations$ from line~\ref{alg:create_sorted_groups} of Algorithm~\ref{alg:insertion_heuristic} with
\begin{align*}
\listItemsOrientations=\left( \underbrace{\overbrace{(\itemIndex_{\ell_1},\orientations_{\ell_1}),(\itemIndex_{\ell_1+1},\orientations_{\ell_1+1}), \ldots}^\text{IdenticalItems 1}, \overbrace{(\itemIndex_{\ell_2},\orientations_{\ell_2}), (\itemIndex_{\ell_2+1},\orientations_{\ell_2+1}), \ldots}^\text{IdenticalItems 2}, \ldots}_\text{SimilarItems 1},
\underbrace{\overbrace{(\itemIndex_{\ell_\iterator},\orientations_{\ell_\iterator}),(\itemIndex_{\ell_\iterator +1},\orientations_{\ell_\iterator+1}), \ldots}^\text{IdenticalItems $\iterator$}, \ldots}_\text{SimilarItems 2}, \ldots \right).
\end{align*}

\subsection{Selecting and potentially moving the next extreme point}\label{sec:choose_extreme_point}

\begin{figure}[H]
	\begin{subfigure}[t]{0.32\textwidth}
		\centering
		\begin{tikzpicture}
			[tdplot_main_coords,
			cube/.style={black},
			scale=0.8]
			\draw[gray] (0,0,0) -- (2,0,0);
			\draw[gray] (0,0,0) -- (0,2.75,0);
			\draw[gray] (0,0,0) -- (0,0,2.25);
			\foreach \x/\y in {0/0.5,0.5/1}{
				\foreach \l/\r in {0/0.5,0.5/1,1/1.5,1.5/2}{
					\mycuboid{0}{\x}{\l}{0.5}{\y}{\r}{white}{1};
				}
			}
			\foreach \l/\r in {0/0.5,0.5/1,1/1.5}{
				\mycuboid{0}{1}{\l}{0.5}{1.5}{\r}{white}{1};
			}
			\draw[gray] (2,0,0) -- (2,2.75,0);
			\draw[gray] (2,0,0) -- (2,0,2.25);
			\draw[gray] (0,0,2.25) -- (2,0,2.25);
			\draw[gray] (0,0,2.25) -- (0,2.75,2.25);
			\draw[gray] (0,2.75,0) -- (2,2.75,0);
			\draw[gray] (0,2.75,0) -- (0,2.75,2.25);
			\draw[gray] (0,2.75,2.25) -- (2,2.75,2.25);
			\draw[gray] (2,0,2.25) -- (2,2.75,2.25);
			\draw[gray] (2,2.75,0) -- (2,2.75,2.25);
		\end{tikzpicture}
		\caption{Load in a manner of stacks by sorting the extreme points in ascending lexicographic order by $x,y,z$.}
	\end{subfigure}
	\hfill
	\begin{subfigure}[t]{0.32\textwidth}
		\centering
		\begin{tikzpicture}
			[tdplot_main_coords,
			cube/.style={black},
			scale=0.8]
			\draw[gray] (0,0,0) -- (2,0,0);
			\draw[gray] (0,0,0) -- (0,2.75,0);
			\draw[gray] (0,0,0) -- (0,0,2.25);
			\foreach \l/\r in {0/0.5,0.5/1,1/1.5}{
				\mycuboid{0}{0}{\l}{0.5}{0.5}{\r}{white}{1};
			}
			\foreach \l/\r in {0/0.5,0.5/1}{
				\mycuboid{0}{0.5}{\l}{0.5}{1}{\r}{white}{1};
			}
			\foreach \x/\y in {0/0.5,0.5/1}{
				\foreach \l/\r in {0/0.5,0.5/1}{
					\mycuboid{0.5}{\x}{\l}{1}{\y}{\r}{white}{1};
				}
			}
			\mycuboid{0}{1}{0}{0.5}{1.5}{0.5}{white}{1};
			\mycuboid{1}{0}{0}{1.5}{0.5}{0.5}{white}{1};
			\draw[gray] (2,0,0) -- (2,2.75,0);
			\draw[gray] (2,0,0) -- (2,0,2.25);
			\draw[gray] (0,0,2.25) -- (2,0,2.25);
			\draw[gray] (0,0,2.25) -- (0,2.75,2.25);
			\draw[gray] (0,2.75,0) -- (2,2.75,0);
			\draw[gray] (0,2.75,0) -- (0,2.75,2.25);
			\draw[gray] (0,2.75,2.25) -- (2,2.75,2.25);
			\draw[gray] (2,0,2.25) -- (2,2.75,2.25);
			\draw[gray] (2,2.75,0) -- (2,2.75,2.25);
		\end{tikzpicture}
		\caption{Load around the origin by selecting the point closest to the origin as next extreme point.}
	\end{subfigure}
	\hfill
	\begin{subfigure}[t]{0.32\textwidth}
		\centering
		\begin{tikzpicture}
			[tdplot_main_coords,
			cube/.style={black},
			scale=0.8]
			\draw[gray] (0,0,0) -- (2,0,0);
			\draw[gray] (0,0,0) -- (0,2.75,0);
			\draw[gray] (0,0,0) -- (0,0,2.25);
			\foreach \x/\y in {0/0.5,0.5/1}{
				\foreach \l/\r in {0/0.5,0.5/1,1/1.5,1.5/2}{
					\mycuboid{\x}{\l}{0}{\y}{\r}{0.5}{white}{1};
				}
			}
			\foreach \l/\r in {0/0.5,0.5/1,1/1.5}{
				\mycuboid{1}{\l}{0}{1.5}{\r}{0.5}{white}{1};
			}
			\draw[gray] (2,0,0) -- (2,2.75,0);
			\draw[gray] (2,0,0) -- (2,0,2.25);
			\draw[gray] (0,0,2.25) -- (2,0,2.25);
			\draw[gray] (0,0,2.25) -- (0,2.75,2.25);
			\draw[gray] (0,2.75,0) -- (2,2.75,0);
			\draw[gray] (0,2.75,0) -- (0,2.75,2.25);
			\draw[gray] (0,2.75,2.25) -- (2,2.75,2.25);
			\draw[gray] (2,0,2.25) -- (2,2.75,2.25);
			\draw[gray] (2,2.75,0) -- (2,2.75,2.25);
		\end{tikzpicture}
		\caption{Load the \ULD\ from bottom to top by sorting the extreme points in ascending lexicographic order by $z,y,x$.}
	\end{subfigure}
	\caption{Different loading strategies.}\label{fig:sorting strategies}
\end{figure}
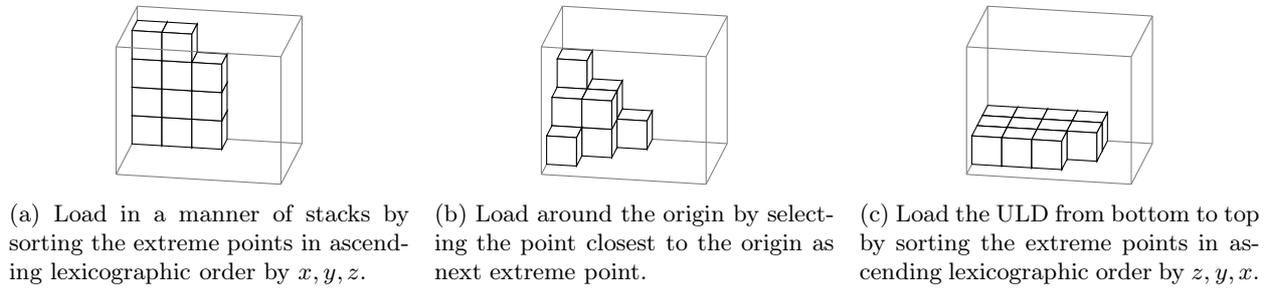

Figure~\ref{fig:sorting strategies} presents three different loading strategies we tested that affect the selection of the next extreme point.
To ensure stable loading of a \ULD\ in accordance with the item selection process from Section~\ref{sec:choose_item}, the loading from bottom to top is preferred over the other ones.
The set of extreme points is therefore sorted in ascending lexicographic order by $z,y,x$ and the next extreme point is chosen by iterating over the sorted list.

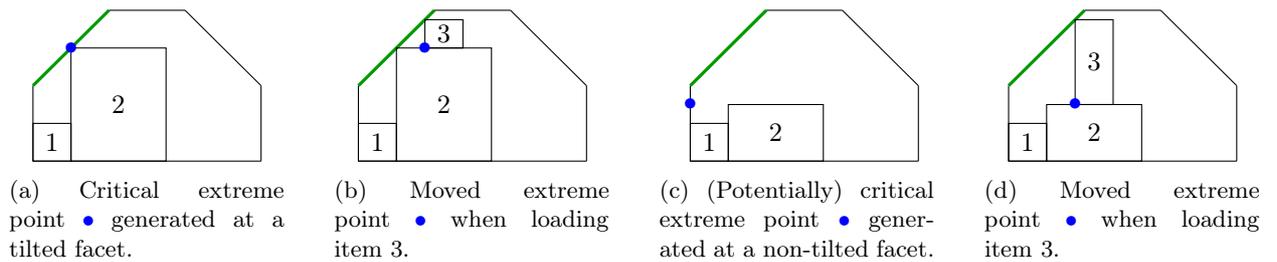
\begin{figure}[H]
	\begin{subfigure}[t]{0.22\textwidth}
		\centering
		\begin{tikzpicture}
			[scale=0.5]
			\draw[draw=black] (0,0) -- (6,0);
			\draw[draw=black] (0,0) -- (0,2);
			\draw[draw=mygreen,very thick] (0,2) -- (2,4);
			\draw[draw=black] (2,4) -- (4,4);
			\draw[draw=black] (4,4) -- (6,2);
			\draw[draw=black] (6,0) -- (6,2);
			\foreach  \x/\y/\w/\h/\n in {0/0/1/1/1,1/0/2.5/3/2}{
				\draw[draw=black] (\x,\y) rectangle ++(\w,\h);
				\node[] at (\x + \w/2,\y+\h/2) {\n};
			}
			\node at (1,3){\blue{\textbullet}};
		\end{tikzpicture}
		\caption{Critical extreme point~\blue{\textbullet} generated at a tilted facet.}
		\label{fig:critical_extreme_points1}
	\end{subfigure}
	\hfill
	\begin{subfigure}[t]{0.22\textwidth}
		\centering
		\begin{tikzpicture}
			[scale=0.5]
			\draw[draw=black] (0,0) -- (6,0);
			\draw[draw=black] (0,0) -- (0,2);
			\draw[draw=mygreen,very thick] (0,2) -- (2,4);
			\draw[draw=black] (2,4) -- (4,4);
			\draw[draw=black] (4,4) -- (6,2);
			\draw[draw=black] (6,0) -- (6,2);
			\foreach  \x/\y/\w/\h/\n in {0/0/1/1/1,1/0/2.5/3/2,1.75/3/1/0.75/3}{
				\draw[draw=black] (\x,\y) rectangle ++(\w,\h);
				\node[] at (\x + \w/2,\y+\h/2) {\n};
			}
			\node at (1.75,3){\blue{\textbullet}};
		\end{tikzpicture}
		\caption{Moved extreme point~\blue{\textbullet} when loading item~3.}
		\label{fig:critical_extreme_points2}
	\end{subfigure}
	\hfill
	\begin{subfigure}[t]{0.22\textwidth}
		\centering
		\begin{tikzpicture}
			[scale=0.5]
			\draw[draw=black] (0,0) -- (6,0);
			\draw[draw=black] (0,0) -- (0,2);
			\draw[draw=mygreen,very thick] (0,2) -- (2,4);
			\draw[draw=black] (2,4) -- (4,4);
			\draw[draw=black] (4,4) -- (6,2);
			\draw[draw=black] (6,0) -- (6,2);
			\foreach  \x/\y/\w/\h/\n in {0/0/1/1/1,1/0/2.5/1.5/2}{
				\draw[draw=black] (\x,\y) rectangle ++(\w,\h);
				\node[] at (\x + \w/2,\y+\h/2) {\n};
			}
			\node at (0,1.5){\blue{\textbullet}};
		\end{tikzpicture}
		\caption{(Potentially) critical extreme point~\blue{\textbullet} generated at a non-tilted facet.}
		\label{fig:critical_extreme_points3}
	\end{subfigure}
	\hfill
	\begin{subfigure}[t]{0.22\textwidth}
		\centering
		\begin{tikzpicture}
			[scale=0.5]
			\draw[draw=black] (0,0) -- (6,0);
			\draw[draw=black] (0,0) -- (0,2);
			\draw[draw=mygreen,very thick] (0,2) -- (2,4);
			\draw[draw=black] (2,4) -- (4,4);
			\draw[draw=black] (4,4) -- (6,2);
			\draw[draw=black] (6,0) -- (6,2);
			\foreach  \x/\y/\w/\h/\n in {0/0/1/1/1,1/0/2.5/1.5/2,1.75/1.5/1/2.25/3}{
				\draw[draw=black] (\x,\y) rectangle ++(\w,\h);
				\node[] at (\x + \w/2,\y+\h/2) {\n};
			}
			\node at (1.75,1.5){\blue{\textbullet}};
		\end{tikzpicture}
		\caption{Moved extreme point~\blue{\textbullet} when loading item~3.}
		\label{fig:critical_extreme_points4}
	\end{subfigure}
	\caption{Extreme points at which (potentially) no item can be loaded (\ref{fig:critical_extreme_points1} and \ref{fig:critical_extreme_points3}) and the corresponding moved extreme points (\ref{fig:critical_extreme_points2} and \ref{fig:critical_extreme_points4}). The critical facet is colored \mygreen{green}/marked in \textbf{bold}.}\label{fig:critical_extreme_points}
\end{figure}

In case the \ULD\ has tilted facets at the top (see Figure~ \ref{fig:ULD_shapes:two_upper_cuts}), we may generate extreme points at which no item can ever be loaded (for a detailed description of how extreme points are generated, we refer to Section~\ref{sec:generating_extreme_points}).
An example is shown in Figure~\ref{fig:critical_extreme_points1}.
After loading item 2, we generate the extreme point~\blue{\textbullet}.
As each item has a non-zero height, clearly no item can ever be loaded at this position.
Similarly, when item 2 is loaded in the way depicted in Figure~\ref{fig:critical_extreme_points3}, only very few items could ever be loaded at the extreme point~\blue{\textbullet}.
We call such an extreme point \emph{critical extreme point} and the corresponding tilted facet limiting the loading space for the extreme point  \emph{critical facet}.

To avoid having these useless extreme points and unused space in the \ULD, we introduce a procedure to move them away from the \ULD\ wall according to the size of the item that is supposed to be loaded at the given point. 

An extreme point is only allowed to be moved if it is directly on a \ULD\  facet (tilted or non-tilted). 
Moreover, we only allow moving away from the origin as the way we generate extreme points (see Section~\ref{sec:generating_extreme_points}) already ensures sufficient alternative extreme points towards the origin.
Considering the relevant \ULD\ shapes (Figure~\ref{fig:ULD_shapes}), moving in $x$-direction would only move an extreme point parallel to any potential critical facet and can be neglected.
Finally, moving an extreme point in $z$-direction (relevant for \ULD\ shapes (e) and (f) in Figure~\ref{fig:ULD_shapes}) very likely results in loading positions at which the item would float.
Hence, the only way in which we allow an extreme point to be moved is away from the origin in $y$-direction.
Again, given the relevant \ULD\ shapes, this implies that an extreme point only needs to be considered for moving if it lies on one of the two facets with critical extreme points depicted in Figures~\ref{fig:critical_extreme_points1} and \ref{fig:critical_extreme_points3}.
More formally, if it lies on a facet with plane equation~$\normalVector_1 x+\normalVector_2 y+\normalVector_3 z=\planeEquationOffset$, where $\normalVector_1=0, \normalVector_2 >0, \normalVector_3 \leq 0$.

If we have such a critical extreme point $\extremePoint=(\extremePoint_1, \extremePoint_2, \extremePoint_3)$, the corresponding critical facet is always given by
$\normalVector_1 x+\normalVector_2 y+\normalVector_3 z=\planeEquationOffset$ with $\normalVector_1=0, \normalVector_2 >0, \normalVector_3 < 0$.
To compute the minimum value~$\shiftValue$ by which we have to move an item of size~$\size$ so that it fits at extreme point~$\extremePoint$, we need to place the top left corner of the item directly at the critical facet (see Figure~\ref{fig:critical_extreme_points2}~and~\ref{fig:critical_extreme_points4}).
Thus, we solve
$\normalVector_1(\extremePoint_1 + \size_1) + \normalVector_2(\extremePoint_2+\shiftValue) + \normalVector_3(\extremePoint_3 + \size_3)=\planeEquationOffset$
and obtain
$$\shiftValue = \frac{\planeEquationOffset - \normalVector_2 \extremePoint_2-\normalVector_3(\extremePoint_3+\size_3)}{\normalVector_2}.$$
Now, if $\shiftValue>0$, we know that we have to move the extreme point in order for the item to fit into the ULD (regarding the critical facet).
The new moved extreme point is $(\extremePoint_1, \extremePoint_2 + \shiftValue, \extremePoint_3)$.

\begin{figure}
	\newcommand\widthFactor{0.47}
	\begin{subfigure}[t]{\widthFactor\textwidth}
		\centering
		\begin{tikzpicture}
			[tdplot_main_coords,
			cube/.style={black},
			scale=0.8]
			\draw[gray] (0,0,0) -- (1.5,0,0);
			\draw[gray] (0,0,0) -- (0,3,0);
			\draw[gray] (0,0,0) -- (0,0,2);
			\draw[mygreen,very thick] (1.5,0,0) -- (2,0.5,0);
			\draw[gray] (2,0.5,0) -- (2,3,0);
			\draw[mygreen,very thick] (2,0.5,0) -- (2,0.5,2);
			\draw[mygreen,very thick] (1.5,0,0) -- (1.5,0,2);
			\draw[gray] (0,0,2) -- (1.5,0,2);
			\draw[mygreen,very thick] (1.5,0,2) -- (2,0.5,2);
			\draw[gray] (0,3,0) -- (2,3,0);
			\draw[gray] (0,3,0) -- (0,3,2);
			\draw[gray] (2,3,0) -- (2,3,2);
			\draw[gray] (0,0,2) -- (0,3,2);
			\draw[gray] (0,3,2) -- (2,3,2);
			\draw[gray] (2,3,2) -- (2,0.5,2);
		\end{tikzpicture}\\
		\caption{Critical facet with $\normalVector_1<0,  \normalVector_2 >0, \normalVector_3=0$ and move direction~$y$.} \label{fig:shifting_extreme_points1}
	\end{subfigure}
	\hfill
	\begin{subfigure}[t]{\widthFactor\textwidth}
		\centering
		\begin{tikzpicture}
			[tdplot_main_coords,
			cube/.style={black},
			scale=0.8]
			\draw[mygreen,very thick] (0.5,3,0) -- (0.5,3,2);
			\draw[mygreen,very thick] (0,2.5,0) -- (0,2.5,2);
			\draw[mygreen,very thick] (0,2.5,0) -- (0.5,3,0);
			\draw[mygreen,very thick] (0,2.5,2) -- (0.5,3,2);
			\draw[gray] (0,0,0) -- (2,0,0);
			\draw[gray] (0,0,0) -- (0,2.5,0);
			\draw[gray] (0,0,2) -- (0,2.5,2);
			\draw[gray] (0,0,0) -- (0,0,2);
			\draw[gray] (2,0,0) -- (2,3,0);
			\draw[gray] (2,0,2) -- (2,3,2);
			\draw[gray] (2,0,0) -- (2,0,2);
			\draw[gray] (0,0,2) -- (2,0,2);
			\draw[gray] (0.5,3,0) -- (2,3,0);
			\draw[gray] (0.5,3,2) -- (2,3,2);
			\draw[gray] (2,3,0) -- (2,3,2);
		\end{tikzpicture}\\
		\caption{Critical facet with $\normalVector_1>0,  \normalVector_2 <0, \normalVector_3=0$ and move direction~$x$.}\label{fig:ULD_shapes:shifting_extreme_points2}
	\end{subfigure}
	\caption{Other shapes where the extreme points can be moved. The critical facet is colored \mygreen{green}/marked in \textbf{bold}.} \label{fig:shifting_extreme_points}
\end{figure}
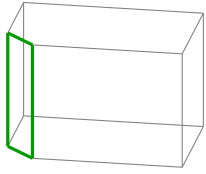
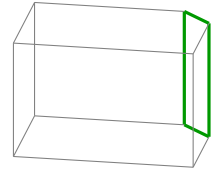

Note that the approach can be extended to other shapes, not relevant for the air freight context, as shown in Figure~\ref{fig:shifting_extreme_points}.

\subsection{Checking whether an item can be loaded at an extreme point}\label{sec:can_be_loaded}

An item can be loaded at an extreme point if it does not exceed the \ULD\ walls, it does not collide with other items, it is supported (non-floating), and it is not stacked on a non-stackable item.
In the following, we will explain how the individual checks are performed.
The section closes with a description of an acceleration technique.

\subsubsection{\ULD\ wall exceedance}
Each \ULD\ facet $\normalVector_1 x+\normalVector_2 y+\normalVector_3 z=\planeEquationOffset$ defines a half-space corresponding to the inside of the \ULD\ via
$\{(x,y,z) \in \mathbb{R}^3: \normalVector_1 x+\normalVector_2 y+\normalVector_3 z \geq \planeEquationOffset\}.$
To check whether an item loaded at a given extreme point lies within the \ULD , we need to check for each \ULD\ facet if all of the item's corner points lie within that half-space.
We say that the item lies on the inside of that facet.
We can apply knowledge about the shape of the \ULD/facet and the way that we generate extreme points to simplify and speed up these checks.

Let's consider an item of size $\size=(\size_1, \size_2, \size_3)$ that is supposed to be inserted at position $\extremePoint=(\extremePoint_1, \extremePoint_2, \extremePoint_3)$.
First, we assume that we always generate valid extreme points, \ie, points within the \ULD.
Hence, we never need to check whether the item's corner point corresponding to $\extremePoint$ lies within the \ULD.
Secondly, we can distinguish between tilted and non-tilted facets.
For the latter case, we can perform the checks for all facets of this kind at once:
We simply compute the \ULD's bounding box, \ie, the smallest cuboid that covers the \ULD.
We assume that, by definition, the \ULD\ is shifted towards the origin
in such a way that the minimum $x$-, $y$-, and $z$-coordinate over all \ULD\ vertices are 0.
Hence, we can store the bounding box as its size $\boundingBox=(\boundingBox_1, \boundingBox_2, \boundingBox_3)$ along the $x$-, $y$-, and $z$-axis.
Now it is sufficient to check whether $\extremePoint_\itemIndex+\size_\itemIndex \leq \boundingBox_\itemIndex$ for $\itemIndex\in\{1,2,3\}$ to ensure that the item lies on the inner side of all non-tilted facets.

For a tilted facet defined by the plane equation~$\normalVector_1 x+\normalVector_2 y+\normalVector_3 z=\planeEquationOffset$, an item lies on the inside of or on the facet if $\normalVector \cdot (\extremePoint_1 + \delta_1 \size_1, \extremePoint_2 + \delta_2 \size_2, \extremePoint_3 + \delta_3 \size_3) \geq \planeEquationOffset$,
where $\delta_\dimension =1$ if $\normalVector_\dimension < 0$ and $\delta_\dimension =0$ otherwise, for $\dimension\in\dimensions$.

\subsubsection{Collision check}

Two items~$\itemIndex$ and $\itemIndexTwo$ of size $\size^\itemIndex$ and $\size^\itemIndexTwo$ loaded at positions $\extremePoint^\itemIndex$ and $\extremePoint^\itemIndexTwo$ collide if $\extremePoint^\itemIndexTwo_\dimension < \extremePoint^\itemIndex_\dimension + \size^\itemIndex_\dimension$ and $\extremePoint^\itemIndex_\dimension < \extremePoint^\itemIndexTwo_\dimension + \size^\itemIndexTwo_\dimension$ for all $\dimension\in\dimensions$.
For the collision check, we just iterate over already loaded items.
If any loaded item~$\itemIndexTwo$ collides with the item~$\itemIndex$ to be loaded, item~$\itemIndex$ cannot be loaded at the position in the given orientation.

\subsubsection{Non-floating and stackability check}

The checks whether an item $\itemIndex$ is floating or placed on a non-stackable item when loaded at position $\extremePoint$ are performed in combination in Algorithm~\ref{alg:floating_stackability}.

The set of loaded items is extended by an artificial item to model the \ULD's base area.
In line~\ref{alg:reducedLoadedItems}, we first compute a reduced set of loaded items $\reducedLoadedItems$ that are relevant for supporting $\itemIndex$: 
These items have to overlap with $\itemIndex$ in the $x$- and $y$-direction, \ie, their base areas overlap.
Moreover, their surface area has to be at an height~$\surfaceAreaHeight$ suitable for direct support or indirect support via padding material, \ie, $\surfaceAreaHeight \in [\extremePoint_3 - \paddingHeight, \extremePoint_3]$.
This includes a dummy item of height~$0$ modelling the \ULD\ floor if applicable according to the previous condition.
Subsequently, we sort this reduced item set by non-ascending $z$ end position.
We then iterate over the ordered set to determine whether the item $i$ has any direct support, its total supported area, and the number of supported corner points.
If the item $\itemIndex$ directly rests on a non-stackable item, stackability requirements are violated (see line~\ref{alg:stackability_violated}) and we can stop.
If it directly rests on a stackable item, the direct support criterion is fulfilled (see line~\ref{alg:directly_supported}).
In lines \ref{alg:floating_stackability:directSupport}--\ref{alg:floating_stackability:total_supported_area}, the total area and the number of corner points supported by stackable items are updated.
The non-floating and stackability check returns 'true', if the item~$\itemIndex$ to load is directly supported and $(i)$ the four bottom corner points are supported or $(ii)$ at least $\minimumItemOverlap$ percent of the item's base area is supported.

In case of a padding height~$\paddingHeight > 0$, calculating the total supported area (see lines~\ref{alg:add_supported_area}-\ref{alg:subtract_supported_area}) results in the following problem.
Since $\paddingHeight > 0$, we can have multiple items in $\reducedLoadedItems$ stacked on top of each other, \ie, we cannot simply compute the supported area of $\itemIndex$ for each individual loaded item $\loadedItemIndex \in \reducedLoadedItems$ and take the sum.
Instead, we would have to calculate the union of the base areas of all stackable items in $\reducedLoadedItems$, subtract the area potentially blocked by non-stackable items in $\reducedLoadedItems$ on top, and then finally compute the overlap with $\itemIndex$.
Solving this problem is computationally expensive.
Therefore, we apply a simplification:
If the current iteration's item $\loadedItemIndex$ is stackable, we calculate the overlap of its base area with $\itemIndex$'s base area and add it to the cumulated supported area (see line~\ref{alg:add_supported_area}).
Subsequently, for each item $\loadedItemIndexTwo$ considered in the previous iterations, we calculate the overlap of the base areas of $\loadedItemIndex$, $\loadedItemIndexTwo$, and $\itemIndex$ and subtract it from the cumulated supported area (see line \ref{alg:subtract_supported_area}).
This way, we assure that item $\loadedItemIndexTwo$'s base area is not counted towards the support for item $\itemIndex$ if there is another item $\loadedItemIndex$ between $\loadedItemIndexTwo$ and $\itemIndex$ (vertically).
However, this procedure only works accurately if there are no more than two potential support items in $\reducedLoadedItems$  stacked on top of each other.

\begin{algorithm}
	\DontPrintSemicolon
	\SetKwInput{Input}{Input}
	\SetKwInput{Output}{Output}
	\newcommand\mycommfont[1]{\footnotesize\ttfamily\textcolor{mygreen}{#1}}
	\SetCommentSty{mycommfont}
	\SetNoFillComment
	\Input{Loaded items~$\loadedItems$, new item~$\itemIndex$ to load at position $\extremePoint^\itemIndex$}
	\Output{Whether the item~$\itemIndex$ is not floating and meets the stackability requirements}
	Add artificial item to $\loadedItems$ to model the \ULD's base area\;
	Define $\reducedLoadedItems=\{ \loadedItemIndex\in\loadedItems : \text{item } \loadedItemIndex \text{ intersects with } (\extremePoint^\itemIndex_1, \extremePoint^\itemIndex_1 + \itemSize_1^i) \times (\extremePoint^\itemIndex_2, \extremePoint^\itemIndex_2 + \itemSize_2^i) \times [\extremePoint^\itemIndex_3 - \paddingHeight, \extremePoint^\itemIndex_3] \} \subseteq \loadedItems$\;\label{alg:reducedLoadedItems}
	Sort loaded items by non-ascending end position $\extremePoint^\loadedItemIndex_3 + \size^\loadedItemIndex_3$ to get $\reducedLoadedItems =[ \loadedItemIndex_1,\loadedItemIndex_2, \ldots ]$\;
	$ \variable{directlySupported} = \text{false}$\;
	$ \variable{totalSupportedArea} = 0$\;
	$ \variable{numberSupportedCornerPoints} = 0$\;
	\For{$\loadedItemIndex \in \reducedLoadedItems$}{
		\If{$\loadedItemIndex$ is not stackable}{
			\If{$\itemIndex$ directly rests on $\loadedItemIndex$}{\Return false \label{alg:stackability_violated}}
			continue}
		\If{$\itemIndex$ directly rests on $\loadedItemIndex$\label{alg:floating_stackability:directSupport}}{
			Update \variable{numberSupportedCornerPoints}\;
			$\variable{directlySupported} = \text{true}$\label{alg:directly_supported}}
		$\variable{additionalSupportedArea} = \text{baseAreaOverlap}(\loadedItemIndex, \itemIndex)$\;\label{alg:add_supported_area}
		\For{$\loadedItemIndexTwo \in \{ \loadedItemIndex_1,\loadedItemIndex_2, \ldots, \loadedItemIndex \}$}{
			$\variable{additionalSupportedArea} \minuseq \text{baseAreaOverlap}(\loadedItemIndex, \loadedItemIndexTwo, \itemIndex)$\;\label{alg:subtract_supported_area}
		}
		\If{$\textnormal{\variable{additionalSupportedArea}} > 0 $}{$\variable{totalSupportedArea} \pluseq \variable{additionalSupportedArea}$}\label{alg:floating_stackability:total_supported_area}
	}
	\Return \variable{directlySupported} and ($\variable{numberSupportedCornerPoints}=4$ or \variable{totalSupportedArea}  $\geq \minimumItemOverlap \cdot \itemSize_1^\itemIndex \cdot \itemSize_2^\itemIndex $)\;
	\caption{Non-floating and stackability check.}\label{alg:floating_stackability}
\end{algorithm}

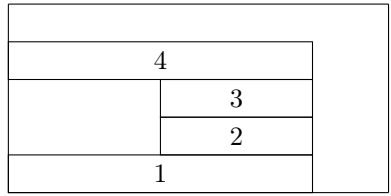
\begin{figure}[]
	\centering
		\begin{tikzpicture}[scale=1.0]
			\draw[draw=black] (0,0) -- (5,0);
			\draw[draw=black] (0,0) -- (0,2.5);
			\draw[draw=black] (0,2.5) -- (5,2.5);
			\draw[draw=black] (5,0) -- (5,2.5);
			\foreach  \x/\y/\w/\h/\n in {0/0/4/0.5/1, 2/0.5/2/0.5/2, 2/1/2/0.5/3, 0/1.5/4/0.5/4}{
				\draw[draw=black] (\x,\y) rectangle ++(\w,\h);
				\node[] at (\x + \w/2,\y+\h/2) {\n};
			}
			
		\end{tikzpicture}
		\caption{Inaccuracy for stacked support items: Item 4 needs to be loaded, items 1-3 are already loaded. The padding height~$\paddingHeight$ is equal to the (vertical) distance between item 4 and item 1.}
		\label{fig:stackability_simplification}
\end{figure}

In Figure~\ref{fig:stackability_simplification}, we show an example where the calculated total supported area is too small.
If calculated accurately, the full base area of item 4 should be considered supported (directly or indirectly).
To calculate the support, we iterate over the items 3, 2, and 1 (in that order).
In the first iteration, we add the base area of item 3.
In the second iteration, we add the base area of item 2 and subtract the overlap of the base areas of item 2 and 3, which amounts to item 3's base area in total.
Finally, in the third iteration, we add the base area of item 1 and subtract the base areas of items 2 and 3.
Hence, we end up with an overall supported area equal to only the base area of item 3.

However, since in realistic instances the padding height~$\paddingHeight$ is typically significantly smaller than the item heights, the scenario outlined above only occurs very rarely.
Moreover, while our simplification might sometimes rule out otherwise feasible solutions, it never allows for infeasible solutions to be accepted:
We only ever overestimate the blocked supported area, \ie, underestimate the actual item support.
Hence, we consider the simplification acceptable.

\subsubsection{Acceleration techniques}
\label{sec:acceleration_techniques}

\begin{figure}
	\centering
	\begin{tikzpicture}
		[tdplot_main_coords,
		cube/.style={black},
		scale=0.5]
		\draw[gray] (0,0,0) -- (5,0,0);
		\draw[gray] (0,0,0) -- (0,5,0);
		\draw[gray] (0,0,0) -- (0,0,5);
		\draw[gray] (5,0,0) -- (5,5,0);
		\draw[gray] (5,0,0) -- (5,0,5);
		\draw[gray] (0,0,5) -- (5,0,5);
		\draw[gray] (0,0,5) -- (0,5,5);
		\draw[gray] (0,5,0) -- (5,5,0);
		\draw[gray] (0,5,0) -- (0,5,5);
		
		\mycuboid{0}{0}{0}{2}{2}{1}{white}{1};
		\mycuboid{2}{0}{0}{3}{1}{1}{white}{1};
		\mycuboid{0}{2}{0}{0.5}{2.5}{0.5}{white}{1};
		\mycuboid{0}{0}{1}{1.5}{1.5}{2}{white}{1};
		
		\draw[gray] (0,5,5) -- (5,5,5);
		\draw[gray] (5,0,5) -- (5,5,5);
		\draw[gray] (5,5,0) -- (5,5,5);
		
		\foreach \x in {0,1,...,5}{
			\draw[gray] (\x,0,5) -- (\x,5,5);
			\draw[gray] (\x,5,0) -- (\x,5,5);
			\draw[gray] (0,\x,5) -- (5,\x,5);
			\draw[gray] (5,\x,0) -- (5,\x,5);
			\draw[gray] (0,5,\x) -- (5,5,\x);
			\draw[gray] (5,0,\x) -- (5,5,\x);
		}
	\end{tikzpicture}
	\caption{A $5 \times 5 \times 5$ grid for a ULD.}\label{fig:grid}
\end{figure}
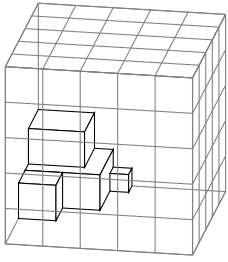

Performing the collision, non-floating, and stackability check for all already loaded items can be very time consuming, especially if the instance consists of many items.
However, for these checks, we only have to consider a small subset of the loaded items within close proximity to the position of the next item to be loaded.
To identify the relevant items quickly, in the preprocessing we divide the \ULD\ space into a grid of equally sized cubes:
First, the average item edge size $\bar{\size}$ (considering all 3 dimensions) is calculated as
\begin{align*}
	\gridCellSize = \frac{1}{3|\itemSet|}\sum_{\itemIndex}^{|\itemSet|} (\itemSize_1^\itemIndex + \itemSize_2^\itemIndex + \itemSize_3^\itemIndex).
\end{align*}
Then, a grid with cells of size $(\gridCellSize, \gridCellSize, \gridCellSize)$ is created.
Here, each cell is defined as a three-dimensional half-open interval.
For example, the grid cell at the origin is defined as $[0,\gridCellSize) \times [0,\gridCellSize) \times [0,\gridCellSize)$.
Figure~\ref{fig:grid} shows an example of such a grid.

Whenever an item item $\itemIndex$ of size $\itemSize^\itemIndex$ is loaded at position $\extremePoint^\itemIndex$, we register it in all grid cells that intersect with $[\extremePoint^\itemIndex_1, \extremePoint^\itemIndex_1 +\size_1^\itemIndex) \times [\extremePoint^\itemIndex_2, \extremePoint^\itemIndex_2+ \size_2^\itemIndex) \times [\extremePoint^\itemIndex_3, \extremePoint^\itemIndex_3 + \size_3^\itemIndex)$. 
Similarly, if we try to load a new item $\itemIndex$ at position $\extremePoint^\itemIndex$, we can easily determine the grid cells relevant for our checks.
For the collision check, these are again all grid cells which the item intersects with if loaded at the given position.
For the non-floating and stackability check, we need to identify all the items below $\itemIndex$ that can potentially provide support considering the maximum padding height~$\paddingHeight$.
These are precisely the items in the grid cells that intersect with $[\extremePoint_1,\extremePoint_1+\size_1^\itemIndex) \times [\extremePoint_2,\extremePoint_2+\size_2^\itemIndex) \times [\extremePoint_3 -\paddingHeight-\varepsilon,\extremePoint_3)$ with $\varepsilon >0$.

Determining the intersecting grid cells is computationally inexpensive:
Given an $\gridCellLimit_1 \times \gridCellLimit_2 \times \gridCellLimit_3$ grid consisting of the cells
\begin{align*}
[\gridCellIndex_1 \gridCellSize - \gridCellSize, \gridCellIndex_1 \gridCellSize ) \times [\gridCellIndex_2 \gridCellSize - \gridCellSize, \gridCellIndex_2 \gridCellSize ) \times [\gridCellIndex_3 \gridCellSize - \gridCellSize, \gridCellIndex_3 \gridCellSize ) \text{ for } \gridCellIndex \in \{1, ..., \gridCellLimit_1\} \times \{1, ..., \gridCellLimit_2\} \times \{1, ..., \gridCellLimit_3\}
\end{align*}
and item~$\itemIndex$, we can calculate for all $\dimension \in \dimensions$
\begin{align*}
\gridCellIndexMin_\dimension = \left \lfloor \frac{\extremePoint^\itemIndex_\dimension}{\gridCellSize} \right \rfloor \text { and }\gridCellIndexMax_\dimension = \left \lfloor \frac{\extremePoint^\itemIndex_\dimension + \size^\itemIndex_\dimension - \varepsilon}{\gridCellSize} \right \rfloor \text { with } \varepsilon > 0.
\end{align*}
The grid cells intersecting with item~$\itemIndex$ are then given by
\begin{align*}
[\gridCellIndex_1 \gridCellSize - \gridCellSize, \gridCellIndex_1 \gridCellSize ) \times [\gridCellIndex_2 \gridCellSize - \gridCellSize, \gridCellIndex_2 \gridCellSize ) \times [\gridCellIndex_3 \gridCellSize - \gridCellSize, \gridCellIndex_3 \gridCellSize ) \text{ for } \gridCellIndex \in \{\gridCellIndexMin_1, ..., \gridCellIndexMax_1\} \times \{\gridCellIndexMin_2, ..., \gridCellIndexMax_2\} \times \{\gridCellIndexMin_3, ..., \gridCellIndexMax_3\}.
\end{align*}
This way, we can significantly reduce the number of items to consider for our checks when loading a new item.

\subsection{Generating extreme points}\label{sec:generating_extreme_points}

\cite{CrainicEtAl2008} present an efficient way to determine candidate points where to place the items.
The calculation takes the already loaded cargo into account and aims to identify sufficiently many and suitable positions such that a dense packing can be achieved and sufficiently few positions to ensure that the problem can be solved within a reasonable time.
This set of extreme points is updated as soon as an item is added to a given packing.
In the following, we describe in detail how extreme points are generated.
We use the main ideas of \cite{CrainicEtAl2008} but adapt their approach with the advantage of generating more extreme points without significantly affecting the runtime.

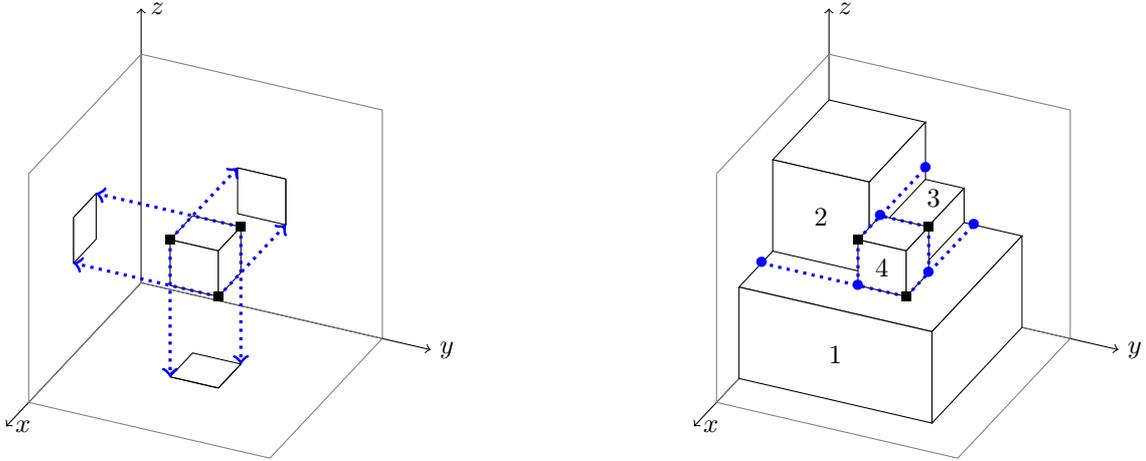
\begin{figure}
	\begin{subfigure}[t]{0.45\textwidth}
		\centering
		\tdplotsetmaincoords{60}{115}
		\begin{tikzpicture}
			[tdplot_main_coords,
			cube/.style={black},
			axis/.style={->,black},
			scale=0.7]
			
			\draw[axis] (0,0,0) -- (6,0,0) node[anchor=west]{$x$};
			\draw[axis] (0,0,0) -- (0,6,0) node[anchor=west]{$y$};
			\draw[axis] (0,0,0) -- (0,0,6) node[anchor=west]{$z$};
			
			\draw[gray] (0,0,0) -- (5,0,0);
			\draw[gray] (0,0,0) -- (0,5,0);
			\draw[gray] (0,0,0) -- (0,0,5);
			\draw[gray] (5,0,0) -- (5,5,0);
			\draw[gray] (5,0,0) -- (5,0,5);
			\draw[gray] (0,0,5) -- (5,0,5);
			\draw[gray] (0,0,5) -- (0,5,5);
			\draw[gray] (0,5,0) -- (5,5,0);
			\draw[gray] (0,5,0) -- (0,5,5);
			
			\mycuboid{2}{2}{0}{3}{3}{0}{white}{1};
			\mycuboid{2}{0}{2}{3}{0}{3}{white}{1};
			\mycuboid{0}{2}{2}{0}{3}{3}{white}{1};
			
			\mycuboid{2}{2}{2}{3}{3}{3}{white}{1};
			
			\draw[blue,dotted, very thick,->] (2, 3, 3) -- (2 ,0, 3);
			\draw[blue,dotted, very thick,->] (2, 3, 3) -- (2 , 3, 0);
			\draw[blue,dotted, very thick,->] (3, 2, 3) -- (0, 2, 3);
			\draw[blue,dotted, very thick,->] (3, 2, 3) -- (3, 2, 0);
			\draw[blue,dotted, very thick,->] (3 ,3 ,2) -- (0, 3, 2);
			\draw[blue,dotted, very thick,->] (3 ,3 ,2) -- (3, 0, 2);
			
			\foreach \Point in {(2, 3, 3), (3, 2, 3), (3 ,3 ,2)}{
				\node at \Point {$\smallFilledSquare$};
			}
		\end{tikzpicture}
		\caption{The six different projections for generating extreme points (for illustration purposes, we ignore that this is invalid as the item is floating). They result from projecting the black corner points~$\smallFilledSquare$ of the item along the orthogonal axes of the \ULD.}\label{fig:projection1}
	\end{subfigure}
	\hfill
	\begin{subfigure}[t]{0.45\textwidth}
		\centering
		\tdplotsetmaincoords{60}{115}
		\begin{tikzpicture}
			[tdplot_main_coords,
			cube/.style={black},
			axis/.style={->,black},
			scale=0.7]
			
			\draw[axis] (0,0,0) -- (6,0,0) node[anchor=west]{$x$};
			\draw[axis] (0,0,0) -- (0,6,0) node[anchor=west]{$y$};
			\draw[axis] (0,0,0) -- (0,0,6) node[anchor=west]{$z$};
			
			\draw[gray] (0,0,0) -- (5,0,0);
			\draw[gray] (0,0,0) -- (0,5,0);
			\draw[gray] (0,0,0) -- (0,0,5);
			\draw[gray] (5,0,0) -- (5,5,0);
			\draw[gray] (5,0,0) -- (5,0,5);
			\draw[gray] (0,0,5) -- (5,0,5);
			\draw[gray] (0,0,5) -- (0,5,5);
			\draw[gray] (0,5,0) -- (5,5,0);
			\draw[gray] (0,5,0) -- (0,5,5);
			
			\mycuboid{0}{0}{0}{4}{4}{2}{white}{1};
			\mycuboid{0}{0}{2}{2.5}{2}{4}{white}{1};
			\mycuboid{0}{2}{2}{2}{2.8}{2.75}{white}{1};
			\mycuboid{2}{2}{2}{3}{3}{3}{white}{1};
			
			\draw[blue,dotted, ,very thick] (2, 3, 3) -- (2 ,2, 3);
			\draw[blue,dotted, very thick ] (2, 3, 3) -- (2 , 3, 2);
			\draw[blue,dotted, very thick] (3, 2, 3) -- (0, 2, 3);
			\draw[blue,dotted, very thick] (3, 2, 3) -- (3, 2, 2);
			\draw[blue,dotted, very thick] (3 ,3 ,2) -- (0, 3, 2);
			\draw[blue,dotted, very thick] (3 ,3 ,2) -- (3, 0, 2);
			
			\foreach \Point in {(2 ,2, 3), (2 , 3, 2), (0, 2, 3), (3, 2, 2), (0, 3, 2), (3, 0, 2)}{
				\node at \Point {\blue{\textbullet}};
			}
			\foreach \Point in {(2, 3, 3), (3, 2, 3), (3 ,3 ,2)}{
				\node at \Point {$\smallFilledSquare$};
			}
			\node[black] at (4, 2, 1) {1};
			\node[black] at (2.5, 1, 3) {2};
			\node[black] at (0.5, 2.4, 2.7) {3};
			\node[black] at (3, 2.5, 2.5) {4};
		\end{tikzpicture}
		\caption{Assuming that there are three already loaded items~1, 2, and 3 in the \ULD, the same loading position for item~4 results in the \blue{blue} extreme points~\blue{\textbullet}.}\label{fig:projection2}
	\end{subfigure}
	\caption{Generating extreme points by projection.}\label{fig:projection}
\end{figure}

Whenever an item is loaded, we are interested in positions which, for at least one dimension, cannot be shifted further towards the \ULD\ origin.
Otherwise, the remaining space would be fragmented unnecessarily.
To obtain these positions (the extreme points), for each newly inserted item,
\cite{CrainicEtAl2008} project certain item corner points along the orthogonal axis of the \ULD.
Here, projecting means to 'walk along a given axis' starting from a given point towards the coordinate origin until another item or a \ULD\ wall is hit as shown in Figure~\ref{fig:projection}.

For a given item $\itemIndex$ of size $\itemSize$ loaded at position $\extremePoint$, \cite{CrainicEtAl2008} suggest the following corner points as starting positions for the projections:
\begin{itemize}[noitemsep,label={--}]
\item $(\extremePoint_1, \extremePoint_2, \extremePoint_3 + \size_3)$ for projection in $x$- and $y$-direction,
\item $(\extremePoint_1, \extremePoint_2 + \size_2, \extremePoint_3)$ for projection in $x$- and $z$-direction, and
\item $(\extremePoint_1 + \size_1, \extremePoint_2, \extremePoint_3)$ for projection in $y$- and $z$-direction.
\end{itemize}

We slightly modify this approach to start the projection further away from the origin and, hence, potentially include more already loaded items in our projection routine (see Figure~\ref{fig:projection}):
\begin{itemize}[noitemsep,label={--}]
\item $(\extremePoint_1 + \size_1, \extremePoint_2 + \size_2, \extremePoint_3)$ for projection in $x$- and $y$-direction,
\item $(\extremePoint_1, \extremePoint_2 + \size_2, \extremePoint_3 + \size_3)$ for projection in $x$- and $z$-direction, and
\item $(\extremePoint_1+ \size_1, \extremePoint_2, \extremePoint_3 + \size_3)$ for projection in $y$- and $z$-direction.
\end{itemize}

\begin{figure}
	\begin{subfigure}[t]{0.45\textwidth}
		\centering
		\begin{tikzpicture}
			\draw[draw=black] (0,0) rectangle ++(5,2);
			\foreach  \x/\y/\w/\h/\n in {0/0/2.7/0.4/1,0/0.4/1/0.5/2,0/0.9/1.2/0.8/3,1/0.4/2/0.5/4}{
				\draw[draw=black] (\x,\y) rectangle ++(\w,\h);
				\node[] at (\x + \w/2,\y+\h/2) {\n};
			}
			\foreach \Point in {(0,1.7),(2.7,0)}{
				\node at \Point {\blue{$\times$}};
			}
			\foreach \Point in {(3, 0)}{
				\node at \Point {\blue{\textbullet}};
			}
			\node at (1,0.9) {$\smallFilledSquare$};
			\draw[blue,dotted, ,very thick] (3,0.4) -- (3, 0);
			\node at (3,0.4) {$\smallFilledSquare$};
		\end{tikzpicture}
		\caption{Extreme points according to \cite{CrainicEtAl2008}.}
	\end{subfigure}
	\hfill
	\begin{subfigure}[t]{0.45\textwidth}
		\centering
		\begin{tikzpicture}
			\draw[draw=black] (0,0) rectangle ++(5,2);
			\foreach  \x/\y/\w/\h/\n in {0/0/2.7/0.4/1,0/0.4/1/0.5/2,0/0.9/1.2/0.8/3,1/0.4/2/0.5/4}{
				\draw[draw=black] (\x,\y) rectangle ++(\w,\h);
				\node[] at (\x + \w/2,\y+\h/2) {\n};
			}
			\foreach \Point in {(0,1.7),(2.7,0)}{
				\node at \Point {\blue{$\times$}};
			}
			\foreach \Point in {(3, 0), (1.2,0.9)}{
				\node at \Point {\blue{\textbullet}};
			}
			\draw[blue,dotted, ,very thick] (3,0.9) -- (1.2,0.9);
			\draw[blue,dotted, ,very thick] (3,0.9) -- (3, 0);
			\node at (3,0.9) {$\smallFilledSquare$};
		\end{tikzpicture}
		\caption{Extreme points according to our procedure.}
	\end{subfigure}
	\caption{Items 1, 2, and 3 are loaded already. \blue{Blue} crosses~\blue{$\times$} denote the already existent, not blocked extreme points. New extreme points created by loading item 4 are highlighted by \blue{\textbullet} and the corresponding projection starting point by $\smallFilledSquare$.}\label{fig:projection_starting_point}
\end{figure}
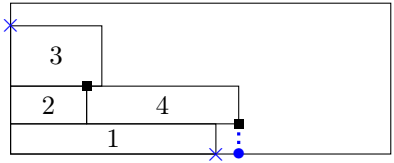
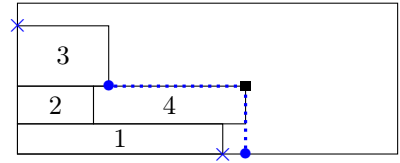

Figure~\ref{fig:projection_starting_point} shows a two-dimensional example for which this modification yields additional viable candidate points.

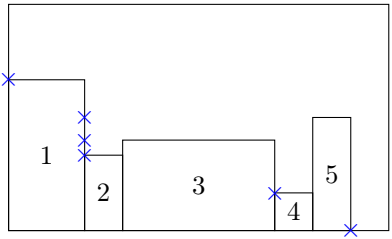
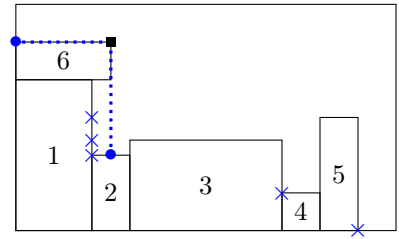
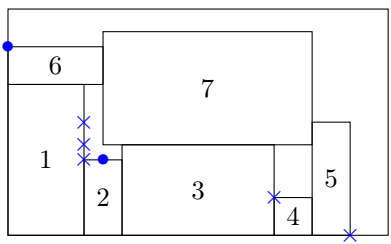
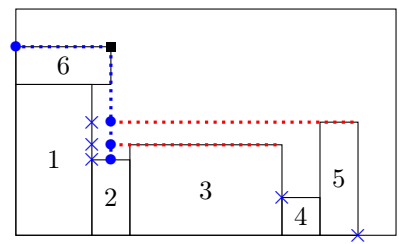
\begin{figure}[]
	\begin{subfigure}[t]{0.45\textwidth}
		\centering
		\begin{tikzpicture}
			\draw[draw=black] (0,0) rectangle ++(5,3);
			\foreach  \x/\y/\w/\h/\n in {0/0/1/2/1,1/0/0.5/1/2,1.5/0/2/1.2/3,3.5/0/0.5/0.5/4,4/0/0.5/1.5/5}{
				\draw[draw=black] (\x,\y) rectangle ++(\w,\h);
				\node[] at (\x + \w/2,\y+\h/2) {\n};
			}
			\foreach \Point in {(1,1), (1,1.2), (3.5,0.5), (4.5,0), (1,1.5), (0,2)}{
				\node at \Point {\blue{$\times$}};
			}
		\end{tikzpicture}
		\caption{Current status. Extreme points are denoted by \blue{$\times$}.}\label{fig:extended_projection1}
	\end{subfigure}
	\hfill
	\begin{subfigure}[t]{0.45\textwidth}
		\centering
		\begin{tikzpicture}
			\draw[draw=black] (0,0) rectangle ++(5,3);
			\foreach  \x/\y/\w/\h/\n in {0/0/1/2/1,1/0/0.5/1/2,1.5/0/2/1.2/3,3.5/0/0.5/0.5/4,4/0/0.5/1.5/5,0/2/1.25/0.5/6}{
				\draw[draw=black] (\x,\y) rectangle ++(\w,\h);
				\node[] at (\x + \w/2,\y+\h/2) {\n};
			}
			\draw[blue,dotted, ,very thick] (1.25,2.5) -- (0,2.5);
			\draw[blue,dotted, ,very thick] (1.25,2.5) -- (1.25,1);
			\node at (1.25,2.5) {$\smallFilledSquare$};
			\foreach \Point in {(1,1), (1,1.2), (3.5,0.5), (4.5,0), (1,1.5)}{
				\node at \Point {\blue{$\times$}};
			}
			\foreach \Point in {(1.25,1), (0,2.5)}{
				\node at \Point {\blue{\textbullet}};
			}
		\end{tikzpicture}
		\caption{Item~6 is loaded. When only considering directly hit items and the \ULD\ wall, only two new extreme points~\blue{\textbullet} are generated when projecting from~$\smallFilledSquare$.}\label{fig:extended_projection2}
	\end{subfigure}
	
	\vspace{3ex}
	\begin{subfigure}[t]{0.45\textwidth}
		\centering
		\begin{tikzpicture}
			\draw[draw=black] (0,0) rectangle ++(5,3);
			\foreach  \x/\y/\w/\h/\n in {0/0/1/2/1,1/0/0.5/1/2,1.5/0/2/1.2/3,3.5/0/0.5/0.5/4,4/0/0.5/1.5/5,0/2/1.25/0.5/6, 1.25/1.2/2.75/1.5/7}{
				\draw[draw=black] (\x,\y) rectangle ++(\w,\h);
				\node[] at (\x + \w/2,\y+\h/2) {\n};
			}
			\foreach \Point in {(1,1), (1,1.2), (3.5,0.5), (4.5,0), (1,1.5)}{
				\node at \Point {\blue{$\times$}};
			}
			\foreach \Point in {(1.25,1), (0,2.5)}{
				\node at \Point {\blue{\textbullet}};
			}
		\end{tikzpicture}
		\caption{Item 7 would fit well at the depicted position but cannot be loaded since no extreme point has been generated.}\label{fig:extended_projection3}
	\end{subfigure}
	\hfill
	\begin{subfigure}[t]{0.45\textwidth}
		\centering
		\begin{tikzpicture}
			\draw[draw=black] (0,0) rectangle ++(5,3);
			\foreach  \x/\y/\w/\h/\n in {0/0/1/2/1,1/0/0.5/1/2,1.5/0/2/1.2/3,3.5/0/0.5/0.5/4,4/0/0.5/1.5/5,0/2/1.25/0.5/6}{
				\draw[draw=black] (\x,\y) rectangle ++(\w,\h);
				\node[] at (\x + \w/2,\y+\h/2) {\n};
			}
			\draw[red,dotted, ,very thick] (1.25,1.2) -- (3.5,1.2);
			\draw[red,dotted, ,very thick] (1.25,1.5) -- (4.5,1.5);
			\draw[blue,dotted, ,very thick] (1.25,2.5) -- (0,2.5);
			\draw[blue,dotted, ,very thick] (1.25,2.5) -- (1.25,1);
			\node at (1.25,2.5) {$\smallFilledSquare$};
			\foreach \Point in {(1,1), (1,1.2), (3.5,0.5), (4.5,0), (1,1.5)}{
				\node at \Point {\blue{$\times$}};
			}
			\foreach \Point in {(1.25,1), (0,2.5), (1.25, 1.2), (1.25, 1.5)}{
				\node at \Point {\blue{\textbullet}};
			}
		\end{tikzpicture}
		\caption{Generate additional extreme points for positions where a newly loaded item might rest on already loaded items~3 or 5.}
		\label{fig:extended_projection4}
	\end{subfigure}
	\caption{Example of additional extreme points when loading a new item for a two-dimensional packing problem.}
	\label{fig:extended_projection}
\end{figure}

For an efficient packing, only creating the extreme points that are obtained from \emph{directly} hitting either the \ULD\ wall or another item is often insufficient:
Figures~\ref{fig:extended_projection1} and \ref{fig:extended_projection2} depict the simple, described way of generating extreme points (as two-dimensional example).
Now, if we try to load item 7, we can observe that it would fit best at the position illustrated in Figure~\ref{fig:extended_projection3}.
However, it cannot be loaded in this spot as no corresponding extreme point has been generated.
We resolve this issue by extending \cite{CrainicEtAl2008} procedure to also generate extreme points for \emph{all} positions at which the projection has not yet directly hit another item,
but if a sufficiently large item were to be loaded at this position, it would rest on another already loaded item.
This is done by artificially extending the item's surface towards the origin along the non-projection directions before performing the projection as illustrated in Figure~\ref{fig:extended_projection4}.

We now formalize the procedure to generate new extreme points.
Algorithm~\ref{alg:calculate_extreme_points} defines the six projections that must be performed and adds an additional extreme point on top a newly loaded stackable item.
Algorithm~\ref{alg:projection} is the projection subroutine that is called by Algorithm~\ref{alg:calculate_extreme_points} in line~\ref{alg:call_projection}.

\begin{algorithm}
	\DontPrintSemicolon
	\SetKw{Continue}{continue}
	\SetKwInput{Input}{Input}
	\SetKwInput{Output}{Output}
	\newcommand\mycommfont[1]{\footnotesize\ttfamily\textcolor{mygreen}{#1}}
	\SetCommentSty{mycommfont}
	\SetNoFillComment
	\Input{Loaded items $\loadedItems$ and an item $\itemIndex$ of size $\size$ loaded at position $\extremePoint$}
	\Output{Set $\extremePoints$ of new extreme points}
	Define set of new extreme points $\extremePoints=\emptyset$\;
	\tcc{Add extreme points via projection}
	Define set of directions $\dimensions=\{1,2,3\}$\;
	\For{$\iterator \in \dimensions$}{
		\If{$\itemIndex$ non-stackable and $\iterator=3$}{\label{alg:scip_for_non_stackable_start}
			\Continue\;
		}\label{alg:scip_for_non_stackable_end}
		\For{$\dimension \in \dimensions$ with $\iterator\neq \dimension$}{
			$\position=\extremePoint$\;\label{alg:starting_point1}
			$\position_\iterator = \position_\iterator + \size_\iterator$\;
			$\position_\dimension = \position_\dimension + \size_\dimension$\;\label{alg:starting_point2}
			$\extremePoints = \extremePoints \cup \text{Projection}(\position,\dimension,\loadedItems)$\;\label{alg:call_projection}
		}
	}
	\tcc{Add extreme point on top of item}
	\If{$\itemIndex$ stackable\label{alg:top_extreme_point_start}}{
		$\extremePoints = \extremePoints \cup \{(\extremePoint_1,\extremePoint_2,\extremePoint_3+\size_3)\}$\;\label{alg:top_extreme_point}
	}\label{alg:top_extreme_point_end}
	\Return $\extremePoints$\;
	\caption{Generate new extreme points.}\label{alg:calculate_extreme_points}
\end{algorithm}

\begin{figure}
	\begin{subfigure}[t]{0.45\textwidth}
		\centering
		\begin{tikzpicture}
			\draw[draw=black] (0,0) rectangle ++(5,2);
			\foreach  \x/\y/\w/\h/\n in {0/0/0.8/0.8/1,0.8/0/1.2/1/2,2/0/3/1.5/3}{
				\draw[draw=black] (\x,\y) rectangle ++(\w,\h);
				\node[] at (\x + \w/2,\y+\h/2) {\n};
			}
			\foreach \Point in {(0,0.8), (0,1), (0,1.5)}{
				\node at \Point {\blue{$\times$}};
			}
			\foreach \Point in {(0.8,1), (2,1.5)}{
				\node at \Point {\blue{$\star$}};
			}
		\end{tikzpicture}
		\caption{Example with three loaded items. Extreme points are denoted by~\blue{$\times$} and extreme points on top of items (added in line~\ref{alg:top_extreme_point}) are denoted by the star~\blue{$\star$}.}
	\end{subfigure}
	\hfill
	\begin{subfigure}[t]{0.45\textwidth}
		\centering
		\begin{tikzpicture}
			\draw[draw=black] (0,0) rectangle ++(5,2);
			\foreach  \x/\y/\w/\h/\n in {0/0/0.8/0.8/1,0.8/0/1.2/1/2,2/0/3/1.5/3,0.8/1/1.2/0.8/4}{
				\draw[draw=black] (\x,\y) rectangle ++(\w,\h);
				\node[] at (\x + \w/2,\y+\h/2) {\n};
			}
			\foreach \Point in {(0,0.8), (0,1), (0,1.5), (0,1.8)}{
				\node at \Point {\blue{$\times$}};
			}
			\foreach \Point in {(0.8,1.8), (2,1.5)}{
				\node at \Point {\blue{$\star$}};
			}
		\end{tikzpicture}
		\caption{Item~4 can only be loaded on top of item~2.}
	\end{subfigure}
	\caption{Example of a special case where the extreme point on top of the item is needed to load item~4 (see Algorithm~\ref{alg:calculate_extreme_points} line~\ref{alg:top_extreme_point}).}\label{fig:special_case}
\end{figure}
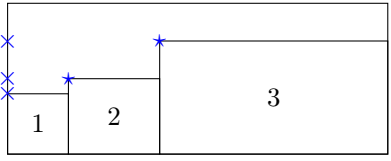
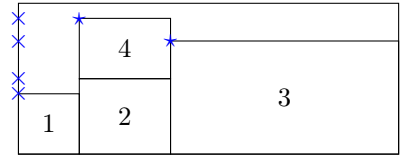

In lines~\ref{alg:starting_point1}-\ref{alg:starting_point2} of Algorithm~\ref{alg:calculate_extreme_points}, we compute the described starting points for our projections depending on the projection direction $\dimension \in \dimensions$.
If the newly inserted item is non-stackable, we do not want to generate any extreme points that would require items loaded at these points to directly rest on the non-stackable item at hand.
Therefore, lines~\ref{alg:scip_for_non_stackable_start}--\ref{alg:scip_for_non_stackable_end} ensure that we do not perform any projection in $x$- or $y$-direction starting from
a corner point on top of the item.
Next to the extreme points of the six projections, the additional extreme point~$(\extremePoint_1,\extremePoint_2,\extremePoint_3+\size_3)$ is added in line~\ref{alg:top_extreme_point} as it turned out to improve the solution quality.
Figure~\ref{fig:special_case} shows a two-dimensional example where this special extreme point is useful.

\begin{algorithm}
	\DontPrintSemicolon
	\SetKw{Continue}{continue}
	\SetKw{Break}{break}
	\SetKwInput{Input}{Input}
	\SetKwInput{Output}{Output}
	\newcommand\mycommfont[1]{\footnotesize\ttfamily\textcolor{mygreen}{#1}}
	\SetCommentSty{mycommfont}
	\SetNoFillComment
	\Input{Projection starting point $\position$, projection direction $\dimension$, and loaded items $\loadedItems$ with loading position~$\loadingPoint^\loadedItemIndex$ and size $\size^\loadedItemIndex$ for $\loadedItemIndex\in\loadedItems$}
	\Output{Set of new extreme points $\extremePoints$}
	Initialize newly generated extreme point $\extremePoint=\position$\;
	Define non-projection directions $\dimensionTwo$ and $\dimensionThree$ such that $\{\dimensionTwo,\dimensionThree,\dimension\} = \{1,2,3\}$\;
	Define set of blocking items $\blockingItemSet=\emptyset$\;
	Define set of extreme points $\extremePoints=\emptyset$\;
	Sort items $\loadedItemIndex \in \loadedItems$ in non-ascending end position order according to the projection direction, \ie, $\loadingPoint_\dimension^\loadedItemIndex + \size_\dimension^\loadedItemIndex$\;\label{alg:item_sorting}
	\For{$\loadedItemIndex \in \loadedItems$}{
		\tcc{Filter out items that are irrelevant for projection}
		\If{$\loadingPoint_\dimension^\loadedItemIndex \geq \position_\dimension$ or $\loadingPoint_\dimensionTwo^\loadedItemIndex + \size_\dimensionTwo^\loadedItemIndex \leq \position_\dimensionTwo$ or $\loadingPoint_\dimensionThree^\loadedItemIndex + \size_\dimensionThree^\loadedItemIndex \leq \position_\dimensionThree$}{\label{alg:line:irrelevant_start}
			\Continue\;
		}\label{alg:line:irrelevant_end}
		\tcc{Filter out items whose projection is blocked by others}
		\If{ ($\loadingPoint_\dimensionTwo^\loadedItemIndex \geq \loadingPoint_\dimensionTwo^\blockingItemIndex$ or $\position_\dimensionTwo \geq \loadingPoint_\dimensionTwo^\blockingItemIndex$) and ($\loadingPoint_\dimensionThree^\loadedItemIndex \geq \loadingPoint_\dimensionThree^\blockingItemIndex$  or $\position_\dimensionThree \geq \loadingPoint_\dimensionThree^\blockingItemIndex$) for any $\blockingItemIndex\in \blockingItemSet$}{\label{alg:line:blocked_items_start}
			\Continue\;
		}\label{alg:line:blocked_items_end}
		\tcc{Add extreme point}
		\If{$\loadingPoint_\dimension^\loadedItemIndex +  \size_\dimension^\loadedItemIndex \leq \position_\dimension$
			and ($\dimension \neq 3$ or $\loadedItemIndex$ is stackable)}{\label{alg:scip_for_non_stackable_projection}
			$\extremePoint_\dimension = \loadingPoint_\dimension^\loadedItemIndex + \size_\dimension^\loadedItemIndex$\;
			$\extremePoints=\extremePoints\cup \{\extremePoint\}$\;\label{alg:line:add_extreme_point}
		}\label{alg:line:end_if}
		\tcc{Stop if projection directly hits the item (all following items are blocked)}
		\If{$\position_\dimensionTwo \geq \loadingPoint_\dimensionTwo^\loadedItemIndex$ and $\position_\dimensionThree \geq \loadingPoint_\dimensionThree^\loadedItemIndex$ \label{alg:line:start_if2}}{
			\Return $\extremePoints$\;\label{alg:line:item_is_hit}
		}\label{alg:line:end_if2}
		\tcc{Add item to list of potentially blocking items}
		$\blockingItemSet=\blockingItemSet\cup \{\loadedItemIndex\}$\;\label{alg:line:add_blocking_items}
	}
	\tcc{Project to \ULD\ wall if projection is not blocked by previous items}
	$\extremePoint_\dimension =0$\;\label{alg:wall_projection}
	$\extremePoints=\extremePoints\cup \{\extremePoint\}$\;\label{alg:projection_to_wall}
	\Return $\extremePoints$\;
	\caption{Projection.}\label{alg:projection}
\end{algorithm}

\begin{figure}
	\centering
	\tdplotsetmaincoords{60}{115}
	\begin{tikzpicture}
		[tdplot_main_coords,
		cube/.style={black},
		axis/.style={->,black}]
		
		\draw[axis] (0,0,0) -- (4,0,0) node[anchor=west]{$x$};
		\draw[axis] (0,0,0) -- (0,5,0) node[anchor=west]{$y$};
		\draw[axis] (0,0,0) -- (0,0,4) node[anchor=west]{$z$};
		
		\mycuboid{0}{0}{0}{2}{3}{1}{white}{1};
		\mycuboid{2}{0}{0}{3}{1}{0.5}{white}{1};	
		\mycuboid{0}{0}{1}{1}{2.5}{2}{white}{1};
		\mycuboid{0}{0}{2}{1.5}{2.5}{3}{white}{1};
		\mycuboid{2}{1}{0}{2.5}{1.5}{1.2}{white}{1};
		\mycuboid{0}{2.5}{1}{3}{3}{3}{gray}{0.6};
		\mycuboid{0}{3}{0}{1}{4}{3}{white}{1};
		\mycuboid{0}{4}{0}{1.2}{4.5}{2.5}{white}{1};
		
		\draw[blue,dotted,very thick] (2.5,1,1.2) -- (1,1,1.2);
		\draw[red,dotted,very thick] (1.2,3,1.2) -- (1.2,4.5,1.2);
		\draw[red,dotted,very thick] (1.5,1,1.2) -- (1.5,1,3);
		
		\foreach \Point in {(1.5,1,1.2), (1,1,1.2)}{
			\node at \Point {\blue{\textbullet}};
		}
		
		\node at (2.5,1,1.2) {$\smallFilledSquare$};
		
		\node[black] at (1,0.25,1.25) {1};
		\node[black] at (2,0.25,0.75) {2};
		\node[black] at (3,0.25,0.25) {3};
		\node[black] at (1,3.25,2.75) {4};
		\node[black] at (1.2,4.25,2.25) {5};
		\node[black] at (1.5,0.25,2.75) {6};
		\node[black] at (3,2.75,2.75) {7};
		\node[black] at (2.5,1.25,0.6) {8};
	\end{tikzpicture}
	\caption{Example of calculating extreme points in $x$-direction when loading the item~8. The starting point of the projection is denoted by the black square~$\smallFilledSquare$. New extreme points are denoted by \blue{blue} bullets~\blue{\textbullet}. Item~7 is a blocking item.}\label{fig:extreme_points_3d}
\end{figure}
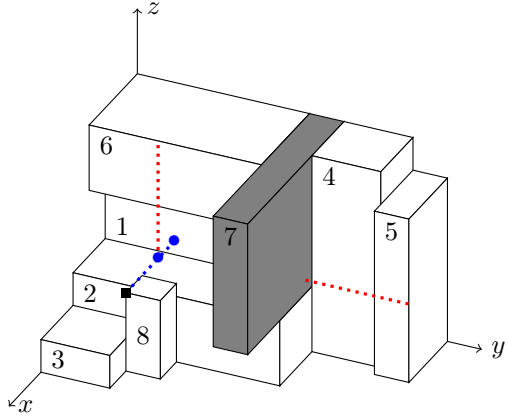

We explain the projection subroutine (Algorithm~\ref{alg:projection}) with the example depicted in Figure~\ref{fig:extreme_points_3d}.
Item~8 has just been loaded and we want to calculate new extreme points by projecting in $x$-direction.
We iterate over all already loaded items in non-ascending order of the items' end positions' $x$-coordinates (see line~\ref{alg:item_sorting}).
In lines~\ref{alg:line:irrelevant_start}--\ref{alg:line:irrelevant_end}, we then filter out items that cannot be hit by the projection (either directly or by surface extension).
These are all items~$\loadedItemIndex$ that are beyond the projection's starting point $\position$ in projection direction ($\loadingPoint_\dimension^\loadedItemIndex \geq \position_\dimension$) or are placed on the origin side of $\position$ in at least one of the two non-projection directions ($\loadingPoint_\dimensionTwo^\loadedItemIndex + \size_\dimensionTwo^\loadedItemIndex \leq \position_\dimensionTwo$ or $\loadingPoint_\dimensionThree^\loadedItemIndex + \size_\dimensionThree^\loadedItemIndex \leq \position_\dimensionThree$).
In the example, item~3 is filtered out as it's on the origin side of $\position$ in $y$-direction.

In the next step (lines~\ref{alg:line:blocked_items_start}--\ref{alg:line:blocked_items_end}), we filter out those items that, in theory, are relevant for the projection but whose surface extension is blocked by other items.
In our example, item 5 is blocked by item 7 and the surface extension is highlighted by the \red{red} dotted line.

Line~\ref{alg:scip_for_non_stackable_projection} ensures that no extreme points resulting from a projection in $z$-direction that -- directly or by extension -- hit a non-stackable item are added.
If an item is not filtered out by any of the above criteria, the corresponding extreme point is added in line~\ref{alg:line:add_extreme_point}.
In the example, extreme points are added via projection onto items 6 and 1.
We terminate the algorithm as soon as an item is hit directly since all other projections are blocked by this item (line~\ref{alg:line:item_is_hit}).
If no item is hit directly, an extreme point at the \ULD\ wall is added (line~\ref{alg:projection_to_wall}).
In the example, the algorithm stops in line~\ref{alg:line:item_is_hit} when item~1 is hit.
Note that the projection to the \ULD\ wall in line~\ref{alg:wall_projection} is only valid for non-tilted facets.
In case the projection hits a tilted facet, the corresponding coordinate can be calculated easily via the corresponding plane equation.

Since \cite{CrainicEtAl2008} only create extreme points that are obtained from directly hitting either the \ULD\ wall or another item, their approach generates a maximum of six extreme points for each newly loaded item.
In our algorithm, (significantly) more extreme points can be generated for each newly loaded item.
Nevertheless, our projection routine is fast.
With the sorting of $\loadedItems$ at the beginning and two nested for-loops, one running over the loaded items~$\loadedItems$ and the other one running over a subset of the loaded items~$\blockingItemSet\subseteq\loadedItems$, the overall time complexity is $\mathcal{O}(\numLoadedItems^2)$.
Note that the set of blocking items is typically small ($\numBlockingItems \ll\numLoadedItems$).

\section{Randomized Greedy Search}\label{sec:loading_ULD}

In order to improve the results of the insertion heuristic and to find a loading pattern for a given \ULD \ that simultaneously maximizes the volume of loaded items and considers the weight balance, we apply a \emph{Randomized Greedy Search} (\RGS).

The \RGS\ repeatedly calls the insertion heuristic with different item sorting criteria~$\sortingCriterion_1,\ldots,\sortingCriterion_5$ (see Section~\ref{sec:choose_item}).
The selection of the sorting criteria is evenly distributed over all $\numRGSiterations$ iterations, \ie, each sorting criterion is performed $\numRGSiterations/5$ times.
If the use of a substructure is allowed, all sorting criteria are applied twice -- once for the \ULD\ with and once for the \ULD\ without the substructure.
After each iteration, the quality of the newly generated solution is evaluated and the best solution is updated accordingly.
Finally, after all iterations have been performed, we improve the load stability of the best found solution by trying to close potential gaps between the loaded items.
The procedure is summarized in Algorithm~\ref{alg:randomize_greedy_search}.

\begin{algorithm}[H]
	\DontPrintSemicolon
	\SetKwInput{Input}{Input}
	\SetKwInput{Output}{Output}
	\Input{%
		Set of items $\itemSet$ and a \ULD~$\ULDIndex\in\ULDSet$, \linebreak
		number~$\numRGSiterations$ of RGS iterations, \linebreak
		different sorting criteria, \linebreak
		degree of randomization~$\degreeOfRandomization$}
	\Output{Loaded \ULD~$\ULDIndex$}
	Best found solution $S_{\text{best}}=\emptyset$\;
	\For{$\useSubstructure \in\{0,1\}$ (whether a substructure is used)}{\label{alg:use_substructure}
		\For{$\sortingCriterion \in\{\sortingCriterion_1,\ldots,\sortingCriterion_5\}$ (iterate over all sorting criteria)}{
			\For{$\iterator \in \{1, 2, ..., \numRGSiterations/5\}$}{
				Set degree of randomization $\selectedDegreeOfRandomization=\degreeOfRandomization$\;
				\If{$j=1$}{$\selectedDegreeOfRandomization=0$\;}
				New solution $\solution=\text{InsertionHeuristic}(\itemSet,\ULDIndex,\sortingCriterion,\selectedDegreeOfRandomization,\useSubstructure)$\;\label{alg:insertion_heuristic_call}
				\If{$\solution$ is better than $\solution_{\text{best}}$}{\label{alg:solution_selection}
					$\solution_{\text{best}} = \solution$\;
				}
			}
		}
	}
	Refine the solution $\solution_{\text{best}}$ by avoiding holes\;\label{alg:avoid_holes}
	\Return $\solution_{\text{best}}$\;
	\caption{Randomized Greedy Search.}\label{alg:randomize_greedy_search}
\end{algorithm}

For simplicity, we assume that the use of a substructure is permitted, \ie, $\substructureAllowed = 1$ and $\useSubstructure \in \{0,1\}$.
Otherwise, we can simply omit the for loop in line~\ref{alg:use_substructure} ($\useSubstructure = 0$).
When calling the insertion heuristic in line~\ref{alg:insertion_heuristic_call}, we ensure that for the first iteration of each sorting policy, we do not apply any randomization.
The user-defined degree of randomization will only be considered in the later iterations.

In the following, we explain the solution selection and the avoidance of holes in the load plan.

\subsection{Solution selection}\label{sec:solution_selection}

The comparison between different solutions is based on a score value that takes the volume utilization, weight balance, and whether the items would also fit into the remaining \ULD s into account.
Besides the maximum center of gravity deviation~$\maximumCenterOfGravityDeviation\in[0,1]$ that defines in which area the \centerOfGravity\ must be (see Figure~\ref{fig:weight_balance_3d}), we introduce the parameter $\weightBalanceImportance\in[0,1]$ which specifies how important the weight balance is during solution selection.
The volume utilization importance is defined by $\volumeUtilizationImportance = 1 - \weightBalanceImportance\in[0,1]$.

For dimension $\dimension\in \{1,2\}$ (only the horizontal weight balance is considered), we define the center of gravity deviation
\begin{align}
	\centerOfGravityDeviation_\dimension =
	\begin{cases}
		\left|\frac{\text{\centerOfGravity}_\dimension - \boundingBox_\dimension / 2}{\boundingBox_\dimension / 2}\right|  &\quad  \text{if }  \left|\frac{\text{\centerOfGravity}_\dimension - \boundingBox_\dimension / 2}{\boundingBox_\dimension / 2}\right| > \maximumCenterOfGravityDeviation \text{,}\\
		0 &\quad \text{otherwise,}
	\end{cases}
\end{align}
where $\boundingBox$ denotes the bounding box of the \ULD.
Then, the weight balance score is 
\begin{align}\label{eq:weight_balance_score}
	\weightBalanceScore = 1 - \frac{\centerOfGravityDeviation_1 + \centerOfGravityDeviation_2}{2}
\end{align}
and the volume utilization score is
\begin{align}
	\volumeUtilizationScore = \frac{\sum_{\itemIndex\in \loadedItems} \itemSize^\itemIndex_1 \itemSize^\itemIndex_2 \itemSize^\itemIndex_3}{\ULDVolume}.
\end{align}
In order to estimate how difficult it will be to load the remaining items in the remaining \ULD s, we define the set~$\remainingItems$ of items that could not be loaded yet.
Moreover, let $\ULDgroups$ be the set of ULD groups (identical \ULD s constitute one \ULDgroup) and $\ULDgroups_{\itemIndex}$ the set of ULD groups into which item~$\itemIndex\in\itemSet$ fits.
The score of a solution is then defined by
\begin{align}
	\solutionScore = \weightBalanceImportance \weightBalanceScore + \volumeUtilizationImportance \volumeUtilizationScore - \frac{\sum_{\itemIndex\in \remainingItems} (|\ULDgroups| - |\ULDgroups_{\itemIndex}|) v_\itemIndex}{\sum_{\itemIndex\in \itemSet} (|\ULDgroups| - |\ULDgroups_{\itemIndex}|) v_\itemIndex}.
\end{align}
Solutions with higher scores are preferred.

\subsection{Avoiding holes}\label{sec:avoiding_holes}

To prevent items from sliding, we try to remove horizontal holes between the loaded items by moving (sets of) items to the center of the \ULD\ in a post-processing step.
For algorithmic performance reasons, the algorithm described below is only performed at the end of the \RGS.
As a result, the weight balance score~\eqref{eq:weight_balance_score} of a solution may improve or even worsen, but this is not taken into account when selecting a solution within the \RGS.

The idea of the algorithm is to perform the two steps 
\begin{enumerate}\setlength\itemsep{-0.5ex}
	\item determine a set of items to be moved together\label{avoid_holes1}
	\item move the items as far as possible towards the center of the ULD\label{avoid_holes2}
\end{enumerate}
repeatedly until no items are moved or a maximum number of iterations is reached.

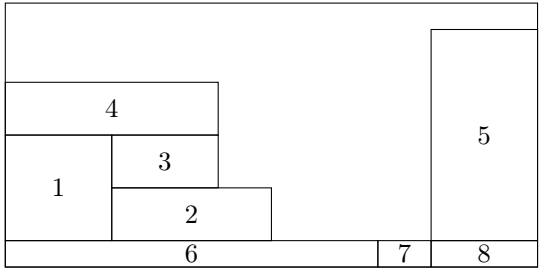
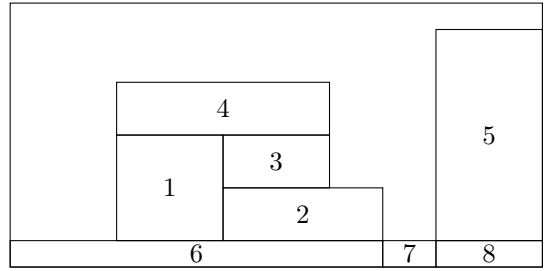
\begin{figure}[htbp]
	\begin{subfigure}[t]{0.45\textwidth}
	\centering
	\begin{tikzpicture}[scale=0.7]
		\draw[draw=black] (0,0) rectangle ++(10,5);
		\foreach  \x/\y/\w/\h/\n in {0/0.5/2/2/1,2/0.5/3/1/2,2/1.5/2/1/3,0/2.5/4/1/4,8/0.5/2/4/5,0/0/7/0.5/6,7/0/1/0.5/7, 8/0/2/0.5/8}{
			\draw[draw=black] (\x,\y) rectangle ++(\w,\h);
			\node[] at (\x + \w/2,\y+\h/2) {\n};
		}
	\end{tikzpicture}
	\caption{The set of items $\neighboringItems=\{1,2,3,4\}$ is moved together to avoid holes. When determining $\neighboringItems$, start items~2 and 3 both would yield the identical set $\neighboringItems$ without lines~\ref{alg:moving_set:return_empty_set_start}-\ref{alg:moving_set:return_empty_set_end} in Algorithm~\ref{alg:moving_set}.}
	\label{fig:set_of_movable_items}
	\end{subfigure}
	\hfill
	\begin{subfigure}[t]{0.45\textwidth}
	\centering
	\begin{tikzpicture}[scale=0.7]
		\draw[draw=black] (0,0) rectangle ++(10,5);
		\foreach  \x/\y/\w/\h/\n in {2/0.5/2/2/1,4/0.5/3/1/2,4/1.5/2/1/3,2/2.5/4/1/4,8/0.5/2/4/5,0/0/7/0.5/6,7/0/1/0.5/7, 8/0/2/0.5/8}{
			\draw[draw=black] (\x,\y) rectangle ++(\w,\h);
			\node[] at (\x + \w/2,\y+\h/2) {\n};
		}
	\end{tikzpicture}
	\caption{New position of the the item set $\neighboringItems$ determined by the binary search. Item 7 is non-stackable. Only considering collisions, $\neighboringItems$ could be moved all the way to item 5. As this would violate non-stackability, the items are only moved to the edge of item 7.}
	\label{fig:move_position}
	\end{subfigure}
	\caption{Example of moving a set of items.}
\end{figure}

In step~\ref{avoid_holes1}, we iterate over all loaded items~$\loadedItemIndex\in\loadedItems$ for which there is a hole in direction $x$ or $y$.
A hole is defined as empty space next to item~$\loadedItemIndex$ in direction of the ULD center with at least one (blocking) item between the item and the \ULD\ wall in this direction.
Afterwards, we determine a set~$\neighboringItems$ of neighboring items that should be moved together with $\loadedItemIndex$.
Note that as illustrated in Figure~\ref{fig:set_of_movable_items}, the set of neighboring items does not have to directly rest on the floor.
Starting with item~$\loadedItemIndex$, the set of neighboring items is determined by Algorithm~\ref{alg:moving_set}.
In line~\ref{alg:bounding_box}, the bounding box~$\boundingBox$ of $\neighboringItems$ is iteratively extended and subsequently, items are added to $\neighboringItems$ that either intersect with $\boundingBox$ or are supported (directly or indirectly) by $\boundingBox$.
The algorithm returns an empty set in line~\ref{alg:moving_set:return_empty_set} to avoid generating several identical sets~$\neighboringItems$ (see Figure~\ref{fig:set_of_movable_items}).

As soon as a non-empty set~$\neighboringItems$ is found, step~\ref{avoid_holes2} is performed.
Only considering collisions with other items, it would be straight forward to simply move the item set as far as required to close the hole.
However, due to non-floating and stackability restrictions, this would frequently result in infeasible load plans.
In these cases, we still aim to move the items as far as possible towards the center of the \ULD\ to reduce the size of the holes and increase load stability.
Hence, we perform a binary search on the interval defined by the farthest position to which $\boundingBox$ can be moved only considering collisions and the current position of the item set bounding box $\boundingBox$ to determine a feasible new position for the item set that minimizes the hole (see Figure~\ref{fig:move_position}).

Note that we only move items in one direction at a time.
In case there are holes for an item in both $x$- and $y$-direction, we first perform the above procedure for one direction and then check, after moving the corresponding items, if there are still holes in the other direction.

\begin{algorithm}
	\DontPrintSemicolon
	\SetKwInput{Input}{Input}
	\SetKwInput{Output}{Output}
	\Input{%
		Set of all loaded items $\loadedItems$ and an item~$\loadedItemIndex\in\loadedItems$ to be moved}
	\Output{A set of movable items $\neighboringItems\subset \loadedItems$}
	Initialize $\neighboringItems=\{ \loadedItemIndex\}$\;
	Define $\position_3$ as the $z$-position at which item $\loadedItemIndex$ is loaded\;
	\Repeat{No item is added to $\neighboringItems$}{
		Determine cuboid bounding box~$\boundingBox$ of $\neighboringItems$\;\label{alg:bounding_box}
		\If{An item~$\loadedItemIndexTwo \in \loadedItems \setminus \neighboringItems$ exists that either intersects with the bounding box~$\boundingBox$ or is supported by $\boundingBox$}{
			\If{$\extremePoint_3 < \position_3$ where $\extremePoint_3$ is the $z$-position at which item $\loadedItemIndexTwo$ is loaded\label{alg:moving_set:return_empty_set_start}}{
				\Return{$\neighboringItems=\emptyset$}\label{alg:moving_set:return_empty_set}
			}\label{alg:moving_set:return_empty_set_end}
			$\neighboringItems = \neighboringItems \cup \{\loadedItemIndexTwo\}$\;
		}
	}
	\Return $\neighboringItems$\;
	\caption{Determine set of movable items.}\label{alg:moving_set}
\end{algorithm}

\section{Loading multiple \ULD s}
\label{sec:heuristic}
To load items into multiple \ULD s, we choose a simple sequential approach (see Algorithm~\ref{alg:heuristic_overview}).
We introduce a \ULD\ selection routine to determine, based on the available items, which \ULD\ should be loaded next/first. 
Once no more items fit into the currently selected \ULD, the remaining items are tried to be loaded into the next \ULD.
As our \ULD\ selection approach generally prefers larger \ULD s over smaller ones,  at the end of the algorithm, the items of the last loaded \ULD\ are tried to be reloaded into a smaller available \ULD\ (line \ref{alg:reload_items}) to avoid unused space.

\begin{algorithm}
	\DontPrintSemicolon
	\SetKwInput{Input}{Input}
	\SetKwInput{Output}{Output}
	\SetNoFillComment
	\Input{Set of items to be loaded $\itemSet$ and set of available ULDs~$\ULDSet$}
	\Output{Feasible solution of \problem}
	\While{$\itemSet \neq \emptyset$ and $\ULDSet \neq \emptyset$}{\label{alg:begin_sequential_loading}
		Select next \ULD~$\ULDIndex\in\ULDSet$\;\label{alg:select_next_ULD}
		Load \ULD~$\ULDIndex$ with available items $\itemSet$\;\label{alg:load_ULD}
		Update $\itemSet$ and $\ULDSet$\;
	}\label{alg:end_sequential_loading}
	\If{$\ULDSet \neq \emptyset$}{
		Reload items of the last loaded \ULD\ into the smallest available \ULD\ in which all these items fit\;\label{alg:reload_items}
	}
	\caption{High-level algorithm overview}\label{alg:heuristic_overview}
\end{algorithm}

We will now explain the selection of the next \ULD\ from line~\ref{alg:select_next_ULD}.
Depending on the instance at hand, certain items are harder to load than others in the sense that they only fit into a limited subset of the available \ULD s.
As we aim to load all items, this means in turn that some of the \ULD s are required to be included in the solution while an efficient packing can render others obsolete.
Now, the idea is to first load the \ULD s that are required anyway and ensure that we utilize them as good as possible:

Let $\ULDgroups$ be the set of all \ULDgroup s (identical \ULD s constitute one \ULDgroup) and $\ULDgroups_{\itemIndex} \subseteq \ULDgroups$ the set of \ULDgroup s into which item~$\itemIndex\in\itemSet$ fits.
Then, the smallest number of fitting \ULDgroup s per item is
\begin{align}
	\smallesNumberFittingULDgroups=\min_{\itemIndex\in\itemSet} \left| \ULDgroups_{\itemIndex} \right|.
\end{align}
We define the set of potential next \ULDgroup s as
\begin{align}
	\mathcal{\ULDgroups}=\bigcup_{\itemIndex\in\itemSet : \left| \ULDgroups_{\itemIndex} \right| =\smallesNumberFittingULDgroups} \ULDgroups_{\itemIndex}.
\end{align}
We select the \ULDgroup\ $\ULDgroupsIndex \in \mathcal{\ULDgroups}$ with the highest cumulative volume of items that can (individually) be loaded in one of the \ULD s belonging to the \ULDgroup.
Ties are broken by the highest \ULD\ volume.

\section{Computational results}
\label{sec:comp_results}

In this section, we first give an overview of the benchmark instances and then describe details of the implementation and the computational setup.
Next, algorithm components are evaluated.
The section closes with a comparison to approaches from the literature and detailed results for an instance set adapted to our presented problem.

\subsection{Benchmark instances}

We consider three types of benchmark sets: three-dimensional Bin Packing Problem instances proposed by \cite{BischoffRatcliff1995}, three-dimensional \problem\ with transportation constraints instances proposed by \cite{PaquayEtAl2018}, and an adaption of the \cite{PaquayEtAl2018} instance set.

The first instance set by \cite{BischoffRatcliff1995} consists of seven subsets of instances grouped by the number of types of different items ($3, 5, 8, 10, 12, 15,  20$), each one containing 100 instances.
For all instances, one type of cuboid \ULD\ (or bin) is available.
To define feasible item orientations, it is specified for each dimension whether the item can be tilted across the corresponding axis.
This yields a more detailed item specification compared to ours, which only allows to define on a general level whether an item can be tilted.
Item weights and \ULD\ weight capacities are not considered.

The second benchmark set proposed by \cite{PaquayEtAl2018} consists of ten subsets of instances grouped by the number of items ($10,20,\ldots,100$), each containing 30 instances.
For all instances, an unlimited number of \ULD s of six types, which can be loaded in a Boeing 777, is available.
The item sets are based on real-world data.
However, the authors had to manipulate orientations and stackablity characteristics because these parameters were missing.
Note that we only consider the final data set and not the training data set.

The \ULD\ definition of \cite{PaquayEtAl2018} differs slightly from ours.
In particular, they do not consider the length/width of the edge on which loading at the base area is prohibited.
As a consequence, loading a substructure is also not reasonable.
Therefore, we additionally test our algorithm on an adapted instance set of the \cite{PaquayEtAl2018} instances.
The item sets remain unchanged and for each \ULD\, we define the length/width of the edge~$\edgeWidth=10$ and the vertical edge offsest~$\verticalEdgeOffset=10$.
Moreover, we not only consider instances with an unlimited number of available \ULD s, but also instances where each \ULD\ type is available only once.

\subsection{Details of the implementation and computational setup}

The \RGS\ is implemented in \texttt{C++} and compiled into a single-thread code.
The computational study is carried out on a MacBook Pro with an Apple M2 Pro chip and 32 GB of RAM.
%

We have conducted extensive preliminary tests to tune our algorithm and determine good parameter settings.
Particularly, we had to pay attention to ensuring that the solution quality was good and that the runtime was short at the same time.
Pretests were conducted on randomly generated and real-world instances.
We were in close contact with our business colleagues to assess the solution quality.

The parameter settings are summarized in Table~\ref{tab:comp_setup}.
We distinguish between algorithm-specific and problem-specific parameters.
Instead of an overall time limit, the number of extreme point checks (see line~\ref{alg:load_item} in Algorithm~\ref{alg:insertion_heuristic}) is limited,
making the stop criterion deterministic.
As most of the runtime can be attributed to these checks, it scales proportionally to the number of extreme point checks.
If a huge number of items fit into a \ULD, the maximum number of extreme point checks can be reached within the first \RGS\ iterations.
In this case, we nevertheless execute at least 10 \RGS\ iterations.
A maximum number of \RGS\ iterations is additionally set to ensure that small instances are solved in short time.
The calculation time for instances with up to 100 items is typically no more than 1 second.

\begin{table}[htbp]
	\centering
	\caption{Overview of parameter settings.}
	\adjustbox{width=0.7\textwidth}{
	\begin{tabular}{lll}
		\toprule
		Type &Parameter & Value\\
		\midrule
		algorithm-specific&maximum number of extreme point checks& 20{,}000{,}000\\
		&minimum number~$\numRGSiterations^{\min}$ of \RGS\ iterations & 10\\
		&maximum number~$\numRGSiterations^{\max}$ of \RGS\ iterations &500\\
		&degree of randomization~$\degreeOfRandomization$& 0.5\\
		&ascending lexicographic extreme point sorting& $z,y,x$\\
		problem-specific&minimum item overlap~$\minimumItemOverlap$& 0.9\\
		&maximum padding height~$\paddingHeight$& 10\\
		&maximum center of gravity deviation~$\maximumCenterOfGravityDeviation$& 0.1\\
		&weight balance importance~$\weightBalanceImportance$& 0.5\\
		\bottomrule
	\end{tabular}%
	}
	\label{tab:comp_setup}%
\end{table}%

\subsection{Evaluation of algorithm components}

In this section, we analyze individual algorithm components.
In particular, we evaluate the proposed grid acceleration, variations of our extreme point generation routine, and different sorting criteria.

The projection in Algorithm~\ref{alg:projection} differs from that in \citep{CrainicEtAl2008} as it not only relies on directly hit items but also items that are hit indirectly via surface extension, resulting in potentially more than one extreme point per projection.
We also consider blocking items and different projection starting points.
Moreover, we allow to move extreme points (see Section~\ref{sec:choose_extreme_point}).
To evaluate these changes and our grid acceleration, we tested our algorithm without the respective extensions on adapted \cite{PaquayEtAl2018} instances, where each \ULD\ type is just available once (cf. \emph{1 \ULD} in Section~\ref{sec:comp_results_paquay}).
In detail, we have changed the algorithm as follows:
\begin{itemize}
\item \emph{No grid acceleration} (\noGrid): We do not use the grid presented in Section~\ref{sec:acceleration_techniques} to accelerate the collision, non-floating, and stackability check when loading the next item.
\item \emph{No set of blocking items} (\noBlocking): In the projection routine, we filter out items that cannot be hit by the projection since they are blocked by others.
We omit this feature by discarding lines \ref{alg:line:blocked_items_start}--\ref{alg:line:blocked_items_end} and \ref{alg:line:add_blocking_items} of Algorithm~\ref{alg:projection}.
\item \emph{No moving of extreme points} (\noExtremePointShift): The moving of extreme points as described in Section~\ref{sec:choose_extreme_point} is not performed.
\item \emph{Mimic extreme point generation of \cite{CrainicEtAl2008}} (\oneExtremePoint): Only one extreme point per iteration is generated by projecting either on an already loaded item or the container wall.
In this variant, lines \ref{alg:line:blocked_items_start}--\ref{alg:line:blocked_items_end} and \ref{alg:line:add_blocking_items} of Algorithm~\ref{alg:projection} are discarded.
Moreover, the \emph{if}-statements in lines~\ref{alg:scip_for_non_stackable_projection}--\ref{alg:line:end_if2} are merged by removing lines \ref{alg:line:end_if}--\ref{alg:line:start_if2} and replacing the \emph{if}-condition in line \ref{alg:scip_for_non_stackable_projection} by \emph{$\position_\dimensionTwo \geq \loadingPoint_\dimensionTwo^\loadedItemIndex$ and $\position_\dimensionThree \geq \loadingPoint_\dimensionThree^\loadedItemIndex$ and $\loadingPoint_\dimension^\loadedItemIndex +  \size_\dimension^\loadedItemIndex \leq \position_\dimension$ and ($\dimension \neq 3$ or $\loadedItemIndex$ is stackable)}.
In addition, the projection start point is adapted and no additional extreme point on top of the newly loaded item is generated by discarding lines \ref{alg:starting_point2} and \ref{alg:top_extreme_point_start}--\ref{alg:top_extreme_point_end} of Algorithm~\ref{alg:calculate_extreme_points}.
The moving of extreme points as described in Section~\ref{sec:choose_extreme_point} is also not performed.
This approach corresponds to the one of \cite{CrainicEtAl2008}.
\end{itemize}

Results are summarized in Table~\ref{tab:analysis_extreme_points}.
For each instance set and algorithm variant, we report the average ratio of the variant runtime to the default runtime~$\solutionTime_{D}$, \ie, our algorithm without any changes.
Similarly, the ratios of the utilization are denoted.
In particular, all values \timeQuotientNoGrid, \timeQuotientNoBlocking, \timeQuotientOneExtremePoint, \timeQuotientNoExtremePointShift > 1 indicate that the default approach is faster.
Moreover, all values \utilizationQuotientNoGrid, \utilizationQuotientNoBlocking, \utilizationQuotientOneExtremePoint, \utilizationQuotientNoExtremePointShift < 1 indicate that the default approach has higher average utilization.

Using a grid accelerates the algorithm by an average of 18.5\%.
The grid is particularly advantageous for instances with many items.
For the variant \noBlocking, the solution time is significantly higher while the improvement in utilization by 0.1\% is negligible.
This shows that the use of the set of blocking items, in fact, prohibits the generation of irrelevant extreme points and allows for lower runtimes without diminishing solution quality.
As expected, the other two variants \noExtremePointShift\ and \oneExtremePoint\ are faster than the default approach since the algorithm has been simplified.
However, the utilization is on average 0.7\% and 2.2\% worse for the variant \noExtremePointShift\ and \oneExtremePoint, respectively.

\begin{table}[htbp]
	\centering
	\caption{Results for different definitions of the set of extreme points and omitted grid acceleration. The solution times $\solutionTime_{NG}$, $\solutionTime_{NB}$, $\solutionTime_{NM}$, $\solutionTime_{CR}$ and utilizations $\utilization_{NG}$, $\utilization_{NB}$, $\utilization_{NM}$, $\utilization_{CR}$ correspond to the variants \noGrid, \noBlocking, \noExtremePointShift, and \oneExtremePoint, respectively.}
	\adjustbox{width=0.6\textwidth}{
	\begin{tabular}{lrrrrrrrr}
		\toprule
		& \multicolumn{2}{c}{\noGrid} & \multicolumn{2}{c}{\noBlocking} & \multicolumn{2}{c}{\noExtremePointShift} & \multicolumn{2}{c}{\oneExtremePoint} \\
		\midrule
		$\numItems$ & \multicolumn{1}{l}{\timeQuotientNoGrid} & \multicolumn{1}{l}{\utilizationQuotientNoGrid} & \multicolumn{1}{l}{\timeQuotientNoBlocking} & \multicolumn{1}{l}{\utilizationQuotientNoBlocking} & \multicolumn{1}{l}{\timeQuotientNoExtremePointShift} & \multicolumn{1}{l}{\utilizationQuotientNoExtremePointShift} & \multicolumn{1}{l}{\timeQuotientOneExtremePoint} & \multicolumn{1}{l}{\utilizationQuotientOneExtremePoint} \\
		\midrule
		\multicolumn{1}{r}{10} & 0.947 & 1.000 & 0.970 & 1.000 & 0.999 & 1.000 & 0.920 & 1.008 \\
		\multicolumn{1}{r}{20} & 1.004 & 1.000 & 1.003 & 1.000 & 0.956 & 0.995 & 0.872 & 0.982 \\
		\multicolumn{1}{r}{30} & 1.068 & 1.000 & 1.043 & 1.012 & 0.940 & 0.997 & 0.777 & 0.956 \\
		\multicolumn{1}{r}{40} & 1.095 & 1.000 & 1.081 & 0.998 & 0.939 & 0.990 & 0.788 & 0.980 \\
		\multicolumn{1}{r}{50} & 1.151 & 1.000 & 1.109 & 0.999 & 0.921 & 0.983 & 0.797 & 0.970 \\
		\multicolumn{1}{r}{60} & 1.203 & 1.000 & 1.139 & 0.998 & 0.898 & 1.001 & 0.775 & 0.956 \\
		\multicolumn{1}{r}{70} & 1.257 & 1.000 & 1.159 & 1.003 & 0.912 & 0.991 & 0.771 & 1.002 \\
		\multicolumn{1}{r}{80} & 1.319 & 1.000 & 1.196 & 1.010 & 0.913 & 0.985 & 0.766 & 0.957 \\
		\multicolumn{1}{r}{90} & 1.377 & 1.000 & 1.228 & 0.999 & 0.919 & 0.986 & 0.757 & 0.986 \\
		\multicolumn{1}{r}{100} & 1.430 & 1.000 & 1.245 & 0.993 & 0.923 & 1.000 & 0.756 & 0.979 \\
		\midrule
		\total & 1.185 & 1.000 & 1.117 & 1.001 & 0.932 & 0.993 & 0.798 & 0.978 \\
		\bottomrule
	\end{tabular}%
	}
	\label{tab:analysis_extreme_points}%
\end{table}%

\begin{sloppypar}
We also evaluated how often which sorting criterion leads to the best solution for a \ULD.
Table~\ref{tab:results_sorting_policies} summarizes results for the adapted \cite{PaquayEtAl2018} benchmark groups \emph{unlimited} and \emph{1 \ULD} (see Section~\ref{sec:comp_results_paquay}).
The sorting policy \sorting{stackability--cumulated volume} leads to the best solution for most of the \ULD s.
Randomly drawing the next item is only useful for 37 \ULD s in total.
However, we point out that in line~\ref{alg:solution_selection} of Algorithm~\ref{alg:randomize_greedy_search}, the best solution is only updated if the new solution is better (and not equal to) the current best solution.
Since we apply the different sorting criteria in the order shown in Table~\ref{tab:results_sorting_policies}, a bias towards the earlier called policies is to be expected.
It is therefore not surprising that \sorting{stackability--cumulative volume} is considered the best for most \ULD s.
In turn, \sorting{random} is better than all other sorting criteria for 37 \ULD s.
Note that we have also tested other common sorting criteria such as the base area.
However, these other sorting criteria did not improve the solution quality.
\end{sloppypar}

\begin{table}[htbp]
	\centering
	\caption{Evaluation of sorting policies.}
	\adjustbox{width=0.5\textwidth}{
	\begin{tabular}{lrr}
		\toprule
		& \multicolumn{1}{l}{unlimited} & \multicolumn{1}{l}{1 \ULD} \\
		\midrule
		\sorting{stackability--cumulated volume} & 347   & 367 \\
		\sorting{stackability--highest volume} & 99    & 123 \\
		\sorting{cumulated volume} & 100   & 99 \\
		\sorting{highest volume} & 93    & 116 \\
		\sorting{random} & 14    & 23 \\
		\bottomrule
	\end{tabular}%
	}
	\label{tab:results_sorting_policies}%
\end{table}%

\subsection{Results for \cite{BischoffRatcliff1995} instances}

In order to study the general packing density of our approach, we compare the \RGS\ with two methods of \cite{BischoffRatcliff1995} for the classical three-dimensional bin packing problem.
Both methods focus on maximizing the volume utilization when loading a single bin.
The first method additionally takes loading stability into account by loading in a kind of layer structure.
The second method addresses multi-drop situations by building stacks.
Their approach is optimized to handle weakly heterogeneous items.

To ensure a fair comparison between the different approaches, we adapt the problem-specific parameters to minimum item overlap~$\minimumItemOverlap=1$, padding height~$\paddingHeight=0$, and weight balance importance~$\weightBalanceImportance=0$.
However, we do not handle tiltable items in the same way.
For items that can only be tilted across one axis, we define that the item is not tiltable at all.
Hence, the solution space considered in our approach is smaller.
In addition, the set of extreme points is sorted in ascending lexicographic order by $x,y,z$ to mimic the multi-drop situation by predominantly building stacks.
Note also that the edge width~$\edgeWidth=0$ and no substructure is allowed ($\substructureAllowed=0$).

Results grouped by the number~\numItemTypes\ of different types of items are summarized in Table~\ref{tab:results_bischoff_ratliff}.
For a single \ULD, we compare the average, minimum, and maximum utilization in percent, denoted by \avgUtilization, \minUtilization, and \maxUtilization, respectively.
The solution time~$\solutionTime$ is reported in milliseconds.
The average number of loaded items is denoted by \avgNumLoadedItems.
The \RGS\ algorithm outperforms the two methods of \cite{BischoffRatcliff1995}.
For the \RGS, the average utilization is $85.0\%$ and decreases when the number~\numItemTypes\ of item types increases.
This is different for the two methods of \cite{BischoffRatcliff1995}.
The average utilization rate of the \RGS\ is higher in all but one case.

Our results lag behind the tree-search-based approach of \cite{ArayaEtAl2017}.
Next to the smaller solution space, the main reason for this performance issues is most likely that our approach is optimized for instances with a more heterogeneous set of items and instances with fewer items because 90 or more items per \ULD\ are very rare in practice.
Note also that the solution time is below 1~second for all but two instances, which is significantly lower than the 30~seconds of \cite{ArayaEtAl2017}.
For these two instances, more than 350 items are loaded into one \ULD.

\begin{table}[htbp]
	\centering
	\caption{Comparison of \cite{BischoffRatcliff1995} (BR) and the proposed \RGS.}
	\begin{tabular}{rrrrrrrrrrrrrr}
		\toprule
		&       & \multicolumn{3}{c}{Method 1 of BR} & \multicolumn{3}{c}{Method 2 of BR} & \multicolumn{6}{c}{\RGS} \\
		\cmidrule(lr){3-5} \cmidrule(lr){6-8} \cmidrule(lr){9-14}
		\multicolumn{1}{l}{\numItemTypes} & \multicolumn{1}{l}{\avgNumItems} & \multicolumn{1}{l}{\avgUtilization} & \multicolumn{1}{l}{\minUtilization} & \multicolumn{1}{l}{\maxUtilization} & \multicolumn{1}{l}{\avgUtilization} & \multicolumn{1}{l}{\minUtilization} & \multicolumn{1}{l}{\maxUtilization} & \multicolumn{1}{l}{\avgUtilization} & \multicolumn{1}{l}{\minUtilization} & \multicolumn{1}{l}{\maxUtilization} & \multicolumn{1}{l}{\avgTime} & \multicolumn{1}{l}{\maxTime} & \multicolumn{1}{l}{\avgNumLoadedItems} \\
		\midrule
		3     & 150.4 & 81.8  & 65.0  & 94.4  & 83.8  & 72.1  & 93.6  & \textbf{87.5} & 75.3  & 93.6  & 218   & 1505  & 125.8 \\
		5     & 136.7 & 81.7  & 66.7  & 93.8  & 84.4  & 68.0  & 91.9  & \textbf{87.3} & 78.3  & 92.2  & 212   & 609   & 112.6 \\
		8     & 134.3 & 83.0  & 66.9  & 92.6  & 83.9  & 75.3  & 89.7  & \textbf{86.2} & 81.2  & 91.1  & 257   & 619   & 108.9 \\
		10    & 132.9 & 82.6  & 66.5  & 88.9  & 83.7  & 73.1  & 90.1  & \textbf{85.0} & 80.0  & 90.2  & 278   & 695   & 106.7 \\
		12    & 132.9 & 82.8  & 70.4  & 90.4  & 83.8  & 74.9  & 89.9  & \textbf{84.2} & 78.7  & 88.4  & 298   & 755   & 106.3 \\
		15    & 131.5 & 81.5  & 64.9  & 89.2  & 82.4  & 72.3  & 88.4  & \textbf{83.0} & 78.4  & 86.8  & 323   & 669   & 102.9 \\
		20    & 130.3 & 80.5  & 70.5  & 88.3  & \textbf{82.0} & 75.6  & 86.9  & 81.5  & 78.2  & 84.9  & 360   & 608   & 100.2 \\
		\bottomrule
	\end{tabular}%
	\label{tab:results_bischoff_ratliff}%
\end{table}%

\subsection{Results for \cite{PaquayEtAl2018} instances}
\label{sec:comp_results_paquay}

For the comparison with the approach of \cite{PaquayEtAl2018}, we use the following problem-specific parameters: padding height~$\paddingHeight=0$ and weight balance importance~$\weightBalanceImportance=100$.
For these instances, the maximum \centerOfGravity\ deviation is between 5\% and 10\% but differs in $x$- and $y$-direction for some containers.
We therefore set the maximum center of gravity deviation~$\maximumCenterOfGravityDeviation=0.05$.
Moreover, we ignore the minimum item overlap~$\minimumItemOverlap$ and only consider items non-floating if the four bottom corner points rest on stackable items or the \ULD\ base area.
Note also that the edge width~$\edgeWidth=0$ and no substructure is allowed ($\substructureAllowed=0$) for all \ULD s.

Table~\ref{tab:results_paquay_et_al} displays the comparison between the two approaches grouped by the number of items~$\numItems$ in an instance.
We report the (average) number of \ULD s~\numULDs\ (\avgNumULDs) and the median of the utilization~\medianUtilization.
For the \RGS, we additionally denote the number~\numULDsCogViolated\ (percentage \percentageCogViolated) of \ULD s for which the \centerOfGravity\ is violated.
Again, the solution time~$\solutionTime$ is reported in milliseconds.

\begin{table}[htbp]
	\centering
	\caption{Comparison of \cite{PaquayEtAl2018} and the \RGS.}
	\begin{tabular}{lrrrrrrrrrrr}
		\toprule
		& \multicolumn{3}{c}{\cite{PaquayEtAl2018}} & \multicolumn{8}{c}{\RGS} \\
		\cmidrule(lr){2-4} \cmidrule(lr){5-12} 
		$\numItems$ & \multicolumn{1}{l}{\numULDs} & \multicolumn{1}{l}{\avgNumULDs} & \multicolumn{1}{l}{\medianUtilization} & \multicolumn{1}{l}{\numULDs} & \multicolumn{1}{l}{\avgNumULDs} & \multicolumn{1}{l}{\medianUtilization} & \multicolumn{1}{l}{\avgUtilization} & \multicolumn{1}{l}{\numULDsCogViolated} & \multicolumn{1}{l}{\percentageCogViolated} & \multicolumn{1}{l}{\avgTime} & \multicolumn{1}{l}{\maxTime} \\
		\midrule
		\multicolumn{1}{r}{10} & 47    & 1.6   & 33.2  & 30    & 1.0   & \textbf{34.3} & 36.1  & 0     & 0.0   & 26.0  & 42 \\
		\multicolumn{1}{r}{20} & 54    & 1.8   & 34.6  & 35    & 1.2   & \textbf{42.1} & 40.3  & 3     & 8.6   & 69.0  & 142 \\
		\multicolumn{1}{r}{30} & 72    & 2.4   & 36.6  & 47    & 1.6   & \textbf{44.7} & 42.7  & 6     & 12.8  & 116.4 & 281 \\
		\multicolumn{1}{r}{40} & 83    & 2.8   & 37.9  & 60    & 2.0   & \textbf{46.5} & 44.1  & 4     & 6.7   & 88.9  & 105 \\
		\multicolumn{1}{r}{50} & 90    & 3.0   & 42.2  & 66    & 2.2   & \textbf{47.5} & 45.7  & 5     & 7.6   & 130.5 & 161 \\
		\multicolumn{1}{r}{60} & 108   & 3.6   & 41.5  & 80    & 2.7   & \textbf{49.6} & 46.4  & 10    & 12.5  & 172.7 & 206 \\
		\multicolumn{1}{r}{70} & 126   & 4.2   & 45.7  & 82    & 2.7   & \textbf{54.0} & 51.6  & 6     & 7.3   & 228.0 & 294 \\
		\multicolumn{1}{r}{80} & 129   & 4.3   & 47.0  & 93    & 3.1   & \textbf{54.5} & 51.1  & 12    & 12.9  & 280.6 & 345 \\
		\multicolumn{1}{r}{90} & 141   & 4.7   & 47.5  & 102   & 3.4   & \textbf{53.5} & 50.3  & 13    & 12.7  & 355.8 & 422 \\
		\multicolumn{1}{r}{100} & 160   & 5.3   & 50.2  & 113   & 3.8   & \textbf{53.0} & 50.3  & 20    & 17.7  & 429.9 & 513 \\
		\midrule
		\total & 1010  &       &       & 708   &       &       &       & 79    &       &       &  \\
		\bottomrule
	\end{tabular}%
	\label{tab:results_paquay_et_al}%
\end{table}%

All in all, the \RGS\ outperforms the two-phase constructive heuristic of \cite{PaquayEtAl2018}.
The median utilization is higher for all instance groups.
Moreover, the number of used \ULD s~\numULDs\ is significantly smaller with 708 instead of 1010 in total.
The maximum solution time of the \RGS\ is 0.5 seconds.
\cite{PaquayEtAl2018} report that the solution time is no longer than 12 seconds.
The comparison of the solution time should be treated with caution, as different machines were used.
For the \RGS, the acceptable \centerOfGravity\ deviation is violated in 79 of 708 \ULD s, which corresponds to 11.2\%.
The approach of \cite{PaquayEtAl2018} is slightly worse in this respect with 12.9\% of the loaded \ULD s violating the acceptable \centerOfGravity\ deviation.
To ensure that the values are comparable, we have analyzed the violation in the same way as \cite{PaquayEtAl2018}, \ie, the maximum \centerOfGravity\ deviation is different for some containers in $x$- and $y$-direction (which was ignored during solution generation).
Note that \cite{PaquayEtAl2018} propose two different merit functions~MF1 and MF2.
Our approach is compared to MF1, while MF2 actually yields better results for the \centerOfGravity\ restriction.
Instance-by-instance results are shown in the Online Appendix.

\subsection{Results for adapted \cite{PaquayEtAl2018} instances}
\label{sec:results_adapted_Paquay}

In order to evaluate our algorithm with all its features, we tested the \RGS\ on adapted \cite{PaquayEtAl2018} instances by setting the edge width~$\edgeWidth=10$ and the vertical edge offset~$\verticalEdgeOffset=10$.
Table~\ref{tab:results_adapted_paquay} reports results for instances with an unlimited number of available \ULD s (denoted by \emph{unlimited}) and instances where each \ULD\ type is available only once (denoted by \emph{1 \ULD}).
The column \numSubstructureUsed\ denotes for how many \ULD s a substructure is used.

The average utilization is between 37.1\% for instances with $\numItems=10$ items and 58.7\% for instances with $\numItems=100$ items.
The acceptable \centerOfGravity\ deviation is violated for 7.2\% of the \ULD s if each \ULD\ type is available indefinitely.
The utilization differs only slightly if each \ULD\ type is available once.
However, the violation of the \centerOfGravity\ deviation increases to 10.1\% of the \ULD s.
Using a substructure is beneficial for around 27\% of the \ULD s.
All instances are solved within 1.3 seconds.
All in all, the results look very promising due to a high utilization of \ULD s.
Instance-by-instance results are shown in the Online Appendix.

\setlength{\tabcolsep}{4.0pt}
\begin{table}[htbp]
	\centering
	\caption{Results for adapated \cite{PaquayEtAl2018} instances with an unlimited number of \ULD s and only one available \ULD\ of each \ULD\ type.}
	\adjustbox{max width=\textwidth}{
	\begin{tabular}{lrrrrrrrrrrrrrrrr}
		\toprule
		& \multicolumn{8}{c}{\RGS\ (unlimited)} & \multicolumn{8}{c}{\RGS\ (1 \ULD)} \\
		\cmidrule(lr){2-9} \cmidrule(lr){10-17} 
		\numItems & \multicolumn{1}{l}{\numULDs} & \multicolumn{1}{l}{\avgNumULDs} & \multicolumn{1}{l}{\avgUtilization} & \multicolumn{1}{l}{\numULDsCogViolated} & \multicolumn{1}{l}{\percentageCogViolated} & \multicolumn{1}{l}{\numSubstructureUsed} & \multicolumn{1}{l}{\avgTime} & \multicolumn{1}{l}{\maxTime} & \multicolumn{1}{l}{\numULDs} & \multicolumn{1}{l}{\avgNumULDs} & \multicolumn{1}{l}{\avgUtilization} & \multicolumn{1}{l}{\numULDsCogViolated} & \multicolumn{1}{l}{\percentageCogViolated} & \multicolumn{1}{l}{\numSubstructureUsed} & \multicolumn{1}{l}{\avgTime} & \multicolumn{1}{l}{\maxTime} \\
		\midrule
		\multicolumn{1}{r}{10} & 30    & 1.0   & 37.1  & 3     & 10.0  & 2     & 68.6  & 107   & 30    & 1.0   & 37.1  & 4     & 13.3  & 2     & 69.8  & 110 \\
		\multicolumn{1}{r}{20} & 33    & 1.1   & 45.9  & 1     & 3.0   & 1     & 158.4 & 293   & 33    & 1.1   & 45.9  & 2     & 6.1   & 2     & 160.6 & 292 \\
		\multicolumn{1}{r}{30} & 45    & 1.5   & 44.3  & 2     & 4.4   & 11    & 274.1 & 627   & 45    & 1.5   & 44.4  & 4     & 8.9   & 9     & 279.2 & 669 \\
		\multicolumn{1}{r}{40} & 57    & 1.9   & 47.3  & 5     & 8.8   & 17    & 237.8 & 860   & 57    & 1.9   & 48.1  & 2     & 3.5   & 14    & 239.5 & 854 \\
		\multicolumn{1}{r}{50} & 60    & 2.0   & 50.7  & 6     & 10.0  & 15    & 303.8 & 548   & 62    & 2.1   & 52.3  & 8     & 12.9  & 14    & 301.6 & 548 \\
		\multicolumn{1}{r}{60} & 72    & 2.4   & 51.5  & 9     & 12.5  & 20    & 392.3 & 465   & 77    & 2.6   & 54.5  & 8     & 10.4  & 27    & 383.4 & 450 \\
		\multicolumn{1}{r}{70} & 78    & 2.6   & 54.9  & 4     & 5.1   & 22    & 495.6 & 686   & 85    & 2.8   & 56.9  & 10    & 11.8  & 24    & 477.4 & 580 \\
		\multicolumn{1}{r}{80} & 86    & 2.9   & 55.1  & 6     & 7.0   & 20    & 620.0 & 850   & 97    & 3.2   & 57.5  & 10    & 10.3  & 41    & 595.1 & 731 \\
		\multicolumn{1}{r}{90} & 92    & 3.1   & 56.1  & 7     & 7.6   & 19    & 763.1 & 960   & 109   & 3.6   & 58.5  & 10    & 9.2   & 35    & 746.5 & 868 \\
		\multicolumn{1}{r}{100} & 100   & 3.3   & 58.7  & 4     & 4.0   & 29    & 923.1 & 1239  & 134   & 4.5   & 56.5  & 16    & 11.9  & 53    & 888.4 & 1049 \\
		\midrule
		\total & 653   &       &       & 47    &       & 156   &       &       & 729   &       &       & 74    &       & 221   &       &  \\
		\bottomrule
	\end{tabular}%
    }
	\label{tab:results_adapted_paquay}%
\end{table}%
\setlength{\tabcolsep}{6.0pt}

\section{Conclusions and outlook}
\label{sec:conclusion}

We have introduced an insertion heuristic embedded into a Randomized Greedy Search to solve a three-dimensional \problem\ with transportation constraints that is relevant in the air freight industry.
In particular, the problem considers load stability, (non-)stackable items, weight distribution, specially shaped \ULD s, padding material, and the usage of a substructure.
The insertion heuristic that is repeatedly called and based on extreme points follows a first fit approach resulting in a very fast algorithm.
An underlying grid structure is introduced for further acceleration.
To ensure a high loading density, the sorting of items with corresponding orientations is based on a grouping procedure that enables both layer and free packing patterns.
Moreover, new findings on extreme points were presented by extending the set of extreme points proposed in the literature \citep{CrainicEtAl2008}.
An algorithm for generating extreme points that almost never reaches the quadratic worst-case complexity in the number of loaded items was proposed.
In a computational study, we demonstrated that our algorithm is competitive with the state of the art approaches and outperforms them on most instances realistic for the air freight industry.
In comparison to the approach of \cite{PaquayEtAl2018}, the median utilization is higher for all and at least five percent higher for all but two instance groups, while the center of gravity deviation is almost the same.
Moreover, results of adapted \cite{PaquayEtAl2018} instances are promising for improving the utilization of \ULD s in practice.
We have also shown that using a substructure is beneficial for around 27\% of the \ULD s for this benchmark set.

Our work, in particular the grid acceleration, the special item sorting, and the definition, generation, and moving of extreme points, can also be extended to many other packing problems, e.g., two-dimensional packing problems or problem variants with costs~\citep{ChanEtAl2006}.
Moreover, the approach can be adapted to assure building stacks, which is relevant in land transportation to easily load and unload (groups of) items.
Another research avenue is to focus not only on a high utilization per \ULD, but to consider several \ULD s simultaneously in order to balance the weight throughout the aircraft.
A further practical requirement to tackle is to ensure that predefined groups of items are loaded into the same \ULD.

\vspace*{-0.3cm}

\section*{Acknowledgement}

The authors would like to thank Michaela Babl and Dennis Schmidt for their detailed answers to all practical questions on the problem and the colleagues from the BinPACKER Air team Ryszard Balewski, Mariusz Ciepluch, Holger K\"{o}hler and Pawe\l\ Krzemi\'{n}ski for fruitful discussions.
Special thanks go to Ivo Hedtke, who provided support with every problem and every question.

\vspace*{-0.3cm}

\bibliographystyle{natbib}

\clearpage

\section*{Online Appendix}
\begin{appendix}
	
	In this Appendix, we present instance-by-instance results.
	Moreover, we present the results for each loaded \ULD.
	The entries in the Tables~\ref{tab:detailed_Paquay}--\ref{tab:detailed_con_1ULD} have the following meaning:
	
	\begin{center}
		\begin{tabular}{rp{.85\textwidth}}
			$\numItems$: & number of items; \\
			\numInstance: & instance number; \\
			$\solutionTime$: & solution time in milliseconds;\\
			$\numLoadedItems$: & number of loaded items;\\
			\numULDs: & number of loaded \ULD s; \\
			\overallUtilization: & percentage of the total utilization (loaded volume / volume of loaded \ULD s);\\
			\containerID: & container id (or type); \\
			$\sortingCriterion$: & sorting criterion of the best solution; \\
			$\utilization$: & utilization (loaded volume / \ULD\ volume); \\
			\CoGDevX: & center of gravity deviation in $x$-direction in percent; \\
			\CoGDevY: & center of gravity deviation in $y$-direction in percent. \\
		\end{tabular}
	\end{center}

	Tables~\ref{tab:detailed_Paquay}--\ref{tab:detailed_Paquay_ULDs} display the results for \cite{PaquayEtAl2018} instances, Tables~\ref{tab:detailed_unlimited}--\ref{tab:detailed_con_unlimited} for adapted \cite{PaquayEtAl2018} instances with an unlimited number of available \ULD s, and Tables~\ref{tab:detailed_1ULD}--\ref{tab:detailed_con_1ULD} for adapted \cite{PaquayEtAl2018} instances where each \ULD\ type is available only once.
	
	\clearpage
	
	\addtocounter{table}{+1}
\begin{multicols*}{4}
	[\begin{center}Table \thetable: Detailed results for the \cite{PaquayEtAl2018} instances.\end{center}]

\setlength{\tabcolsep}{3pt}
\tiny

\centering

\tablefirsthead{%
	\toprule
	\multicolumn{1}{l}{$\numItems$} & \multicolumn{1}{l}{\numInstance} & \multicolumn{1}{l}{$\solutionTime$} & \multicolumn{1}{l}{$\numLoadedItems$} & \multicolumn{1}{l}{\numULDs} & \multicolumn{1}{l}{\overallUtilization} \\
	\midrule}
\tablehead{%
	\toprule
	\multicolumn{1}{l}{$\numItems$} & \multicolumn{1}{l}{\numInstance} & \multicolumn{1}{l}{$\solutionTime$} & \multicolumn{1}{l}{$\numLoadedItems$} & \multicolumn{1}{l}{\numULDs} & \multicolumn{1}{l}{\overallUtilization} \\
	\midrule}
\tabletail{%
	\midrule
	\multicolumn{6}{r}{\textit{Continued on next column}}\\}
\tablelasttail{\bottomrule}

\TrickSupertabularIntoMulticols

\begin{supertabular}{rrrrrr}
		10    & 0     & 15    & 10    & 1     & 24.9 \\
		10    & 1     & 38    & 10    & 1     & 27.9 \\
		10    & 2     & 30    & 10    & 1     & 41.6 \\
		10    & 3     & 34    & 10    & 1     & 27.4 \\
		10    & 4     & 21    & 10    & 1     & 56.1 \\
		10    & 5     & 8     & 10    & 1     & 13.1 \\
		10    & 6     & 23    & 10    & 1     & 32.7 \\
		10    & 7     & 22    & 10    & 1     & 58.5 \\
		10    & 8     & 31    & 10    & 1     & 38.3 \\
		10    & 9     & 30    & 10    & 1     & 48.2 \\
		10    & 10    & 20    & 10    & 1     & 26.4 \\
		10    & 11    & 42    & 10    & 1     & 26.0 \\
		10    & 12    & 21    & 10    & 1     & 46.6 \\
		10    & 13    & 14    & 10    & 1     & 55.0 \\
		10    & 14    & 22    & 10    & 1     & 62.6 \\
		10    & 15    & 23    & 10    & 1     & 48.1 \\
		10    & 16    & 38    & 10    & 1     & 35.9 \\
		10    & 17    & 34    & 10    & 1     & 25.8 \\
		10    & 18    & 34    & 10    & 1     & 29.7 \\
		10    & 19    & 21    & 10    & 1     & 42.5 \\
		10    & 20    & 24    & 10    & 1     & 25.0 \\
		10    & 21    & 34    & 10    & 1     & 38.4 \\
		10    & 22    & 16    & 10    & 1     & 15.8 \\
		10    & 23    & 35    & 10    & 1     & 27.9 \\
		10    & 24    & 35    & 10    & 1     & 26.1 \\
		10    & 25    & 16    & 10    & 1     & 23.3 \\
		10    & 26    & 38    & 10    & 1     & 16.1 \\
		10    & 27    & 23    & 10    & 1     & 51.8 \\
		10    & 28    & 22    & 10    & 1     & 50.5 \\
		10    & 29    & 16    & 10    & 1     & 41.9 \\
		20    & 0     & 142   & 20    & 1     & 50.8 \\
		20    & 1     & 25    & 20    & 1     & 39.4 \\
		20    & 2     & 48    & 20    & 1     & 42.4 \\
		20    & 3     & 108   & 20    & 1     & 47.2 \\
		20    & 4     & 76    & 20    & 1     & 48.1 \\
		20    & 5     & 25    & 20    & 1     & 33.3 \\
		20    & 6     & 43    & 20    & 1     & 27.2 \\
		20    & 7     & 118   & 20    & 1     & 61.5 \\
		20    & 8     & 69    & 20    & 1     & 43.4 \\
		20    & 9     & 117   & 20    & 1     & 42.4 \\
		20    & 10    & 34    & 20    & 2     & 42.4 \\
		20    & 11    & 39    & 20    & 1     & 47.2 \\
		20    & 12    & 25    & 20    & 1     & 40.0 \\
		20    & 13    & 24    & 20    & 2     & 39.1 \\
		20    & 14    & 75    & 20    & 1     & 35.0 \\
		20    & 15    & 38    & 20    & 1     & 30.5 \\
		20    & 16    & 108   & 20    & 1     & 44.5 \\
		20    & 17    & 110   & 20    & 1     & 55.9 \\
		20    & 18    & 47    & 20    & 1     & 35.7 \\
		20    & 19    & 24    & 20    & 2     & 41.6 \\
		20    & 20    & 27    & 20    & 1     & 46.8 \\
		20    & 21    & 54    & 20    & 1     & 46.0 \\
		20    & 22    & 32    & 20    & 2     & 38.3 \\
		20    & 23    & 120   & 20    & 1     & 24.1 \\
		20    & 24    & 116   & 20    & 1     & 37.8 \\
		20    & 25    & 120   & 20    & 1     & 50.0 \\
		20    & 26    & 24    & 20    & 1     & 30.5 \\
		20    & 27    & 128   & 20    & 1     & 46.1 \\
		20    & 28    & 127   & 20    & 1     & 48.5 \\
		20    & 29    & 28    & 20    & 2     & 37.3 \\
		30    & 0     & 53    & 30    & 2     & 46.0 \\
		30    & 1     & 281   & 30    & 1     & 42.0 \\
		30    & 2     & 247   & 30    & 1     & 61.0 \\
		30    & 3     & 55    & 30    & 2     & 39.5 \\
		30    & 4     & 49    & 30    & 2     & 50.9 \\
		30    & 5     & 54    & 30    & 2     & 50.7 \\
		30    & 6     & 48    & 30    & 2     & 49.5 \\
		30    & 7     & 250   & 30    & 1     & 52.7 \\
		30    & 8     & 249   & 30    & 1     & 50.3 \\
		30    & 9     & 48    & 30    & 2     & 47.0 \\
		30    & 10    & 52    & 30    & 1     & 38.3 \\
		30    & 11    & 246   & 30    & 1     & 34.8 \\
		30    & 12    & 52    & 30    & 2     & 44.9 \\
		30    & 13    & 254   & 30    & 1     & 36.1 \\
		30    & 14    & 49    & 30    & 2     & 48.6 \\
		30    & 15    & 262   & 30    & 1     & 42.2 \\
		30    & 16    & 54    & 30    & 2     & 39.7 \\
		30    & 17    & 51    & 30    & 2     & 42.2 \\
		30    & 18    & 51    & 30    & 2     & 47.8 \\
		30    & 19    & 60    & 30    & 2     & 39.8 \\
		30    & 20    & 47    & 30    & 2     & 49.9 \\
		30    & 21    & 58    & 30    & 2     & 35.7 \\
		30    & 22    & 151   & 30    & 1     & 45.5 \\
		30    & 23    & 52    & 30    & 1     & 45.1 \\
		30    & 24    & 247   & 30    & 1     & 45.1 \\
		30    & 25    & 52    & 30    & 2     & 50.1 \\
		30    & 26    & 56    & 30    & 2     & 39.8 \\
		30    & 27    & 52    & 30    & 2     & 46.5 \\
		30    & 28    & 262   & 30    & 1     & 64.6 \\
		30    & 29    & 50    & 30    & 1     & 43.0 \\
		40    & 0     & 100   & 40    & 2     & 48.5 \\
		40    & 1     & 79    & 40    & 2     & 50.3 \\
		40    & 2     & 95    & 40    & 2     & 46.7 \\
		40    & 3     & 80    & 40    & 2     & 37.9 \\
		40    & 4     & 86    & 40    & 2     & 43.8 \\
		40    & 5     & 77    & 40    & 2     & 50.8 \\
		40    & 6     & 91    & 40    & 2     & 45.3 \\
		40    & 7     & 95    & 40    & 2     & 35.3 \\
		40    & 8     & 99    & 40    & 2     & 57.4 \\
		40    & 9     & 90    & 40    & 2     & 57.0 \\
		40    & 10    & 81    & 40    & 2     & 41.4 \\
		40    & 11    & 96    & 40    & 2     & 45.1 \\
		40    & 12    & 89    & 40    & 2     & 41.5 \\
		40    & 13    & 86    & 40    & 2     & 47.9 \\
		40    & 14    & 102   & 40    & 2     & 42.6 \\
		40    & 15    & 83    & 40    & 2     & 46.0 \\
		40    & 16    & 92    & 40    & 2     & 38.4 \\
		40    & 17    & 88    & 40    & 2     & 52.4 \\
		40    & 18    & 78    & 40    & 2     & 45.0 \\
		40    & 19    & 104   & 40    & 2     & 41.3 \\
		40    & 20    & 85    & 40    & 2     & 49.7 \\
		40    & 21    & 80    & 40    & 2     & 48.2 \\
		40    & 22    & 97    & 40    & 2     & 53.0 \\
		40    & 23    & 87    & 40    & 2     & 49.9 \\
		40    & 24    & 80    & 40    & 2     & 47.9 \\
		40    & 25    & 105   & 40    & 2     & 43.5 \\
		40    & 26    & 95    & 40    & 2     & 52.5 \\
		40    & 27    & 68    & 40    & 2     & 48.2 \\
		40    & 28    & 99    & 40    & 2     & 55.4 \\
		40    & 29    & 80    & 40    & 2     & 41.1 \\
		50    & 0     & 144   & 50    & 2     & 54.3 \\
		50    & 1     & 130   & 50    & 2     & 53.8 \\
		50    & 2     & 115   & 50    & 3     & 46.4 \\
		50    & 3     & 143   & 50    & 2     & 50.0 \\
		50    & 4     & 118   & 50    & 3     & 43.2 \\
		50    & 5     & 112   & 50    & 3     & 43.3 \\
		50    & 6     & 122   & 50    & 2     & 51.5 \\
		50    & 7     & 118   & 50    & 3     & 46.4 \\
		50    & 8     & 143   & 50    & 2     & 50.4 \\
		50    & 9     & 126   & 50    & 2     & 47.4 \\
		50    & 10    & 148   & 50    & 2     & 47.0 \\
		50    & 11    & 122   & 50    & 2     & 45.1 \\
		50    & 12    & 125   & 50    & 2     & 52.4 \\
		50    & 13    & 136   & 50    & 2     & 48.8 \\
		50    & 14    & 131   & 50    & 2     & 45.5 \\
		50    & 15    & 138   & 50    & 2     & 44.2 \\
		50    & 16    & 138   & 50    & 2     & 50.1 \\
		50    & 17    & 128   & 50    & 2     & 53.2 \\
		50    & 18    & 135   & 50    & 2     & 45.1 \\
		50    & 19    & 161   & 50    & 2     & 54.5 \\
		50    & 20    & 134   & 50    & 2     & 43.5 \\
		50    & 21    & 130   & 50    & 2     & 51.9 \\
		50    & 22    & 124   & 50    & 3     & 40.6 \\
		50    & 23    & 136   & 50    & 2     & 51.2 \\
		50    & 24    & 128   & 50    & 2     & 51.0 \\
		50    & 25    & 122   & 50    & 3     & 44.8 \\
		50    & 26    & 116   & 50    & 2     & 51.3 \\
		50    & 27    & 137   & 50    & 2     & 50.8 \\
		50    & 28    & 133   & 50    & 2     & 50.9 \\
		50    & 29    & 123   & 50    & 2     & 47.7 \\
		60    & 0     & 205   & 60    & 3     & 44.1 \\
		60    & 1     & 177   & 60    & 3     & 44.5 \\
		60    & 2     & 177   & 60    & 3     & 47.3 \\
		60    & 3     & 150   & 60    & 3     & 49.7 \\
		60    & 4     & 206   & 60    & 3     & 48.0 \\
		60    & 5     & 170   & 60    & 3     & 50.9 \\
		60    & 6     & 140   & 60    & 3     & 48.8 \\
		60    & 7     & 184   & 60    & 3     & 44.1 \\
		60    & 8     & 174   & 60    & 2     & 48.9 \\
		60    & 9     & 178   & 60    & 3     & 49.8 \\
		60    & 10    & 161   & 60    & 3     & 45.4 \\
		60    & 11    & 180   & 60    & 2     & 56.9 \\
		60    & 12    & 139   & 60    & 3     & 51.2 \\
		60    & 13    & 193   & 60    & 3     & 49.8 \\
		60    & 14    & 158   & 60    & 2     & 50.4 \\
		60    & 15    & 142   & 60    & 3     & 56.4 \\
		60    & 16    & 141   & 60    & 3     & 55.7 \\
		60    & 17    & 192   & 60    & 2     & 52.3 \\
		60    & 18    & 185   & 60    & 3     & 45.8 \\
		60    & 19    & 190   & 60    & 2     & 48.0 \\
		60    & 20    & 157   & 60    & 3     & 48.3 \\
		60    & 21    & 199   & 60    & 2     & 46.6 \\
		60    & 22    & 190   & 60    & 2     & 49.6 \\
		60    & 23    & 195   & 60    & 2     & 60.2 \\
		60    & 24    & 189   & 60    & 2     & 51.3 \\
		60    & 25    & 146   & 60    & 3     & 53.3 \\
		60    & 26    & 177   & 60    & 2     & 55.0 \\
		60    & 27    & 182   & 60    & 3     & 43.0 \\
		60    & 28    & 147   & 60    & 3     & 52.3 \\
		60    & 29    & 156   & 60    & 3     & 52.2 \\
		70    & 0     & 218   & 70    & 3     & 49.0 \\
		70    & 1     & 260   & 70    & 2     & 63.4 \\
		70    & 2     & 226   & 70    & 3     & 52.3 \\
		70    & 3     & 259   & 70    & 3     & 46.2 \\
		70    & 4     & 211   & 70    & 3     & 54.0 \\
		70    & 5     & 223   & 70    & 3     & 49.8 \\
		70    & 6     & 264   & 70    & 3     & 52.2 \\
		70    & 7     & 294   & 70    & 2     & 61.4 \\
		70    & 8     & 205   & 70    & 3     & 50.6 \\
		70    & 9     & 210   & 70    & 3     & 54.3 \\
		70    & 10    & 223   & 70    & 3     & 51.6 \\
		70    & 11    & 216   & 70    & 3     & 55.7 \\
		70    & 12    & 193   & 70    & 3     & 51.5 \\
		70    & 13    & 227   & 70    & 3     & 53.2 \\
		70    & 14    & 257   & 70    & 2     & 53.0 \\
		70    & 15    & 210   & 70    & 3     & 48.3 \\
		70    & 16    & 245   & 70    & 3     & 52.0 \\
		70    & 17    & 194   & 70    & 3     & 54.1 \\
		70    & 18    & 234   & 70    & 2     & 55.2 \\
		70    & 19    & 244   & 70    & 3     & 46.0 \\
		70    & 20    & 190   & 70    & 3     & 53.2 \\
		70    & 21    & 272   & 70    & 2     & 56.1 \\
		70    & 22    & 238   & 70    & 3     & 54.1 \\
		70    & 23    & 211   & 70    & 3     & 52.3 \\
		70    & 24    & 215   & 70    & 2     & 54.2 \\
		70    & 25    & 245   & 70    & 2     & 60.4 \\
		70    & 26    & 188   & 70    & 3     & 53.7 \\
		70    & 27    & 232   & 70    & 3     & 51.0 \\
		70    & 28    & 177   & 70    & 3     & 56.2 \\
		70    & 29    & 259   & 70    & 2     & 51.8 \\
		80    & 0     & 268   & 80    & 3     & 55.7 \\
		80    & 1     & 264   & 80    & 3     & 62.1 \\
		80    & 2     & 261   & 80    & 3     & 55.0 \\
		80    & 3     & 300   & 80    & 3     & 55.1 \\
		80    & 4     & 256   & 80    & 3     & 55.9 \\
		80    & 5     & 233   & 80    & 4     & 56.1 \\
		80    & 6     & 300   & 80    & 3     & 51.6 \\
		80    & 7     & 262   & 80    & 4     & 45.4 \\
		80    & 8     & 269   & 80    & 3     & 53.8 \\
		80    & 9     & 308   & 80    & 3     & 52.7 \\
		80    & 10    & 300   & 80    & 3     & 52.0 \\
		80    & 11    & 289   & 80    & 3     & 59.5 \\
		80    & 12    & 345   & 80    & 3     & 43.8 \\
		80    & 13    & 252   & 80    & 3     & 60.3 \\
		80    & 14    & 314   & 80    & 2     & 57.5 \\
		80    & 15    & 260   & 80    & 3     & 56.0 \\
		80    & 16    & 278   & 80    & 3     & 50.5 \\
		80    & 17    & 288   & 80    & 3     & 49.5 \\
		80    & 18    & 240   & 80    & 3     & 59.7 \\
		80    & 19    & 272   & 80    & 3     & 58.3 \\
		80    & 20    & 334   & 80    & 3     & 54.2 \\
		80    & 21    & 267   & 80    & 4     & 50.4 \\
		80    & 22    & 297   & 80    & 3     & 51.2 \\
		80    & 23    & 271   & 80    & 3     & 51.7 \\
		80    & 24    & 327   & 80    & 3     & 52.5 \\
		80    & 25    & 307   & 80    & 3     & 49.1 \\
		80    & 26    & 279   & 80    & 3     & 55.2 \\
		80    & 27    & 254   & 80    & 3     & 51.9 \\
		80    & 28    & 285   & 80    & 3     & 51.4 \\
		80    & 29    & 239   & 80    & 4     & 50.0 \\
		90    & 0     & 397   & 90    & 3     & 49.0 \\
		90    & 1     & 360   & 90    & 3     & 49.8 \\
		90    & 2     & 379   & 90    & 4     & 54.7 \\
		90    & 3     & 405   & 90    & 3     & 52.5 \\
		90    & 4     & 333   & 90    & 4     & 48.5 \\
		90    & 5     & 422   & 90    & 3     & 52.3 \\
		90    & 6     & 385   & 90    & 4     & 52.3 \\
		90    & 7     & 327   & 90    & 3     & 59.3 \\
		90    & 8     & 391   & 90    & 3     & 50.1 \\
		90    & 9     & 285   & 90    & 3     & 60.0 \\
		90    & 10    & 364   & 90    & 4     & 51.4 \\
		90    & 11    & 333   & 90    & 3     & 54.3 \\
		90    & 12    & 328   & 90    & 3     & 57.3 \\
		90    & 13    & 417   & 90    & 3     & 54.9 \\
		90    & 14    & 368   & 90    & 4     & 54.8 \\
		90    & 15    & 342   & 90    & 3     & 56.9 \\
		90    & 16    & 348   & 90    & 3     & 55.1 \\
		90    & 17    & 359   & 90    & 3     & 57.6 \\
		90    & 18    & 329   & 90    & 3     & 54.0 \\
		90    & 19    & 404   & 90    & 3     & 57.6 \\
		90    & 20    & 323   & 90    & 4     & 47.7 \\
		90    & 21    & 288   & 90    & 4     & 51.5 \\
		90    & 22    & 372   & 90    & 3     & 56.1 \\
		90    & 23    & 333   & 90    & 3     & 56.2 \\
		90    & 24    & 326   & 90    & 4     & 49.6 \\
		90    & 25    & 378   & 90    & 3     & 51.3 \\
		90    & 26    & 399   & 90    & 4     & 46.4 \\
		90    & 27    & 341   & 90    & 4     & 51.7 \\
		90    & 28    & 303   & 90    & 4     & 51.2 \\
		90    & 29    & 336   & 90    & 4     & 52.6 \\
		100   & 0     & 484   & 100   & 3     & 53.2 \\
		100   & 1     & 406   & 100   & 4     & 55.4 \\
		100   & 2     & 385   & 100   & 4     & 50.8 \\
		100   & 3     & 407   & 100   & 4     & 52.9 \\
		100   & 4     & 413   & 100   & 4     & 47.6 \\
		100   & 5     & 513   & 100   & 3     & 49.5 \\
		100   & 6     & 362   & 100   & 4     & 58.2 \\
		100   & 7     & 419   & 100   & 3     & 62.0 \\
		100   & 8     & 462   & 100   & 4     & 52.9 \\
		100   & 9     & 394   & 100   & 4     & 50.5 \\
		100   & 10    & 436   & 100   & 4     & 48.5 \\
		100   & 11    & 494   & 100   & 4     & 47.2 \\
		100   & 12    & 399   & 100   & 4     & 55.1 \\
		100   & 13    & 490   & 100   & 3     & 58.4 \\
		100   & 14    & 448   & 100   & 4     & 52.9 \\
		100   & 15    & 444   & 100   & 4     & 49.8 \\
		100   & 16    & 429   & 100   & 4     & 51.7 \\
		100   & 17    & 413   & 100   & 3     & 59.1 \\
		100   & 18    & 457   & 100   & 3     & 56.4 \\
		100   & 19    & 455   & 100   & 4     & 50.5 \\
		100   & 20    & 431   & 100   & 3     & 56.1 \\
		100   & 21    & 435   & 100   & 4     & 51.8 \\
		100   & 22    & 440   & 100   & 4     & 48.4 \\
		100   & 23    & 394   & 100   & 3     & 65.3 \\
		100   & 24    & 404   & 100   & 4     & 50.6 \\
		100   & 25    & 409   & 100   & 4     & 49.3 \\
		100   & 26    & 414   & 100   & 4     & 52.3 \\
		100   & 27    & 443   & 100   & 4     & 53.1 \\
		100   & 28    & 398   & 100   & 4     & 56.2 \\
		100   & 29    & 420   & 100   & 5     & 43.7 \\
\end{supertabular}{\makeatletter\def\@currentlabel{\thetable}\label{tab:detailed_Paquay}}
\setlength{\tabcolsep}{6pt}
\end{multicols*}
	\addtocounter{table}{+1}
\begin{multicols*}{4}
	[\begin{center}Table \thetable: Detailed \ULD\ results for \cite{PaquayEtAl2018} instances.\end{center}]

\setlength{\tabcolsep}{1pt}
\tiny

\centering

\tablefirsthead{%
	\toprule
	\multicolumn{1}{l}{$\numItems$} & \multicolumn{1}{l}{\numInstance} & \containerID & $\sortingCriterion$ & \multicolumn{1}{l}{$\utilization$} & \multicolumn{1}{l}{\CoGDevX} & \multicolumn{1}{l}{\CoGDevY} \\
	\midrule}
\tablehead{%
	\toprule
	\multicolumn{1}{l}{$\numItems$} & \multicolumn{1}{l}{\numInstance} & \containerID & $\sortingCriterion$ & \multicolumn{1}{l}{$\utilization$} & \multicolumn{1}{l}{\CoGDevX} & \multicolumn{1}{l}{\CoGDevY} \\
	\midrule}
\tabletail{%
	\midrule
	\multicolumn{7}{r}{\textit{Continued on next column}}\\}
\tablelasttail{\bottomrule}

\TrickSupertabularIntoMulticols

\begin{supertabular}{rrllrrr}
		10    & 0     & 5     & \stackCumVol & 24.9  & 1.4   & 0.7 \\
		10    & 1     & 4     & \stackCumVol & 27.9  & 0.7   & 0.5 \\
		10    & 2     & 4     & \stackCumVol & 41.6  & 3.4   & 1.2 \\
		10    & 3     & 4     & \stackCumVol & 27.4  & 1.4   & 3.2 \\
		10    & 4     & 3     & \stackCumVol & 56.1  & 2.6   & 3.7 \\
		10    & 5     & 5     & \stackCumVol & 13.1  & 1.5   & 0.0 \\
		10    & 6     & 3     & \stackCumVol & 32.7  & 2.0   & 1.7 \\
		10    & 7     & 3     & \stackCumVol & 58.5  & 5.2   & 0.8 \\
\end{supertabular}{\makeatletter\def\@currentlabel{\thetable}\label{tab:detailed_Paquay_ULDs}}
\setlength{\tabcolsep}{6pt}
\end{multicols*}
	
	\addtocounter{table}{+1}
\begin{multicols*}{4}
	[\begin{center}Table \thetable: Detailed results for adapted \cite{PaquayEtAl2018} instances with an unlimited number of available \ULD s.\end{center}]

\setlength{\tabcolsep}{3pt}
\tiny

\centering

\tablefirsthead{%
	\toprule
	\multicolumn{1}{l}{$\numItems$} & \multicolumn{1}{l}{\numInstance} & \multicolumn{1}{l}{$\solutionTime$} & \multicolumn{1}{l}{$\numLoadedItems$} & \multicolumn{1}{l}{\numULDs} & \multicolumn{1}{l}{\overallUtilization} \\
	\midrule}
\tablehead{%
	\toprule
	\multicolumn{1}{l}{$\numItems$} & \multicolumn{1}{l}{\numInstance} & \multicolumn{1}{l}{$\solutionTime$} & \multicolumn{1}{l}{$\numLoadedItems$} & \multicolumn{1}{l}{\numULDs} & \multicolumn{1}{l}{\overallUtilization} \\
	\midrule}
\tabletail{%
	\midrule
	\multicolumn{6}{r}{\textit{Continued on next column}}\\}
\tablelasttail{\bottomrule}

\TrickSupertabularIntoMulticols

\begin{supertabular}{rrrrrr}
		10    & 0     & 41    & 10    & 1     & 24.9 \\
		10    & 1     & 98    & 10    & 1     & 27.9 \\
		10    & 2     & 82    & 10    & 1     & 41.6 \\
		10    & 3     & 92    & 10    & 1     & 27.4 \\
		10    & 4     & 57    & 10    & 1     & 56.1 \\
		10    & 5     & 22    & 10    & 1     & 13.1 \\
		10    & 6     & 44    & 10    & 1     & 48.5 \\
		10    & 7     & 60    & 10    & 1     & 58.5 \\
		10    & 8     & 82    & 10    & 1     & 38.3 \\
		10    & 9     & 78    & 10    & 1     & 48.2 \\
		10    & 10    & 55    & 10    & 1     & 26.4 \\
		10    & 11    & 107   & 10    & 1     & 26.0 \\
		10    & 12    & 58    & 10    & 1     & 46.6 \\
		10    & 13    & 40    & 10    & 1     & 55.0 \\
		10    & 14    & 62    & 10    & 1     & 62.6 \\
		10    & 15    & 63    & 10    & 1     & 48.1 \\
		10    & 16    & 97    & 10    & 1     & 35.9 \\
		10    & 17    & 93    & 10    & 1     & 25.8 \\
		10    & 18    & 93    & 10    & 1     & 29.7 \\
		10    & 19    & 58    & 10    & 1     & 42.5 \\
		10    & 20    & 43    & 10    & 1     & 37.1 \\
		10    & 21    & 89    & 10    & 1     & 38.4 \\
		10    & 22    & 42    & 10    & 1     & 15.8 \\
		10    & 23    & 94    & 10    & 1     & 27.9 \\
		10    & 24    & 96    & 10    & 1     & 26.1 \\
		10    & 25    & 43    & 10    & 1     & 23.3 \\
		10    & 26    & 102   & 10    & 1     & 16.1 \\
		10    & 27    & 62    & 10    & 1     & 51.8 \\
		10    & 28    & 61    & 10    & 1     & 50.5 \\
		10    & 29    & 43    & 10    & 1     & 41.9 \\
		20    & 0     & 264   & 20    & 1     & 58.9 \\
		20    & 1     & 60    & 20    & 1     & 39.4 \\
		20    & 2     & 118   & 20    & 1     & 42.4 \\
		20    & 3     & 258   & 20    & 1     & 47.2 \\
		20    & 4     & 180   & 20    & 1     & 48.1 \\
		20    & 5     & 62    & 20    & 1     & 33.3 \\
		20    & 6     & 106   & 20    & 1     & 27.2 \\
		20    & 7     & 284   & 20    & 1     & 61.5 \\
		20    & 8     & 173   & 20    & 1     & 43.4 \\
		20    & 9     & 219   & 20    & 1     & 67.3 \\
		20    & 10    & 74    & 20    & 2     & 55.3 \\
		20    & 11    & 93    & 20    & 1     & 47.2 \\
		20    & 12    & 61    & 20    & 1     & 40.0 \\
		20    & 13    & 293   & 20    & 1     & 62.6 \\
		20    & 14    & 119   & 20    & 1     & 40.5 \\
		20    & 15    & 96    & 20    & 1     & 30.5 \\
		20    & 16    & 260   & 20    & 1     & 44.5 \\
		20    & 17    & 265   & 20    & 1     & 55.9 \\
		20    & 18    & 113   & 20    & 1     & 35.7 \\
		20    & 19    & 59    & 20    & 2     & 41.6 \\
		20    & 20    & 66    & 20    & 1     & 46.8 \\
		20    & 21    & 128   & 20    & 1     & 46.0 \\
		20    & 22    & 67    & 20    & 2     & 49.9 \\
		20    & 23    & 227   & 20    & 1     & 38.2 \\
		20    & 24    & 276   & 20    & 1     & 51.8 \\
		20    & 25    & 232   & 20    & 1     & 57.9 \\
		20    & 26    & 62    & 20    & 1     & 30.5 \\
		20    & 27    & 240   & 20    & 1     & 53.4 \\
		20    & 28    & 235   & 20    & 1     & 56.2 \\
		20    & 29    & 62    & 20    & 1     & 46.6 \\
		30    & 0     & 124   & 30    & 2     & 49.2 \\
		30    & 1     & 627   & 30    & 1     & 57.6 \\
		30    & 2     & 556   & 30    & 1     & 61.0 \\
		30    & 3     & 132   & 30    & 2     & 39.5 \\
		30    & 4     & 113   & 30    & 2     & 50.9 \\
		30    & 5     & 128   & 30    & 2     & 50.7 \\
		30    & 6     & 113   & 30    & 2     & 49.5 \\
		30    & 7     & 571   & 30    & 1     & 52.7 \\
		30    & 8     & 566   & 30    & 1     & 50.3 \\
		30    & 9     & 114   & 30    & 2     & 50.2 \\
		30    & 10    & 122   & 30    & 1     & 38.3 \\
		30    & 11    & 556   & 30    & 1     & 47.8 \\
		30    & 12    & 118   & 30    & 2     & 47.1 \\
		30    & 13    & 577   & 30    & 1     & 49.4 \\
		30    & 14    & 117   & 30    & 2     & 48.6 \\
		30    & 15    & 479   & 30    & 1     & 67.0 \\
		30    & 16    & 119   & 30    & 2     & 51.7 \\
		30    & 17    & 115   & 30    & 1     & 52.7 \\
		30    & 18    & 123   & 30    & 2     & 47.8 \\
		30    & 19    & 143   & 30    & 2     & 39.8 \\
		30    & 20    & 110   & 30    & 2     & 53.3 \\
		30    & 21    & 131   & 30    & 2     & 46.5 \\
		30    & 22    & 346   & 30    & 1     & 45.5 \\
		30    & 23    & 122   & 30    & 1     & 45.1 \\
		30    & 24    & 574   & 30    & 1     & 45.1 \\
		30    & 25    & 122   & 30    & 2     & 46.8 \\
		30    & 26    & 126   & 30    & 2     & 51.9 \\
		30    & 27    & 481   & 30    & 1     & 54.3 \\
		30    & 28    & 580   & 30    & 1     & 64.6 \\
		30    & 29    & 119   & 30    & 1     & 43.0 \\
		40    & 0     & 217   & 40    & 2     & 51.5 \\
		40    & 1     & 183   & 40    & 2     & 50.3 \\
		40    & 2     & 217   & 40    & 2     & 46.7 \\
		40    & 3     & 184   & 40    & 2     & 49.5 \\
		40    & 4     & 197   & 40    & 2     & 46.9 \\
		40    & 5     & 184   & 40    & 2     & 50.8 \\
		40    & 6     & 199   & 40    & 2     & 48.1 \\
		40    & 7     & 860   & 40    & 1     & 57.6 \\
		40    & 8     & 205   & 40    & 2     & 60.9 \\
		40    & 9     & 212   & 40    & 2     & 60.5 \\
		40    & 10    & 536   & 40    & 1     & 51.8 \\
		40    & 11    & 209   & 40    & 2     & 45.1 \\
		40    & 12    & 199   & 40    & 2     & 54.1 \\
		40    & 13    & 201   & 40    & 2     & 47.9 \\
		40    & 14    & 220   & 40    & 2     & 52.3 \\
		40    & 15    & 186   & 40    & 2     & 46.0 \\
		40    & 16    & 197   & 40    & 2     & 50.1 \\
		40    & 17    & 210   & 40    & 2     & 52.4 \\
		40    & 18    & 176   & 40    & 2     & 58.8 \\
		40    & 19    & 378   & 40    & 1     & 51.6 \\
		40    & 20    & 205   & 40    & 2     & 49.7 \\
		40    & 21    & 191   & 40    & 2     & 48.2 \\
		40    & 22    & 211   & 40    & 2     & 56.2 \\
		40    & 23    & 195   & 40    & 2     & 49.9 \\
		40    & 24    & 195   & 40    & 2     & 47.9 \\
		40    & 25    & 204   & 40    & 2     & 53.4 \\
		40    & 26    & 216   & 40    & 2     & 52.5 \\
		40    & 27    & 160   & 40    & 2     & 51.6 \\
		40    & 28    & 207   & 40    & 2     & 58.8 \\
		40    & 29    & 181   & 40    & 2     & 53.6 \\
		50    & 0     & 330   & 50    & 2     & 62.8 \\
		50    & 1     & 293   & 50    & 2     & 57.5 \\
		50    & 2     & 286   & 50    & 2     & 50.3 \\
		50    & 3     & 318   & 50    & 2     & 50.0 \\
		50    & 4     & 265   & 50    & 3     & 52.4 \\
		50    & 5     & 249   & 50    & 3     & 52.5 \\
		50    & 6     & 268   & 50    & 2     & 64.1 \\
		50    & 7     & 304   & 50    & 2     & 50.3 \\
		50    & 8     & 325   & 50    & 2     & 50.4 \\
		50    & 9     & 280   & 50    & 1     & 59.3 \\
		50    & 10    & 295   & 50    & 2     & 57.7 \\
		50    & 11    & 297   & 50    & 2     & 52.1 \\
		50    & 12    & 280   & 50    & 2     & 56.1 \\
		50    & 13    & 309   & 50    & 2     & 56.4 \\
		50    & 14    & 303   & 50    & 2     & 45.5 \\
		50    & 15    & 287   & 50    & 2     & 47.0 \\
		50    & 16    & 288   & 50    & 2     & 53.2 \\
		50    & 17    & 268   & 50    & 2     & 65.3 \\
		50    & 18    & 310   & 50    & 2     & 45.1 \\
		50    & 19    & 323   & 50    & 2     & 57.8 \\
		50    & 20    & 548   & 50    & 1     & 54.3 \\
		50    & 21    & 289   & 50    & 2     & 51.9 \\
		50    & 22    & 343   & 50    & 2     & 61.7 \\
		50    & 23    & 308   & 50    & 2     & 51.2 \\
		50    & 24    & 305   & 50    & 2     & 51.0 \\
		50    & 25    & 327   & 50    & 2     & 50.4 \\
		50    & 26    & 259   & 50    & 2     & 51.3 \\
		50    & 27    & 320   & 50    & 2     & 50.8 \\
		50    & 28    & 264   & 50    & 2     & 62.5 \\
		50    & 29    & 274   & 50    & 2     & 47.7 \\
		60    & 0     & 421   & 60    & 3     & 53.5 \\
		60    & 1     & 386   & 60    & 3     & 52.0 \\
		60    & 2     & 465   & 60    & 2     & 59.3 \\
		60    & 3     & 412   & 60    & 2     & 55.9 \\
		60    & 4     & 417   & 60    & 3     & 49.8 \\
		60    & 5     & 458   & 60    & 2     & 63.8 \\
		60    & 6     & 310   & 60    & 3     & 50.6 \\
		60    & 7     & 389   & 60    & 3     & 53.5 \\
		60    & 8     & 391   & 60    & 2     & 48.9 \\
		60    & 9     & 457   & 60    & 2     & 62.4 \\
		60    & 10    & 358   & 60    & 3     & 53.0 \\
		60    & 11    & 389   & 60    & 2     & 56.9 \\
		60    & 12    & 370   & 60    & 2     & 64.2 \\
		60    & 13    & 370   & 60    & 3     & 49.8 \\
		60    & 14    & 353   & 60    & 2     & 50.4 \\
		60    & 15    & 328   & 60    & 3     & 58.6 \\
		60    & 16    & 325   & 60    & 3     & 57.2 \\
		60    & 17    & 431   & 60    & 2     & 60.4 \\
		60    & 18    & 388   & 60    & 3     & 53.5 \\
		60    & 19    & 394   & 60    & 2     & 58.9 \\
		60    & 20    & 411   & 60    & 2     & 52.3 \\
		60    & 21    & 427   & 60    & 2     & 57.2 \\
		60    & 22    & 403   & 60    & 2     & 49.6 \\
		60    & 23    & 437   & 60    & 2     & 60.2 \\
		60    & 24    & 424   & 60    & 2     & 51.3 \\
		60    & 25    & 375   & 60    & 2     & 57.8 \\
		60    & 26    & 368   & 60    & 2     & 58.4 \\
		60    & 27    & 384   & 60    & 3     & 50.3 \\
		60    & 28    & 388   & 60    & 2     & 56.6 \\
		60    & 29    & 339   & 60    & 3     & 54.1 \\
		70    & 0     & 455   & 70    & 3     & 57.2 \\
		70    & 1     & 522   & 70    & 2     & 67.3 \\
		70    & 2     & 581   & 70    & 2     & 58.8 \\
		70    & 3     & 560   & 70    & 3     & 54.0 \\
		70    & 4     & 450   & 70    & 3     & 54.0 \\
		70    & 5     & 473   & 70    & 3     & 58.2 \\
		70    & 6     & 530   & 70    & 3     & 54.2 \\
		70    & 7     & 686   & 70    & 2     & 53.1 \\
		70    & 8     & 442   & 70    & 3     & 52.5 \\
		70    & 9     & 458   & 70    & 3     & 54.3 \\
		70    & 10    & 476   & 70    & 3     & 53.6 \\
		70    & 11    & 544   & 70    & 2     & 62.7 \\
		70    & 12    & 412   & 70    & 3     & 60.3 \\
		70    & 13    & 469   & 70    & 3     & 55.2 \\
		70    & 14    & 565   & 70    & 2     & 61.3 \\
		70    & 15    & 431   & 70    & 3     & 55.1 \\
		70    & 16    & 604   & 70    & 2     & 65.2 \\
		70    & 17    & 452   & 70    & 3     & 54.1 \\
		70    & 18    & 496   & 70    & 2     & 63.8 \\
		70    & 19    & 490   & 70    & 3     & 55.8 \\
		70    & 20    & 403   & 70    & 3     & 62.2 \\
		70    & 21    & 532   & 70    & 2     & 56.1 \\
		70    & 22    & 505   & 70    & 3     & 54.1 \\
		70    & 23    & 463   & 70    & 3     & 61.1 \\
		70    & 24    & 445   & 70    & 2     & 54.2 \\
		70    & 25    & 509   & 70    & 2     & 64.2 \\
		70    & 26    & 423   & 70    & 3     & 53.7 \\
		70    & 27    & 568   & 70    & 2     & 57.4 \\
		70    & 28    & 394   & 70    & 3     & 58.3 \\
		70    & 29    & 531   & 70    & 2     & 54.9 \\
		80    & 0     & 705   & 80    & 2     & 62.7 \\
		80    & 1     & 534   & 80    & 3     & 64.5 \\
		80    & 2     & 572   & 80    & 3     & 57.0 \\
		80    & 3     & 670   & 80    & 3     & 57.2 \\
		80    & 4     & 560   & 80    & 3     & 55.9 \\
		80    & 5     & 565   & 80    & 3     & 65.1 \\
		80    & 6     & 664   & 80    & 3     & 58.8 \\
		80    & 7     & 558   & 80    & 4     & 50.7 \\
		80    & 8     & 559   & 80    & 3     & 62.9 \\
		80    & 9     & 676   & 80    & 3     & 52.7 \\
		80    & 10    & 628   & 80    & 3     & 52.0 \\
		80    & 11    & 605   & 80    & 3     & 59.5 \\
		80    & 12    & 801   & 80    & 2     & 57.6 \\
		80    & 13    & 539   & 80    & 3     & 60.3 \\
		80    & 14    & 663   & 80    & 2     & 66.5 \\
		80    & 15    & 567   & 80    & 3     & 56.0 \\
		80    & 16    & 585   & 80    & 3     & 50.5 \\
		80    & 17    & 601   & 80    & 3     & 57.9 \\
		80    & 18    & 515   & 80    & 3     & 62.0 \\
		80    & 19    & 615   & 80    & 3     & 58.3 \\
		80    & 20    & 850   & 80    & 2     & 61.0 \\
		80    & 21    & 630   & 80    & 3     & 53.3 \\
		80    & 22    & 627   & 80    & 3     & 58.4 \\
		80    & 23    & 546   & 80    & 3     & 59.0 \\
		80    & 24    & 762   & 80    & 2     & 59.0 \\
		80    & 25    & 663   & 80    & 3     & 57.4 \\
		80    & 26    & 699   & 80    & 2     & 62.1 \\
		80    & 27    & 552   & 80    & 3     & 51.9 \\
		80    & 28    & 572   & 80    & 3     & 58.7 \\
		80    & 29    & 517   & 80    & 4     & 57.3 \\
		90    & 0     & 823   & 90    & 3     & 55.9 \\
		90    & 1     & 744   & 90    & 3     & 56.9 \\
		90    & 2     & 836   & 90    & 3     & 65.9 \\
		90    & 3     & 835   & 90    & 3     & 52.5 \\
		90    & 4     & 706   & 90    & 4     & 55.6 \\
		90    & 5     & 828   & 90    & 3     & 61.1 \\
		90    & 6     & 883   & 90    & 3     & 60.7 \\
		90    & 7     & 715   & 90    & 3     & 59.3 \\
		90    & 8     & 820   & 90    & 3     & 60.8 \\
		90    & 9     & 613   & 90    & 3     & 60.0 \\
		90    & 10    & 772   & 90    & 4     & 51.4 \\
		90    & 11    & 724   & 90    & 3     & 54.3 \\
		90    & 12    & 698   & 90    & 3     & 59.4 \\
		90    & 13    & 960   & 90    & 2     & 61.7 \\
		90    & 14    & 816   & 90    & 3     & 57.9 \\
		90    & 15    & 680   & 90    & 3     & 56.9 \\
		90    & 16    & 735   & 90    & 3     & 55.1 \\
		90    & 17    & 736   & 90    & 3     & 63.4 \\
		90    & 18    & 709   & 90    & 3     & 59.4 \\
		90    & 19    & 879   & 90    & 3     & 57.6 \\
		90    & 20    & 732   & 90    & 3     & 65.8 \\
		90    & 21    & 739   & 90    & 3     & 62.3 \\
		90    & 22    & 780   & 90    & 3     & 56.1 \\
		90    & 23    & 706   & 90    & 3     & 56.2 \\
		90    & 24    & 695   & 90    & 3     & 61.3 \\
		90    & 25    & 784   & 90    & 3     & 60.0 \\
		90    & 26    & 854   & 90    & 3     & 57.3 \\
		90    & 27    & 747   & 90    & 3     & 62.3 \\
		90    & 28    & 655   & 90    & 3     & 63.3 \\
		90    & 29    & 688   & 90    & 4     & 54.0 \\
		100   & 0     & 1239  & 100   & 2     & 59.9 \\
		100   & 1     & 935   & 100   & 3     & 66.0 \\
		100   & 2     & 866   & 100   & 3     & 61.5 \\
		100   & 3     & 843   & 100   & 4     & 59.1 \\
		100   & 4     & 910   & 100   & 3     & 63.3 \\
		100   & 5     & 1045  & 100   & 3     & 56.5 \\
		100   & 6     & 743   & 100   & 4     & 65.0 \\
		100   & 7     & 873   & 100   & 3     & 62.0 \\
		100   & 8     & 1020  & 100   & 3     & 63.8 \\
		100   & 9     & 891   & 100   & 3     & 61.1 \\
		100   & 10    & 1005  & 100   & 3     & 67.0 \\
		100   & 11    & 1070  & 100   & 4     & 52.7 \\
		100   & 12    & 830   & 100   & 4     & 60.7 \\
		100   & 13    & 1014  & 100   & 3     & 58.4 \\
		100   & 14    & 1022  & 100   & 3     & 57.2 \\
		100   & 15    & 897   & 100   & 4     & 57.0 \\
		100   & 16    & 925   & 100   & 4     & 59.3 \\
		100   & 17    & 823   & 100   & 3     & 61.4 \\
		100   & 18    & 973   & 100   & 3     & 65.9 \\
		100   & 19    & 970   & 100   & 3     & 60.1 \\
		100   & 20    & 903   & 100   & 3     & 64.0 \\
		100   & 21    & 871   & 100   & 4     & 59.5 \\
		100   & 22    & 859   & 100   & 4     & 59.6 \\
		100   & 23    & 834   & 100   & 3     & 65.3 \\
		100   & 24    & 910   & 100   & 3     & 67.3 \\
		100   & 25    & 782   & 100   & 4     & 60.7 \\
		100   & 26    & 884   & 100   & 4     & 52.3 \\
		100   & 27    & 946   & 100   & 3     & 63.2 \\
		100   & 28    & 904   & 100   & 3     & 66.9 \\
		100   & 29    & 907   & 100   & 4     & 55.7 \\
\end{supertabular}{\makeatletter\def\@currentlabel{\thetable}\label{tab:detailed_unlimited}}
\setlength{\tabcolsep}{6pt}
\end{multicols*}
	\addtocounter{table}{+1}
\begin{multicols*}{4}
	[\begin{center}Table \thetable: Detailed \ULD\ results for adapted \cite{PaquayEtAl2018} instances with an unlimited number of available \ULD s.\end{center}]

\setlength{\tabcolsep}{1pt}
\tiny

\centering

\tablefirsthead{%
	\toprule
	\multicolumn{1}{l}{$\numItems$} & \multicolumn{1}{l}{\numInstance} & \containerID & $\sortingCriterion$ & \multicolumn{1}{l}{$\utilization$} & \multicolumn{1}{l}{\CoGDevX} & \multicolumn{1}{l}{\CoGDevY} \\
	\midrule}
\tablehead{%
	\toprule
	\multicolumn{1}{l}{$\numItems$} & \multicolumn{1}{l}{\numInstance} & \containerID & $\sortingCriterion$ & \multicolumn{1}{l}{$\utilization$} & \multicolumn{1}{l}{\CoGDevX} & \multicolumn{1}{l}{\CoGDevY} \\
	\midrule}
\tabletail{%
	\midrule
	\multicolumn{7}{r}{\textit{Continued on next column}}\\}
\tablelasttail{\bottomrule}

\TrickSupertabularIntoMulticols

\begin{supertabular}{rrllrrr}
		10    & 0     & 5     & \stackCumVol & 24.9  & 3.0   & 2.6 \\
		10    & 1     & 4     & \stackCumVol & 27.9  & 6.5   & 6.6 \\
		10    & 2     & 4     & \stackCumVol & 41.6  & 13.9  & 8.9 \\
		10    & 3     & 4     & \stackCumVol & 27.4  & 1.8   & 6.8 \\
		10    & 4     & 3     & \stackCumVol & 56.1  & 0.2   & 2.9 \\
		10    & 5     & 5     & \stackCumVol & 13.1  & 1.8   & 0.0 \\
		10    & 6     & 1     & \cumVol & 48.5  & 9.3   & 24.7 \\
		10    & 7     & 3     & \stackCumVol & 58.5  & 1.5   & 6.6 \\
\end{supertabular}{\makeatletter\def\@currentlabel{\thetable}\label{tab:detailed_con_unlimited}}
\setlength{\tabcolsep}{6pt}
\end{multicols*}
	
	\addtocounter{table}{+1}
\begin{multicols*}{4}
	[\begin{center}Table \thetable: Detailed results for adapted \cite{PaquayEtAl2018} instances where each \ULD\ type is available only once.\end{center}]

\setlength{\tabcolsep}{3pt}
\tiny

\centering

\tablefirsthead{%
	\toprule
	\multicolumn{1}{l}{$\numItems$} & \multicolumn{1}{l}{\numInstance} & \multicolumn{1}{l}{$\solutionTime$} & \multicolumn{1}{l}{$\numLoadedItems$} & \multicolumn{1}{l}{\numULDs} & \multicolumn{1}{l}{\overallUtilization} \\
	\midrule}
\tablehead{%
	\toprule
	\multicolumn{1}{l}{$\numItems$} & \multicolumn{1}{l}{\numInstance} & \multicolumn{1}{l}{$\solutionTime$} & \multicolumn{1}{l}{$\numLoadedItems$} & \multicolumn{1}{l}{\numULDs} & \multicolumn{1}{l}{\overallUtilization} \\
	\midrule}
\tabletail{%
	\midrule
	\multicolumn{6}{r}{\textit{Continued on next column}}\\}
\tablelasttail{\bottomrule}

\TrickSupertabularIntoMulticols

\begin{supertabular}{rrrrrr}
		10    & 0     & 41    & 10    & 1     & 24.9 \\
		10    & 1     & 101   & 10    & 1     & 27.9 \\
		10    & 2     & 81    & 10    & 1     & 41.6 \\
		10    & 3     & 94    & 10    & 1     & 27.4 \\
		10    & 4     & 57    & 10    & 1     & 56.1 \\
		10    & 5     & 22    & 10    & 1     & 13.1 \\
		10    & 6     & 44    & 10    & 1     & 48.5 \\
		10    & 7     & 60    & 10    & 1     & 58.5 \\
		10    & 8     & 94    & 10    & 1     & 38.3 \\
		10    & 9     & 80    & 10    & 1     & 48.2 \\
		10    & 10    & 56    & 10    & 1     & 26.4 \\
		10    & 11    & 110   & 10    & 1     & 26.0 \\
		10    & 12    & 59    & 10    & 1     & 46.6 \\
		10    & 13    & 40    & 10    & 1     & 55.0 \\
		10    & 14    & 61    & 10    & 1     & 62.6 \\
		10    & 15    & 64    & 10    & 1     & 48.1 \\
		10    & 16    & 98    & 10    & 1     & 35.9 \\
		10    & 17    & 96    & 10    & 1     & 25.8 \\
		10    & 18    & 93    & 10    & 1     & 29.7 \\
		10    & 19    & 58    & 10    & 1     & 42.5 \\
		10    & 20    & 44    & 10    & 1     & 37.1 \\
		10    & 21    & 90    & 10    & 1     & 38.4 \\
		10    & 22    & 43    & 10    & 1     & 15.8 \\
		10    & 23    & 96    & 10    & 1     & 27.9 \\
		10    & 24    & 97    & 10    & 1     & 26.1 \\
		10    & 25    & 44    & 10    & 1     & 23.3 \\
		10    & 26    & 103   & 10    & 1     & 16.1 \\
		10    & 27    & 63    & 10    & 1     & 51.8 \\
		10    & 28    & 61    & 10    & 1     & 50.5 \\
		10    & 29    & 43    & 10    & 1     & 41.9 \\
		20    & 0     & 265   & 20    & 1     & 58.9 \\
		20    & 1     & 60    & 20    & 1     & 39.4 \\
		20    & 2     & 121   & 20    & 1     & 42.4 \\
		20    & 3     & 259   & 20    & 1     & 47.2 \\
		20    & 4     & 183   & 20    & 1     & 48.1 \\
		20    & 5     & 64    & 20    & 1     & 33.3 \\
		20    & 6     & 104   & 20    & 1     & 27.2 \\
		20    & 7     & 287   & 20    & 1     & 61.5 \\
		20    & 8     & 175   & 20    & 1     & 43.4 \\
		20    & 9     & 222   & 20    & 1     & 67.3 \\
		20    & 10    & 73    & 20    & 2     & 55.3 \\
		20    & 11    & 96    & 20    & 1     & 47.2 \\
		20    & 12    & 62    & 20    & 1     & 40.0 \\
		20    & 13    & 292   & 20    & 1     & 62.6 \\
		20    & 14    & 121   & 20    & 1     & 40.5 \\
		20    & 15    & 98    & 20    & 1     & 30.5 \\
		20    & 16    & 267   & 20    & 1     & 44.5 \\
		20    & 17    & 268   & 20    & 1     & 55.9 \\
		20    & 18    & 116   & 20    & 1     & 35.7 \\
		20    & 19    & 59    & 20    & 2     & 41.6 \\
		20    & 20    & 67    & 20    & 1     & 46.8 \\
		20    & 21    & 130   & 20    & 1     & 46.0 \\
		20    & 22    & 69    & 20    & 2     & 49.9 \\
		20    & 23    & 230   & 20    & 1     & 38.2 \\
		20    & 24    & 280   & 20    & 1     & 51.8 \\
		20    & 25    & 238   & 20    & 1     & 57.9 \\
		20    & 26    & 62    & 20    & 1     & 30.5 \\
		20    & 27    & 241   & 20    & 1     & 53.4 \\
		20    & 28    & 244   & 20    & 1     & 56.2 \\
		20    & 29    & 64    & 20    & 1     & 46.6 \\
		30    & 0     & 129   & 30    & 2     & 49.2 \\
		30    & 1     & 669   & 30    & 1     & 57.6 \\
		30    & 2     & 570   & 30    & 1     & 61.0 \\
		30    & 3     & 135   & 30    & 2     & 39.5 \\
		30    & 4     & 112   & 30    & 2     & 50.9 \\
		30    & 5     & 130   & 30    & 2     & 50.7 \\
		30    & 6     & 115   & 30    & 2     & 49.5 \\
		30    & 7     & 579   & 30    & 1     & 52.7 \\
		30    & 8     & 570   & 30    & 1     & 50.3 \\
		30    & 9     & 115   & 30    & 2     & 50.2 \\
		30    & 10    & 124   & 30    & 1     & 38.3 \\
		30    & 11    & 573   & 30    & 1     & 47.8 \\
		30    & 12    & 120   & 30    & 2     & 47.1 \\
		30    & 13    & 579   & 30    & 1     & 49.4 \\
		30    & 14    & 117   & 30    & 2     & 48.6 \\
		30    & 15    & 485   & 30    & 1     & 67.0 \\
		30    & 16    & 121   & 30    & 2     & 51.7 \\
		30    & 17    & 116   & 30    & 1     & 52.7 \\
		30    & 18    & 124   & 30    & 2     & 47.8 \\
		30    & 19    & 144   & 30    & 2     & 39.8 \\
		30    & 20    & 110   & 30    & 2     & 53.3 \\
		30    & 21    & 132   & 30    & 2     & 46.5 \\
		30    & 22    & 348   & 30    & 1     & 45.5 \\
		30    & 23    & 123   & 30    & 1     & 45.1 \\
		30    & 24    & 583   & 30    & 1     & 45.1 \\
		30    & 25    & 123   & 30    & 2     & 46.8 \\
		30    & 26    & 127   & 30    & 2     & 51.9 \\
		30    & 27    & 488   & 30    & 1     & 54.3 \\
		30    & 28    & 595   & 30    & 1     & 64.6 \\
		30    & 29    & 119   & 30    & 1     & 43.0 \\
		40    & 0     & 221   & 40    & 2     & 51.5 \\
		40    & 1     & 183   & 40    & 2     & 50.3 \\
		40    & 2     & 214   & 40    & 2     & 46.7 \\
		40    & 3     & 186   & 40    & 2     & 49.5 \\
		40    & 4     & 195   & 40    & 2     & 46.9 \\
		40    & 5     & 185   & 40    & 2     & 50.8 \\
		40    & 6     & 200   & 40    & 2     & 48.1 \\
		40    & 7     & 854   & 40    & 1     & 57.6 \\
		40    & 8     & 211   & 40    & 2     & 60.9 \\
		40    & 9     & 215   & 40    & 2     & 60.5 \\
		40    & 10    & 543   & 40    & 1     & 51.8 \\
		40    & 11    & 206   & 40    & 2     & 45.1 \\
		40    & 12    & 205   & 40    & 2     & 54.1 \\
		40    & 13    & 203   & 40    & 2     & 47.9 \\
		40    & 14    & 225   & 40    & 2     & 52.3 \\
		40    & 15    & 183   & 40    & 2     & 46.0 \\
		40    & 16    & 198   & 40    & 2     & 50.1 \\
		40    & 17    & 211   & 40    & 2     & 52.4 \\
		40    & 18    & 178   & 40    & 2     & 58.8 \\
		40    & 19    & 384   & 40    & 1     & 51.6 \\
		40    & 20    & 208   & 40    & 2     & 49.7 \\
		40    & 21    & 190   & 40    & 2     & 48.2 \\
		40    & 22    & 215   & 40    & 2     & 56.2 \\
		40    & 23    & 197   & 40    & 2     & 49.9 \\
		40    & 24    & 196   & 40    & 2     & 47.9 \\
		40    & 25    & 204   & 40    & 2     & 53.4 \\
		40    & 26    & 221   & 40    & 2     & 52.5 \\
		40    & 27    & 161   & 40    & 2     & 51.6 \\
		40    & 28    & 211   & 40    & 2     & 58.8 \\
		40    & 29    & 183   & 40    & 2     & 53.6 \\
		50    & 0     & 310   & 50    & 2     & 62.8 \\
		50    & 1     & 295   & 50    & 2     & 57.5 \\
		50    & 2     & 265   & 50    & 3     & 50.8 \\
		50    & 3     & 324   & 50    & 2     & 50.0 \\
		50    & 4     & 269   & 50    & 3     & 57.4 \\
		50    & 5     & 253   & 50    & 3     & 57.6 \\
		50    & 6     & 273   & 50    & 2     & 64.1 \\
		50    & 7     & 292   & 50    & 2     & 58.2 \\
		50    & 8     & 335   & 50    & 2     & 50.4 \\
		50    & 9     & 285   & 50    & 1     & 59.3 \\
		50    & 10    & 297   & 50    & 2     & 57.7 \\
		50    & 11    & 281   & 50    & 2     & 52.1 \\
		50    & 12    & 284   & 50    & 2     & 56.1 \\
		50    & 13    & 296   & 50    & 2     & 56.4 \\
		50    & 14    & 309   & 50    & 2     & 45.5 \\
		50    & 15    & 291   & 50    & 2     & 47.0 \\
		50    & 16    & 293   & 50    & 2     & 53.2 \\
		50    & 17    & 273   & 50    & 2     & 65.3 \\
		50    & 18    & 318   & 50    & 2     & 45.1 \\
		50    & 19    & 329   & 50    & 2     & 54.5 \\
		50    & 20    & 548   & 50    & 1     & 54.3 \\
		50    & 21    & 295   & 50    & 2     & 51.9 \\
		50    & 22    & 318   & 50    & 2     & 61.7 \\
		50    & 23    & 316   & 50    & 2     & 51.2 \\
		50    & 24    & 307   & 50    & 2     & 51.0 \\
		50    & 25    & 265   & 50    & 3     & 50.9 \\
		50    & 26    & 259   & 50    & 2     & 51.3 \\
		50    & 27    & 324   & 50    & 2     & 50.8 \\
		50    & 28    & 263   & 50    & 2     & 62.5 \\
		50    & 29    & 280   & 50    & 2     & 47.7 \\
		60    & 0     & 434   & 60    & 3     & 58.6 \\
		60    & 1     & 398   & 60    & 3     & 49.6 \\
		60    & 2     & 439   & 60    & 2     & 59.3 \\
		60    & 3     & 335   & 60    & 3     & 56.5 \\
		60    & 4     & 450   & 60    & 3     & 45.8 \\
		60    & 5     & 394   & 60    & 3     & 58.2 \\
		60    & 6     & 307   & 60    & 3     & 56.4 \\
		60    & 7     & 396   & 60    & 3     & 58.6 \\
		60    & 8     & 394   & 60    & 2     & 48.9 \\
		60    & 9     & 433   & 60    & 2     & 62.4 \\
		60    & 10    & 361   & 60    & 3     & 60.3 \\
		60    & 11    & 396   & 60    & 2     & 56.9 \\
		60    & 12    & 349   & 60    & 2     & 64.2 \\
		60    & 13    & 363   & 60    & 3     & 55.6 \\
		60    & 14    & 355   & 60    & 2     & 50.4 \\
		60    & 15    & 336   & 60    & 3     & 62.3 \\
		60    & 16    & 327   & 60    & 3     & 65.0 \\
		60    & 17    & 412   & 60    & 2     & 60.4 \\
		60    & 18    & 387   & 60    & 3     & 51.0 \\
		60    & 19    & 405   & 60    & 2     & 58.9 \\
		60    & 20    & 355   & 60    & 3     & 52.9 \\
		60    & 21    & 439   & 60    & 2     & 57.2 \\
		60    & 22    & 414   & 60    & 2     & 49.6 \\
		60    & 23    & 420   & 60    & 2     & 60.2 \\
		60    & 24    & 431   & 60    & 2     & 51.3 \\
		60    & 25    & 326   & 60    & 3     & 60.9 \\
		60    & 26    & 378   & 60    & 2     & 58.4 \\
		60    & 27    & 387   & 60    & 3     & 48.0 \\
		60    & 28    & 336   & 60    & 3     & 59.7 \\
		60    & 29    & 344   & 60    & 3     & 59.3 \\
		70    & 0     & 461   & 70    & 3     & 65.1 \\
		70    & 1     & 523   & 70    & 2     & 67.3 \\
		70    & 2     & 507   & 70    & 3     & 59.4 \\
		70    & 3     & 556   & 70    & 3     & 61.4 \\
		70    & 4     & 442   & 70    & 3     & 71.7 \\
		70    & 5     & 478   & 70    & 3     & 55.5 \\
		70    & 6     & 500   & 70    & 4     & 56.3 \\
		70    & 7     & 573   & 70    & 3     & 56.0 \\
		70    & 8     & 435   & 70    & 3     & 58.5 \\
		70    & 9     & 453   & 70    & 3     & 60.6 \\
		70    & 10    & 459   & 70    & 3     & 59.7 \\
		70    & 11    & 455   & 70    & 3     & 63.3 \\
		70    & 12    & 414   & 70    & 3     & 57.5 \\
		70    & 13    & 473   & 70    & 3     & 61.5 \\
		70    & 14    & 545   & 70    & 2     & 61.3 \\
		70    & 15    & 425   & 70    & 3     & 61.4 \\
		70    & 16    & 580   & 70    & 2     & 65.2 \\
		70    & 17    & 434   & 70    & 3     & 60.3 \\
		70    & 18    & 483   & 70    & 2     & 63.8 \\
		70    & 19    & 503   & 70    & 3     & 61.1 \\
		70    & 20    & 403   & 70    & 3     & 59.3 \\
		70    & 21    & 462   & 70    & 3     & 59.2 \\
		70    & 22    & 488   & 70    & 3     & 60.3 \\
		70    & 23    & 478   & 70    & 3     & 58.2 \\
		70    & 24    & 456   & 70    & 2     & 54.2 \\
		70    & 25    & 518   & 70    & 2     & 64.2 \\
		70    & 26    & 394   & 70    & 4     & 54.1 \\
		70    & 27    & 504   & 70    & 3     & 58.0 \\
		70    & 28    & 385   & 70    & 3     & 65.0 \\
		70    & 29    & 536   & 70    & 2     & 54.9 \\
		80    & 0     & 594   & 80    & 3     & 63.3 \\
		80    & 1     & 506   & 80    & 3     & 71.9 \\
		80    & 2     & 582   & 80    & 3     & 62.5 \\
		80    & 3     & 664   & 80    & 3     & 63.8 \\
		80    & 4     & 539   & 80    & 3     & 62.3 \\
		80    & 5     & 513   & 80    & 4     & 68.2 \\
		80    & 6     & 643   & 80    & 3     & 65.6 \\
		80    & 7     & 555   & 80    & 5     & 56.4 \\
		80    & 8     & 564   & 80    & 3     & 59.9 \\
		80    & 9     & 671   & 80    & 3     & 58.7 \\
		80    & 10    & 628   & 80    & 3     & 57.9 \\
		80    & 11    & 602   & 80    & 3     & 66.3 \\
		80    & 12    & 731   & 80    & 3     & 58.2 \\
		80    & 13    & 555   & 80    & 3     & 66.1 \\
		80    & 14    & 644   & 80    & 2     & 66.5 \\
		80    & 15    & 568   & 80    & 3     & 62.4 \\
		80    & 16    & 579   & 80    & 3     & 56.3 \\
		80    & 17    & 601   & 80    & 3     & 55.2 \\
		80    & 18    & 515   & 80    & 3     & 67.9 \\
		80    & 19    & 603   & 80    & 3     & 64.9 \\
		80    & 20    & 700   & 80    & 3     & 64.3 \\
		80    & 21    & 577   & 80    & 4     & 61.3 \\
		80    & 22    & 598   & 80    & 4     & 60.8 \\
		80    & 23    & 541   & 80    & 3     & 65.8 \\
		80    & 24    & 689   & 80    & 3     & 62.2 \\
		80    & 25    & 665   & 80    & 3     & 54.8 \\
		80    & 26    & 600   & 80    & 3     & 62.7 \\
		80    & 27    & 548   & 80    & 3     & 57.8 \\
		80    & 28    & 570   & 80    & 4     & 59.1 \\
		80    & 29    & 509   & 80    & 5     & 62.2 \\
		90    & 0     & 829   & 90    & 3     & 62.3 \\
		90    & 1     & 743   & 90    & 3     & 63.4 \\
		90    & 2     & 816   & 90    & 3     & 73.5 \\
		90    & 3     & 815   & 90    & 3     & 58.5 \\
		90    & 4     & 699   & 90    & 5     & 60.4 \\
		90    & 5     & 842   & 90    & 3     & 58.3 \\
		90    & 6     & 806   & 90    & 4     & 63.5 \\
		90    & 7     & 694   & 90    & 3     & 66.0 \\
		90    & 8     & 837   & 90    & 3     & 66.6 \\
		90    & 9     & 584   & 90    & 4     & 62.5 \\
		90    & 10    & 765   & 89    & 5     & 62.2 \\
		90    & 11    & 727   & 90    & 3     & 60.5 \\
		90    & 12    & 710   & 90    & 3     & 65.1 \\
		90    & 13    & 866   & 90    & 3     & 65.0 \\
		90    & 14    & 803   & 90    & 5     & 61.2 \\
		90    & 15    & 726   & 90    & 3     & 63.5 \\
		90    & 16    & 744   & 90    & 3     & 61.4 \\
		90    & 17    & 692   & 90    & 4     & 68.4 \\
		90    & 18    & 667   & 90    & 4     & 62.2 \\
		90    & 19    & 868   & 90    & 3     & 64.2 \\
		90    & 20    & 674   & 90    & 4     & 66.3 \\
		90    & 21    & 645   & 90    & 5     & 64.0 \\
		90    & 22    & 767   & 90    & 3     & 62.6 \\
		90    & 23    & 693   & 90    & 3     & 62.6 \\
		90    & 24    & 657   & 90    & 4     & 61.8 \\
		90    & 25    & 779   & 90    & 3     & 57.2 \\
		90    & 26    & 840   & 90    & 4     & 57.8 \\
		90    & 27    & 734   & 90    & 4     & 64.9 \\
		90    & 28    & 653   & 90    & 4     & 65.9 \\
		90    & 29    & 721   & 90    & 5     & 58.6 \\
		100   & 0     & 1039  & 100   & 3     & 60.5 \\
		100   & 1     & 883   & 100   & 5     & 61.7 \\
		100   & 2     & 817   & 100   & 5     & 63.2 \\
		100   & 3     & 806   & 98    & 5     & 63.9 \\
		100   & 4     & 847   & 100   & 4     & 66.3 \\
		100   & 5     & 1041  & 100   & 3     & 62.9 \\
		100   & 6     & 740   & 100   & 6     & 68.4 \\
		100   & 7     & 909   & 100   & 3     & 59.1 \\
		100   & 8     & 956   & 100   & 4     & 64.3 \\
		100   & 9     & 816   & 100   & 4     & 70.3 \\
		100   & 10    & 981   & 100   & 4     & 67.5 \\
		100   & 11    & 1049  & 99    & 5     & 58.0 \\
		100   & 12    & 787   & 95    & 6     & 64.9 \\
		100   & 13    & 997   & 100   & 3     & 65.1 \\
		100   & 14    & 947   & 99    & 5     & 59.8 \\
		100   & 15    & 868   & 100   & 5     & 61.9 \\
		100   & 16    & 924   & 97    & 5     & 62.2 \\
		100   & 17    & 813   & 100   & 3     & 68.4 \\
		100   & 18    & 1018  & 100   & 3     & 62.8 \\
		100   & 19    & 947   & 100   & 4     & 62.9 \\
		100   & 20    & 867   & 100   & 4     & 64.6 \\
		100   & 21    & 864   & 100   & 6     & 62.6 \\
		100   & 22    & 829   & 100   & 5     & 66.4 \\
		100   & 23    & 808   & 100   & 3     & 72.8 \\
		100   & 24    & 832   & 100   & 5     & 64.8 \\
		100   & 25    & 780   & 98    & 5     & 63.9 \\
		100   & 26    & 862   & 99    & 6     & 60.4 \\
		100   & 27    & 914   & 100   & 4     & 68.2 \\
		100   & 28    & 828   & 100   & 5     & 64.4 \\
		100   & 29    & 883   & 92    & 6     & 57.0 \\
\end{supertabular}{\makeatletter\def\@currentlabel{\thetable}\label{tab:detailed_1ULD}}
\setlength{\tabcolsep}{6pt}
\end{multicols*}
	\addtocounter{table}{+1}
\begin{multicols*}{4}
	[\begin{center}Table \thetable: Detailed \ULD\ results for adapted \cite{PaquayEtAl2018} instances where each \ULD\ type is available only once.\end{center}]

\setlength{\tabcolsep}{1pt}
\tiny

\centering

\tablefirsthead{%
	\toprule
	\multicolumn{1}{l}{$\numItems$} & \multicolumn{1}{l}{\numInstance} & \containerID & $\sortingCriterion$ & \multicolumn{1}{l}{$\utilization$} & \multicolumn{1}{l}{\CoGDevX} & \multicolumn{1}{l}{\CoGDevY} \\
	\midrule}
\tablehead{%
	\toprule
	\multicolumn{1}{l}{$\numItems$} & \multicolumn{1}{l}{\numInstance} & \containerID & $\sortingCriterion$ & \multicolumn{1}{l}{$\utilization$} & \multicolumn{1}{l}{\CoGDevX} & \multicolumn{1}{l}{\CoGDevY} \\
	\midrule}
\tabletail{%
	\midrule
	\multicolumn{7}{r}{\textit{Continued on next column}}\\}
\tablelasttail{\bottomrule}

\TrickSupertabularIntoMulticols

\begin{supertabular}{rrllrrr}
		10    & 0     & 5     & \stackCumVol & 24.9  & 3.0   & 2.6 \\
		10    & 1     & 4     & \stackCumVol & 27.9  & 6.5   & 6.6 \\
		10    & 2     & 4     & \stackCumVol & 41.6  & 13.9  & 8.9 \\
		10    & 3     & 4     & \stackCumVol & 27.4  & 1.8   & 6.8 \\
		10    & 4     & 3     & \stackCumVol & 56.1  & 0.2   & 2.9 \\
		10    & 5     & 5     & \stackCumVol & 13.1  & 1.8   & 0.0 \\
		10    & 6     & 1     & \cumVol & 48.5  & 9.3   & 24.7 \\
		10    & 7     & 3     & \stackCumVol & 58.5  & 1.5   & 6.6 \\
\end{supertabular}{\makeatletter\def\@currentlabel{\thetable}\label{tab:detailed_con_1ULD}}
\setlength{\tabcolsep}{6pt}
\end{multicols*}
	
\end{appendix}

\end{document}